\documentclass{article}
\usepackage{mathrsfs}
\usepackage{amsmath}
\usepackage{amssymb, amsmath}
\usepackage{graphicx}
\usepackage{fancyhdr}
\usepackage [latin1]{inputenc}
\catcode`\@=11 \@addtoreset{equation}{section}

\catcode`\@=12
\usepackage{color}

\usepackage[colorlinks,
            linkcolor=blue,
            anchorcolor=blue,
            citecolor=red,
            ]{hyperref}

\DeclareFontFamily{U}{mathx}{}
\DeclareFontShape{U}{mathx}{m}{n}{ <-> mathx10 }{}
\DeclareSymbolFont{mathx}{U}{mathx}{m}{n}
\DeclareFontSubstitution{U}{mathx}{m}{n}

\DeclareMathAccent{\widecheck}{0}{mathx}{"71}

\makeatletter
\newcommand{\wcheck}[1]{\mathpalette\wcheck@{#1}}
\newcommand{\wcheck@}[2]{%
  \begingroup
  \edef\wcheck@font{\the
    \ifx#1\displaystyle\textfont\else\ifx#1\textstyle\textfont
    \else\ifx#1\scriptstyle\scriptfont\else\scriptscriptfont\fi\fi\fi\@ne
  }%
  \sbox\z@{\wcheck@font\mbox{#2}\mbox{\char\the\skewchar\font}}%
  \sbox\tw@{\wcheck@font#2\char\the\skewchar\font}%
  \dimen@=\dimexpr\wd\tw@-\wd\z@\relax
  {\,\kern2\dimen@\widecheck{\!\kern-2\dimen@#2\!}\,}%
  \endgroup
}
\makeatother

\newtheorem{Theorem}{Theorem}[section]
\newtheorem{Proposition}{Proposition}[section]
\newtheorem{Lemma}{Lemma}[section]
\newtheorem{Corollary}{Corollary}[section]

\newtheorem{Remark}{Remark}[section]

\newcommand{\newcom}{\newcommand}
\newcommand{\bTheorem}[1]{
\begin{Theorem} \label{T#1} }
\newcommand{\eT}{\end{Theorem}}

\newcommand{\bProposition}[1]{
\begin{Proposition} \label{P#1}}
\newcommand{\eP}{\end{Proposition}}

\newcommand{\bLemma}[1]{
\begin{Lemma} \label{L#1} }
\newcommand{\eL}{\end{Lemma}}

\newcommand{\bCorollary}[1]{
\begin{Corollary} \label{C#1} }
\newcommand{\eC}{\end{Corollary}}

\newcommand{\beq}{\begin{equation}}
\newcommand{\eeq}{\end{equation}}
\newcom{\ben}{\begin{eqnarray}}
\newcom{\een}{\end{eqnarray}}
\newcom{\beno}{\begin{eqnarray*}}
\newcom{\eeno}{\end{eqnarray*}}
\newcom{\bali}{\begin{aligned}}
\newcom{\eali}{\end{aligned}}

\newcommand{\bFormula}[1]{
\begin{equation} \label{#1}}
\newcommand{\eF}{\end{equation}}

\newcommand{\De}{\Delta}
\newcommand{\f}{\frac}

\newcommand{\p}{\partial}

\newcommand{\vr}{\varrho}
\newcommand{\vae}{a^\ep}
\newcommand{\Td}{\mathbb{T}_b^d}
\newcommand{\T}{|\mathbb{T}_b^d|}

\newcommand{\vre}{\vr^\ep}

\newcommand{\vue}{\vu^\ep}

\newcommand{\vVe}{\vV^\ep}
\newcommand{\vWe}{\vW^\ep}
\newcommand{\vZe}{\vZ^\ep}

\newcommand{\va}{\vc{a}}
\newcommand{\vu}{\vc{u}}

\newcommand{\vZ}{\vc{Z}}

\newcommand{\vv}{\mathbf{v}}
\newcommand{\vV}{\vc{V}}

\newcommand{\vW}{\vc{W}}

\newcommand{\vc}[1]{{\boldsymbol #1}}

\newcommand{\Div}{{\rm div}}
\newcommand{\Grad}{\nabla}

\newcommand{\dx}{\,{\rm d} x}

\newcommand{\dt}{\,{\rm d} t }

\newcommand{\intTd}[1]{\int_{\mathbb{T}_b^d} #1 \ \ \dx}

\newcommand{\bProof}{{\bf Proof: }}

\newcommand{\ep}{\varepsilon}

\newcommand\Cbox[2]{%
    \newbox\contentbox%
    \newbox\bkgdbox%
    \setbox\contentbox\hbox to \hsize{%
        \vtop{
            \kern\columnsep
            \hbox to \hsize{%
                \kern\columnsep%
                \advance\hsize by -2\columnsep%
                \setlength{\textwidth}{\hsize}%
                \vbox{
                    \parskip=\baselineskip
                    \parindent=0bp
                    #2
                }%
                \kern\columnsep%
            }%
            \kern\columnsep%
        }%
    }%
    \setbox\bkgdbox\vbox{
        \color{#1}
        \hrule width  \wd\contentbox %
               height \ht\contentbox %
               depth  \dp\contentbox
        \color{black}
    }%
    \wd\bkgdbox=0bp%
    \vbox{\hbox to \hsize{\box\bkgdbox\box\contentbox}}%
    \vskip\baselineskip%
}

\begin{document}


\pagestyle{fancy} \lhead{\color{blue}{Low Mach number limits on the tours}} \rhead{\emph{S.Li}}

\title{\bf Low Mach number limit of  the compressible Navier-Stokes system for  large initial date with critical regularity on the torus} \rhead{\emph{S.Li}}

\author{
Sai Li    \\
Department of Applied Mathematics, \\ Nanjing Forestry University, Nanjing, 210037, People's Republic of China \\ Email: lsmath@njfu.edu.cn
}

\maketitle

{\centerline {\bf Abstract }}
\vspace{2mm}
{We study the low Mach number limit of the compressible Navier-Stokes equations on the torus. For large initial data with critical regularity, we prove that solutions to the compressible Navier-Stokes system exist as long as the corresponding solutions to the incompressible Navier-Stokes system exist, provided that the Mach number is sufficiently small. Furthermore, we establish the convergence of solutions of the compressible system to those of the incompressible system as the Mach number tends to zero. Our approach combines high-medium-low frequency analysis of density and velocity with the solution filtering technique via acoustic wave groups. This work provides an affirmative answer to the  problem posed by Danchin [\emph{Amer. J. Math.}, 124 (2002), 1153-1219]:\textquotedblleft Does convergence hold for large data with critical regularity ?\textquotedblright
}

\vspace{2mm}


\vspace{2mm}


\tableofcontents

\section{Introduction}
We focus on the mathematical theory of low Mach number flows and investigate following compressible Navier-Stokes system on the torus $\mathbb{T}_b^d  (d\geq2)$:
\begin{equation}\label{equ1}
  \left\{\bali
  &\p_t\vre+\Div(\vre\vue)=0,\\
  &\p_t(\vre\vue)+\Div(\vre\vue\otimes\vue)-\mu\De\vue-(\lambda+\mu)\Grad\Div\vue+\f{\Grad P(\vre)}{\ep^2}=\vre f^\ep,\\
  &(\vre,\vue)|_{t=0}=(\vr^\ep_0,\vu^\ep_0).
  \eali
  \right.
\end{equation}

Here, the unknowns $\vre(t,x)$ and $\vue(t,x)$ represent the density of fluids and the velocity field where $t\geq0$ is the time coordinate and the $x\in\mathbb{R}^d$ is spatial coordinate. $\vre, \vue$ are assumed to be $2\pi b_h$-periodic ($b_h>0$) in variable $x_h$ ($1\leq h\leq d$). The viscosity coefficients $\mu$ and $\lambda$ satisfy $\mu>0,\nu:=2\mu+\lambda>0$. $0<\ep\leq1$ is the Mach number. $P(\cdot)$ is a suitable smooth function which represents the scalar pressure and $f^\ep$ is the external force.

A fundamental question is:\textquotedblleft How does the solution of system \eqref{equ1}  depend on the Mach number $\ep$ ?\textquotedblright  Over the years, rigorous mathematical analysis has demonstrated that under appropriate conditions, as the parameter $\ep$ tends to zero, the solutions of \eqref{equ1}  converge to those of the incompressible Navier-Stokes equations.

The rigorous theoretical study of the low Mach number limit was initiated by Klainerman and Majda \cite{KM} (for the inviscid case, see Ebin \cite{EDG} and Klainerman and Majda \cite{KM1}), who first established the convergence theory for local strong solutions of system \eqref{equ1}. Under the assumption that the limiting incompressible Navier-Stokes equations admit global solutions, Hagstrom and Lorenz \cite{HLo} proved the existence of global strong solutions to system \eqref{equ1} for certain initial data and rigorously justified the incompressible limit (Hoff \cite{DH} also established analogous global convergence results). Lions and Masmoudi \cite{LM} were the first to investigate the low Mach number limit within the framework of global weak solutions; Desjardins and Grenier \cite{DG} subsequently simplified the proof in \cite{LM} by employing dispersion estimates in the whole-space setting. Furthermore, Desjardins et al. \cite{DGLM} analyzed the case of bounded domains with Dirichlet boundary conditions, while Feireisl et al. \cite{Fe3}  further explored the scenario where the domain depends on the Mach number $\ep$. Due to the scope of this article, only a concise review is provided here. For a more systematic and in-depth exposition, readers are referred to the monograph by Feireisl and  Novotn\'{y} \cite{Fe4}, the review article by Jiang and Masmoudi \cite{JN}, and the survey by Schochet \cite{SS1}.

This paper investigates the convergence of system \eqref{equ1} to the incompressible Navier-Stokes system in critical spaces \footnote{It is clear that \eqref{equ1} is invariant by the transform
\begin{align}\label{trans}
(\vre(t,x),\ \ \vue(t,x))\rightarrow(\vre(l^2t,lx),\ \ l\vue(l^2t,lx)),\ \ l>0,
\end{align}
provided the pressure term $P$ has been changed to $l^2P$. Critical spaces for system \eqref{equ1} are norm invariant for all
$l$ (up to an irrelevant constant) with respect to the transform \eqref{trans}.}, where the well-posedness in such spaces was established by Danchin \cite{D4,D5}.
Without loss of generality, we normalize the reference density to $1$ and assume the pressure law satisfies $P'(1)=1$.
We set $\vr^\ep_0=1+\ep a^\ep_0$ and introduce
$\vae:=\f{\vre-1}\ep$, then $(\vae,\vue)$ satisfies
\begin{equation}\label{equ2}
  \left\{\bali
  &\p_t\vae+\f{\Div\vue}\ep=-\Div(\vae\vue),\\
  &\p_t\vue+\vue\cdot\Grad\vue-\mathcal{A}\vue+\f{\Grad\vae}\ep=-\kappa\vae\Grad\vae-\widetilde{K}(\ep\vae)\vae\Grad\vae-I(\ep\vae)\mathcal{A}\vue+f^\ep,\\
  &(\vae,\vue)\mid_{t=0}=(a^\ep_0,\vue_0),
  \eali
  \right.
\end{equation}
where
\[
\mathcal{A}\vu:=\mu\Delta {\vu} +{(\mu+\lambda) }\nabla \Div{\vu},\ \ I(a):=\f a {1+a}, \ \ \f{P'(1+a)}{1+a}=1+\kappa a+a\widetilde{K}(a) \text{ with } \widetilde{K}(0)=0.
\]
For simplicity, we assume hereafter that the initial data $(a^\ep_0,\vue_0)$ and external force $f^\ep$  are independent of $\ep$, and we
denote them by $(a_0,\vu_0)$ and $f$, respectively. The \emph{critical regularity} condition on the initial data is given by
\begin{align}\label{cond1}
(a_0,\vu_0)\in \underline{B}_{2,1}^{\f{d} 2}(\mathbb{T}_b^d)\times\left(B_{2,1}^{\f{d} 2-1}(\mathbb{T}_b^d)\right)^d,
\end{align}
where $\underline{B}_{2,1}^{\f{d} 2}(\mathbb{T}_b^d)=\{g\in B_{2,1}^{\f{d} 2}(\mathbb{T}^d):\intTd g=0\}$. 

Let $\mathbb{P}$ and $\mathbb{Q}$ denote the  projection operators onto divergence-free and gradient vector fields respectively, defined by
\[
\mathbb{Q}:=-\Grad(-\Delta)^{-1}\Div,\ \ \mathbb{P}:=I-\mathbb{Q}.
\]
One may guess that $(\vae,\vue)$ converges to $(0,\vv)$, where $\vv$ satisfies following incompressible Navier-Stokes system
\begin{equation}\label{equ3}
  \left\{\bali
  &\p_t\vv-\mu\Delta\vv+\vv\cdot\Grad\vv+\Grad \pi=h,\ \ \Div \vv=0\\
  &\vv\mid_{t=0}=\vv_0,
  \eali
  \right.
\end{equation}
with  the initial data $\vv_0=\mathbb{P}\vu_0$ and the external force $h=\mathbb{P}f$.

Danchin \cite{D1, D3} established that for both the torus and the whole space cases, whenever the solution $\vv$ exists, the solution $(\vae,\vue)$ also exists and converges to $(0,\vv)$, provided that the initial data are either small in the sense that
\begin{align} \label{critical}
  \ep\|a_0\|_{\dot{B}^{\f d 2}_{2,1}}+\|a_0\|_{\dot{B}^{\f d 2-1}_{2,1}}+\|\vu_0\|_{\dot{B}^{\f d 2-1}_{2,1}}\ll1,
 \end{align}
 or large with subcritical regularity.  Subsequently, Danchin and He  \cite{D2}  extended the results of  \cite{D3}  to critical $L^p$ spaces, and Li and Mu \cite{Lif}  generalized the theorems of  \cite{D1} to the magnetohydrodynamic equations.
However, the convergence remains open for large data satisfying \eqref{cond1}  on the torus (see \cite{D1}, p. 1162), and for large initial data satisfying
$(a_0,\vu_0)\in \left(B_{2,1}^{\f{d} 2}(\mathbb{R}^d)\cap B_{2,1}^{\f{d} 2-1}(\mathbb{R}^d)\right)\times\left(B_{2,1}^{\f{d} 2-1}(\mathbb{R}^d)\right)^d$
on the whole space (see \cite{D3}, Remark 0.5) . In a recent advancement, Fujii \cite{F} resolved this problem for the whole-space case through the development of an innovative high-medium-low frequency framework integrated with dispersive estimate techniques. This methodology was subsequently extended in Fujii and Li \cite{F1} to examine the low Mach number limit of global large solutions to the Navier-Stokes-Korteweg system on $\mathbb{R}^2$ in critical $\widehat{L^p}$ spaces. Furthermore, Li \cite{LS} extended the results of \cite{F} to the magnetohydrodynamic equations on $\mathbb{R}^d$, and Ogino \cite{O} established the diverging  of  $\p_t\mathbb{Q}\vue$ in system  \eqref{equ1}  as the Mach number tends to zero in the whole-space setting.

The torus configuration demonstrates qualitatively distinct dynamical behavior compared to the whole space case, particularly in terms of acoustic wave propagation characteristics. Upon application of the projection operator $\mathbb{Q}$
subsystem \eqref{equ2}$_2$, we obtain the following coupled system governing acoustic wave dynamics:
\begin{equation}
  \left\{\bali
  &\p_t\vae+\f{\Div\mathbb{Q}\vue}\ep=-F,\\
  &\p_t\mathbb{Q}\vue+\f{\Grad\vae}\ep=\mathbb{Q}G\\
  \eali
  \right.\nonumber
\end{equation}
 where
 \[
 F=\Div(\vae\vue),\ \ G=\mathcal{A}\vue-\vue\cdot\Grad\vue-\kappa\vae\Grad\vae-\widetilde{K}(\ep\vae)\vae\Grad\vae-I(\ep\vae)\mathcal{A}\vue+f.
 \]
This system generates highly oscillatory acoustic waves. In the whole space setting, these waves exhibit natural dispersion phenomena that facilitate asymptotic energy dissipation through spatial spreading. However, on the torus, the boundedness of the domain fundamentally suppresses dispersion mechanisms, resulting in persistent acoustic oscillations that consequently yields only weak convergence of the solutions. This fundamental difference in wave propagation dynamics renders the approach of Fujii \cite{F}, which crucially relies on dispersive estimates, inapplicable to the torus case.

To investigate these persistent acoustic oscillations, we adopt the filtering technique, which has been successfully implemented in various contexts: Schochet \cite{SS2} for hyperbolic system;  Gallagher \cite{G1} for  parabolic system;  Grenier \cite{Gr}, Danchin \cite{D1}, and Masmoudi \cite{M} for Navier-Stokes system, among other references. We define the filtering operator $\mathcal{L}(\tau)=e^{-\tau L}$ where $L$ acts on $\underline{\mathcal{S}}'(\mathbb{T}_b^d)\times \Big(\mathcal{S}'(\mathbb{T}_b^d)\Big)^d$ as
\[
L\left(
\begin{array}{l}
a\\
\vu
\end{array}
\right)
=\left(
\begin{array}{l}
\Div\vu\\
\Grad a
\end{array}
\right).
\]
Letting  $\mathcal{L}^1(\tau)$ and $\mathcal{L}^2(\tau)$  denote the first and remaining
$d$ components of $\mathcal{L}(\tau)$ respectively. We set $\vVe:=\mathcal{L}\Big(-\f t {\ep}\Big)\left(
\begin{array}{l}
\ \ \vae\\
\mathbb{Q}\vue
\end{array}
\right)$. Applying the operator $\mathcal{L}(-\f t {\ep})$ to system \eqref{equ2} yields the following filtered system,
\begin{align}\label{equ5}
&\p_t\vVe+\mathcal{Q}_1^\ep(\mathbb{P}\vue,\vVe)+\mathcal{Q}_2^\ep(\vVe,\vVe)-\nu\mathcal{A}_2^\ep(D)\vVe \nonumber\\
&=\mathcal{L}\Big(-\f t {\ep}\Big)\left(\begin{array}{l}
\qquad\qquad\qquad\qquad 0 \\
\mathbb{Q}(f-\mathbb{P}\vue\cdot\Grad\mathbb{P}\vue-\widetilde{K}(\ep\vae)\vae\Grad\vae-I(\ep\vae)\mathcal{A}\vue)
\end{array}\right),
\end{align}
where
\[
\mathcal{A}_2^\ep(D)B:=\mathcal{L}\Big(-\f t {\ep}\Big)\left(
\begin{array}{l}
\qquad  0\\
\De(\mathcal{L}^2\Big(\f t {\ep}\Big)B)
\end{array}
\right),
\]
\[
\mathcal{Q}_1^\ep(\vu,B):=\mathcal{L}\Big(-\f t {\ep}\Big)\left(
\begin{array}{l}
\qquad\qquad\Div\left(\vu\mathcal{L}^1\Big(\f t {\ep}\Big)B\right)\\
\mathbb{Q}\left(\vu\cdot\Grad\mathcal{L}^2\Big(\f t {\ep}\Big)B+\mathcal{L}^2\Big(\f t {\ep}\Big)B\cdot\Grad\vu\right)
\end{array}
\right),
\]
\[
\mathcal{Q}_2^\ep(A,B):=\f1 2\mathcal{L}\Big(-\f t {\ep}\Big)\left(
\begin{array}{l}
\qquad\Div\left(\mathcal{L}^1\Big(\f t {\ep}\Big)A\mathcal{L}^2\Big(\f t {\ep}\Big)B+\mathcal{L}^2\Big(\f t {\ep}\Big)A\mathcal{L}^1\Big(\f t {\ep}\Big)B\right)\\
\Grad\left(\mathcal{L}^2\Big(\f t {\ep}\Big)A\cdot\mathcal{L}^2\Big(\f t {\ep}\Big)B\right)+\kappa\Grad\left(\mathcal{L}^1\Big(\f t {\ep}\Big)A\cdot\mathcal{L}^1\Big(\f t {\ep}\Big)B\right)
\end{array}
\right).
\]
Formally, $\vVe$ converges to the solution $\vV$ of the limit system
\begin{equation}\label{equ4}
  \left\{\bali
  &\p_t\vV+\mathcal{Q}_1(\vv,\vV)+\mathcal{Q}_2(\vV,\vV)-\f\nu2\De\vV=0,\\
  &\vV|_{t=0}=\vV_0.
  \eali
  \right.
\end{equation}
with the initial data $\vV_0=
\left(
\begin{array}{l}
\ \ a_0\\
\mathbb{Q}\vu_0
\end{array}
\right)$.  The complete derivation of this limiting system appears in Section \ref{Sec3}.

 \emph{Our primary objective is to resolve Danchin's open problem for the torus case.} We expect to prove that for arbitrarily large initial data satisfying \eqref{cond1}, $(\vae, \vue)$ exists whenever $\vv$ exists and $(\mathbb{P}\vue,\vVe)$ strongly converges to $(\vv,\vV)$. The strong convergence of $\vVe$ to $\vV$, combined with the representation$\left(
\begin{array}{l}
\ \ \vae\\
\mathbb{Q}\vue
\end{array}
\right)=\mathcal{L}\Big(\f t {\ep}\Big)(\vVe-\vV)+\mathcal{L}\Big(\f t {\ep}\Big)\vV$ implies the weak convergence of $(\vae,\mathbb{Q}\vue)$ to zero.

\medskip

The remainder of this paper is organized as follows. Section \ref{Sec2} introduces the basic notation and states the main theorem. In Section \ref{newsec}, we outline the proof strategy, discuss the difficulties encountered, and present the methods used to address them. Section \ref{complete} provides a complete derivation of Theorem \ref{main}. Section \ref{Sec3} is devoted to analyzing the decay properties of $\vVe - \vV$ and $\mathbb{P}\vue - \vv$. In Section \ref{Sec4}, we establish control relations among several energy functionals and prove key estimates for the high and medium frequencies. Finally, Section \ref{Sec5} presents a collection of fundamental lemmas essential to our analysis.

Throughout this paper,  the letter $C$ denotes a generic positive constant that is independent of $\ep$ and whose value may vary from line to line. We write $C(a_1,a_2,\cdots,a_n)$ to
indicate a positive constant that depends only on $a_1,a_2,\cdots,a_n$.  The notation $a_1\lesssim  a_2$  means that $a_1\leq Ca_2$, and  $a_1\approx a_2$ indicates that $C^{-1}a_1\leq a_2\leq Ca_1$. For two Banach spaces $B_1$ and $B_2$, the norm on their intersection is defined is defined by $\|\cdot\|_{B_1\cap B_2}=\|\cdot\|_{B_1}+\|\cdot\|_{B_2}$.

\section{Notations and the main theorem} \label{Sec2}
\subsection{Notations}
Let $\varphi$ be a symmetric, non-negative smooth bump function supported in $\{\xi\in\mathbb{R}:\f3 4\leq|\xi|\leq \f8 3\}$  satisfying the dyadic partition of unity condition:
\[
\sum_{j\in\mathbb{Z}}\varphi(2^{-j}\xi)=1\text{ for any } \xi\in\mathbb{R}\backslash\{0\}.
\]
(see \cite{BCD}  for construction of such functions).

For $g\in\mathcal{S}'(\mathbb{T}_b^d)$, $k \in \widetilde{\mathbb{Z}}^d$ and $j\in\mathbb{Z}$, we define:
\[
\hat{g}_k := \left< g,\f{e^{-ik\cdot x}} {\sqrt{|\mathbb{T}_b^d|}} \right>_{\mathcal{S}'(\mathbb{T}_b^d) \times\mathcal{S}(\mathbb{T}_b^d)},
\]
\[
\Delta_j g  := \sum\limits_{k \in \widetilde{\mathbb{Z}}^d}\varphi(2^{-j}|k|)\hat{g}_k\f{e^{ik\cdot x}} {\sqrt{|\mathbb{T}_b^d|}},\ \  S_j g := \f{\hat{g}_0}{\sqrt{|\mathbb{T}_b^d|}}+\sum\limits_{j'\leq j-1}\Delta_{j'} g,
\]
where  $\mathcal{S}(\mathbb{T}_b^d)$ and $\mathcal{S}'(\mathbb{T}_b^d)$ denote the space of smooth functions which are $2\pi b_h$-periodic ($1\leq h\leq d$) in the $h$th variable and its dual space, $\widetilde{\mathbb{Z}}^d:=\mathbb{Z}/b_1\times \mathbb{Z}/b_2\times\cdots\times\mathbb{Z}/b_d$ and $|\mathbb{T}_b^d|:=2\pi b_1\cdot 2\pi b_2\cdots 2\pi b_d$. This yields the decomposition:
\[
g=\sum\limits_{k \in \widetilde{\mathbb{Z}}^d} {{\hat g}_k}\,\f{e^{ ik \cdot x}}{\sqrt{\T}}=\f{\hat{g}_0}{\sqrt{|\mathbb{T}_b^d|}}+\sum\limits_{j\in \mathbb{Z}} \Delta_j g.
\]

We set $j_b:=\max\{j\in \mathbb{Z}:2^{-j}\cdot \min\{\f1 {b_1},\cdots,\f1 {b_d}\}\geq \f 8 3\}$, then
\[
\varphi(2^{-j}|k|)=0,\, \text{ if } j\leq j_b,\ \ k \in \widetilde{\mathbb{Z}}^d.
\]
It follows that
\beq\label{zero}
\Delta_j g\equiv0,\, \text{ if } j\leq j_b.
\eeq

 We define the zero-mean subspace $\underline{\mathcal{S}}'(\mathbb{T}_b^d):=\{g\in \mathcal{S}'(\mathbb{T}_b^d): \hat{g}_0=0\}$ with projection $\underline{g}:=g-\f{\hat{g}_0}{\sqrt{|\mathbb{T}_b^d|}}$.

 For $1\leq p,r\leq \infty $ and $s\in\mathbb{R}$, we introduce:
\begin{align*}
\intTd g &:= \int_{[0,2\pi b_1)\times [0,2\pi b_2)\times\cdots\times [0,2\pi b_d)}g\ \ \dx,\ \ \|g\|_{L^p(\mathbb{T}_b^d)}:=\left(\intTd {|g|^p}\right)^{\f1 p},\\
\|g\|_{B_{p,r}^s}  &
:=
\left(\sum\limits_{j\in\mathbb{Z}}2^{jsr}\|\Delta_{j} g\|^{r}_{L^p(\mathbb{T}_b^d)}+\left\|\f{\hat{g}_0}{\sqrt{\T}}\right\|^{r}_{L^p(\mathbb{T}_b^d)}\right)^{\f1 r},\\
 B_{p,r}^s(\mathbb{T}_b^d) & :=
 \left\{g\in\mathcal{S}'(\mathbb{T}_b^d):\|g\|_{B_{p,r}^s}<\infty \right\},\ \  \underline{B}_{p,r}^s(\mathbb{T}_b^d):=B_{p,r}^s(\mathbb{T}_b^d)\cap \underline{\mathcal{S}}'(\mathbb{T}_b^d),
 \end{align*}
 \begin{align*}
 H^s(\mathbb{T}_b^d) & :=
 \left\{g\in\mathcal{S}'(\mathbb{T}_b^d):\|g\|_{H^s}:\left(|{\hat{g}_0}|^2+\sum\limits_{k\neq0}|k|^{2s}|{\hat{g}_k}|^2\right)^{\f1 2}<\infty \right\}.
\end{align*}
It is clear that for $s\in\mathbb{R}$, $B_{2,2}^s(\Td)=H^s(\Td)$ and $\|g\|_{B_{2,2}^s}\approx \|g\|_{H^s}$.

For $1\leq q,p,r\leq\infty$, $s\in\mathbb{R}$, $0<T\leq\infty$ and a time-dependent functions $\Psi$ taking values in Besov spaces, we define:
\begin{align*}
&\|\Psi\|_{L^q_T(B_{p,r}^s)} :=\left\|\|\Psi(t,\cdot)\|_{B_{p,r}^s}\right\|_{L^q(0,T)},\\ &\|\Psi\|_{\widetilde{L}^{q}_T(B_{p,r}^s)}:=\left(\sum\limits_{j\in\mathbb{Z}}2^{jsr}\|\Delta_{j} \Psi\|^{r}_{L^q(0,T;L^p)}+\left\|\f{\hat{\Psi}_0}{\sqrt{\T}}\right\|^{r}_{L^q(0,T;L^p)}\right)^{\f1 r}.
\nonumber
\end{align*}
The second norms was first introduced by Chemin and Lerner in \cite{CL}. The corresponding zero-mean versions are defined as:
\[
\|\Psi\|_{L^q_T(\underline{B}_{p,r}^s)}:=\|\underline{\Psi}\|_{L^q_T(B_{p,r}^s)},\ \ \|\Psi\|_{L^q_T(\underline{H}^s)}:=\|\underline{\Psi}\|_{L^q_T(H^s)},\ \ \|\Psi\|_{\widetilde{L}^{q}_T(\underline{B}_{p,r}^s)}:=\|\underline{\Psi}\|_{\widetilde{L}^{q}_T(B_{p,r}^s)}.
\]

From Minkowski inequality, there holds
\begin{align}\label{minski}
\|\Psi\|_{L^q_T(B_{p,r}^s)}\leq\|\Psi\|_{\widetilde{L}^{q}_T(B_{p,r}^s)}, \ \ \text{if}\ \ r\leq q;\quad \|\Psi\|_{L^q_T(B_{p,r}^s)}\geq\|\Psi\|_{\widetilde{L}^{q}_T(B_{p,r}^s)}, \ \ \text{if}\ \ r\geq q.
\end{align}

Thanks to \eqref{zero}, we have that for $1\leq p,q,r_1,r_2\leq\infty$ and $-\infty<s_1< s_2<+\infty$,
\begin{align}
 &\max\{\|\Psi\|_{L^q_T(B_{p,r_1}^{s_1})},\|\Psi\|_{\widetilde{L}^{q}_T(B_{p,r_1}^{s_1})}\}
 \lesssim\min\{\|\Psi\|_{L^{q}_T(B_{p,r_2}^{s_2})},\|\Psi\|_{\widetilde{L}^{q}_T(B_{p,r_2}^{s_2})}\} \label{minski2}
\end{align}

 For $0\leq\zeta < \eta<\infty$,  we define
\[
\|g\|^{h;\eta}_{B_{p,r}^s}:=\left(\sum\limits_{2^j\geq\eta}2^{jsr}\|\Delta_{j} g\|^{r}_{L^p(\mathbb{T}_b^d)}\right)^{\f1 r},\ \ \|g\|^{m;\zeta,\eta}_{B_{p,r}^s}:=\left(\sum\limits_{\zeta\leq2^j<\eta}2^{jsr}\|\Delta_{j} g\|^{r}_{L^p(\mathbb{T}_b^d)}\right)^{\f1 r},
\]
\[
\|g\|^{l;\zeta}_{B_{p,r}^s}:=\left(\sum\limits_{2^j<\zeta}2^{jsr}\|\Delta_{j} g\|^{r}_{L^p(\mathbb{T}_b^d)}+\left\|\f{\hat{g}_0}{\sqrt{\T}}\right\|^{r}_{L^p(\mathbb{T}_b^d)}\right)^{\f1 r},
\]
\[
\|\Psi\|^{h;\eta}_{L^q_T(\dot{B}_{p,r}^s)}:=\left\|\|\Psi\|^{h;\eta}_{B_{p,r}^s}\right\|_{L^q(0,T)},
\]
\[ \|\Psi\|^{m;\zeta,\eta}_{L^q_T(B_{p,r}^s)}:=\left\|\|\Psi\|^{m;\zeta,\eta}_{B_{p,r}^s}\right\|_{L^q(0,T)}, \ \ \|\Psi\|^{l;\zeta}_{L^q_T(B_{p,r}^s)}:=\left\|\|\Psi\|^{l;\zeta}_{B_{p,r}^s}\right\|_{L^q(0,T)},
\]
\[
\|\Psi\|^{h;\eta}_{\widetilde{L}^q_T(\dot{B}_{p,r}^s)}:=\left(\sum\limits_{2^j\geq\eta}2^{jsr}\|\Delta_{j} \Psi\|^{r}_{L^q(0,T;L^p)}\right)^{\f1 r},
\]
\[
\|\Psi\|^{m;\zeta,\eta}_{\widetilde{L}^q_T(\dot{B}_{p,r}^s)}:=\left(\sum\limits_{\zeta\leq2^j<\eta}2^{jsr}\|\Delta_{j} \Psi\|^{r}_{L^q(0,T;L^p)}\right)^{\f1 r},
\]
\[
\|\Psi\|^{l;\zeta}_{\widetilde{L}^q_T(B_{p,r}^s)}
:=
\left(\sum\limits_{2^j<\zeta}2^{jsr}\|\Delta_{j} \Psi\|^{r}_{L^q(0,T;L^p)}+\left\|\f{\hat{\Psi}_0}{\sqrt{\T}}\right\|^{r}_{L^q(0,T;L^p)}\right)^{\f1 r}.
\]
The corresponding zero-mean versions are defined as:
\[
 \|g\|^{l;\zeta}_{\underline{B}_{p,r}^s}:=\|\underline{g}\|^{l;\zeta}_{B_{p,r}^s},\ \ \|\Psi\|^{l;\zeta}_{\widetilde{L}^q_T(\underline{B}_{p,r}^s)}:=\|\underline{\Psi}\|^{l;\zeta}_{\widetilde{L}^q_T(B_{p,r}^s)}.
\]

Finally, we introduce Bony's para-product decomposition. For $g,h\in\mathcal{S}'(\mathbb{T}_b^d)$, we have
\begin{align}
gh=T_{g}h+T_{h}g+R(g,h)+\f{\hat{g}_0}{\sqrt{|\mathbb{T}_b^d|}}\cdot\f{\hat{h}_0}{\sqrt{|\mathbb{T}_b^d|}}
\nonumber
\end{align}
where
\begin{align}
\quad T_{g}h:=\sum\limits_{j\in\mathbb{Z}}S_{j-2} g\Delta_{j} h,\quad
R(g,h):=\sum\limits_{|j-j'|\leq2}\Delta_{j} g\Delta_{j'} h.
\nonumber
\end{align}
We additionally define $T'_{g}h:=T_{h}g+R(g,h)$.

The following almost orthogonality properties are well-known:
\beq\label{ao1}
\Delta_{j'}\Delta_{j} g\equiv0, \text{ if }\, |j'-j|\geq2, \ \   \Delta_{j'}( S_{j-2}h\Delta_{j} g)\equiv0,\text{ if }\, |j'-j|\geq3.
\eeq
and
\begin{align}\label{ao2}
\Delta_{j''}(\Delta_{j} g\Delta_{j'} h)\equiv0, \text{ if } j''-j\geq5,\,|j'-j|\leq2.
\end{align}

\begin{Remark} We introduce the norms $\|\cdot\|_{L^q_T(\underline{B}_{p,r}^s)}$, $\|\cdot\|_{\widetilde{L}^{q}_T(\underline{B}_{p,r}^s)}$ and $\|\cdot\|^{l;\zeta}_{\widetilde{L}^q_T(\underline{B}_{p,r}^s)}$ with the primary aim of eliminating the influence of the zero mode. As evident from their definitions, these norms satisfy the following relations:
\[
\|\Psi\|_{L^q_T(\underline{B}_{p,r}^s)}=\left\|\left(\sum\limits_{j\in\mathbb{Z}}2^{jsr}\|\Delta_{j} \Psi(\cdot)\|^{r}_{L^p(\mathbb{T}_b^d)}\right)^{\f1 r}\right\|_{L^q(0,T)},
\]
\[
\|\Psi\|_{\widetilde{L}^{q}_T(\underline{B}_{p,r}^s)}=\left(\sum\limits_{j\in\mathbb{Z}}2^{jsr}\|\Delta_{j} \Psi\|^{r}_{L^q(0,T;L^p)}\right)^{\f1 r},\ \
\|\Psi\|^{l;\zeta}_{\widetilde{L}^q_T(\underline{B}_{p,r}^s)}=\left(\sum\limits_{2^j<\zeta}2^{jsr}\|\Delta_{j} \Psi\|^{r}_{L^q(0,T;L^p)}\right)^{\f1 r}.
\]
\end{Remark}

\subsection{The main theorem}
Let $0 < T \leq \infty$ and  let $B$ be a Besov space. We denote by $C_T(B) $ the space of continuous bounded $B$-valued functions defined on $[0, T]$ when $0 < T < \infty$, and on $[0, \infty)$ when $T = \infty$. Furthermore,we define  $\widetilde{C}_T(B) := C_T(B)\cap \widetilde{L}^{\infty}_T(B)$. 
From the definition of $L$, $(\text{Ker}L)^\perp=\{(a, \Grad g): a\in L^2(\mathbb{T}_b^d), \intTd a=0 \text{ and } g\in H^1(\mathbb{T}_b^d)\}$.
Our main result reads
\begin{Theorem}\label{main} Let $0<T_0\leq\infty$, $0<\vartheta<\f1 2$. Consider the initial data $(a_0,\vu_0)$ satisfy \eqref{cond1}, with external force $f\in L^1_{T_0}(B^{\f d 2-1}_{2,1}(\mathbb{T}_b^d))$ such that
\[
 \mathbb{Q}f\in C_{T_0}(H^{-S}(\mathbb{T}_b^d)) \text{ and }
 \p_t\mathbb{Q}f\in L^1_{T_0}(H^{-S}(\mathbb{T}_b^d)) \text{ for some } S>0.
\]
Assume system \eqref{equ3} with initial data $\vv_0=\mathbb{P}\vu_0$ and forcing term  $h=\mathbb{P}f$ admits a solution $\vv\in \left(\widetilde{C}_{T_0}(B^{\f d 2-1}_{2,1}(\mathbb{T}_b^d)) \cap  L^1_{T_0}(\underline{B}^{\f d 2+1}_{2,1}(\mathbb{T}_b^d))\right)^d$. Then:
\begin{enumerate}
  \item {\bf{Existence:}} System \eqref{equ4} with initial data $\vV_0=
\left(
\begin{array}{l}
\ \ a_0\\
\mathbb{Q}\vu_0
\end{array}
\right)$  has a solution $\vV\in (\text{Ker}L)^\perp\cap\left(\widetilde{C}_{T_0}(B^{\f d 2-1}_{2,1}(\mathbb{T}_b^d)) \cap  L^1_{T_0}(B^{\f d 2+1}_{2,1}(\mathbb{T}_b^d))\right)^{d+1}$.
  \item {\bf{Uniform Existence, Uniform Energy Estimates and Uniform Bounds on Density:}} There exist a
  $\ep_0(d,\mu,\nu, a_0, \vu_0, \vv,\vV,f)>0$ such that for any $\ep\leq\ep_0$, system \eqref{equ2} posses a unique solution $(\vae,\vue)$ satisfying $\ep\|\vae\|_{L^\infty((0,T_0)\times \mathbb{R}^d)}\leq \f1 2$ and
  \[
  (\vae,\vue)\in \widetilde{C}_{T_0}(\underline{B}^{\f d 2}_{2,1}(\mathbb{T}_b^d))\times \left(\widetilde{C}_{T_0}(B^{\f d 2-1}_{2,1}(\mathbb{T}_b^d)) \cap  L^1_{T_0}(\underline{B}^{\f d 2+1}_{2,1}(\mathbb{T}_b^d))\right)^d,
  \]
with the following estimates:
\begin{align*}
 \ep\|\vae\|^{h;\f{\eta_0}\ep}_{\widetilde{L}^\infty_T(B^{\f d 2}_{2,1})}+\f1 \ep \|\vae\|^{h;\f{\eta_0}\ep}_{L^1_T(B^{\f d 2}_{2,1})}&+\|\vae\|^{l;\f{\eta_0}\ep}_{\widetilde{L}^\infty_T(B^{\f d 2-1}_{2,1})\cap L^1_T(B^{\f d 2+1}_{2,1})}\\
 &+\|\vue\|_{\widetilde{L}^\infty_T(B^{\f d 2-1}_{2,1})\cap L^1_T(\underline{B}^{\f d 2+1}_{2,1})}\leq C(d,\mu,\nu,\vv,\vV,f),
\end{align*}
where the positive constant $\eta_0$, which depends on $d$ and $\nu$, is the one introduced in Lemma \ref{lema4.2}.
  \item {\bf{Decay Estimates:}} The solution satisfies:
  \begin{align*}
  \limsup\limits_{\ep\rightarrow0} \Bigg(&\f{\|\vVe-\vV\|_{\widetilde{L}^\infty_T(H^{\f d 2-1-\vartheta})\cap L^2_T(H^{\f d 2-\vartheta})}}{\tau_{\theta}(\ep)}\\
  &\quad+\f{\|\mathbb{P}\vue-\vv\|_{\widetilde{L}^\infty_{T_0}(B^{\f d 2-1-\vartheta}_{2,1})\cap L^1_{T_0}(\underline{B}^{\f d 2+1-\vartheta}_{2,1})}}{\ep^{\f{\theta}{1+\theta}}}\Bigg) \leq C(\theta,d,\mu,\nu,\vv,\vV,f),
  \end{align*}
where $\tau_{\vartheta}(\ep)$ is monotonically increasing with respect to $\ep$, and $\tau_{\vartheta}(\ep)$ tends to zero as $\ep$ tends to zero\footnote{Danchin \cite{D1} pointed out that  $\tau_{\vartheta}(\ep)$ may tend to zero slower than any power of $\ep$, depending on the periodic parameter vector $b=(b_1,b_2,\cdots,b_d)$.}.
  \item {\bf{Vanishing Limits:}} Moreover, we have:
  \[
 \lim_{\ep\rightarrow0}\Bigg(\ep\|\vae\|_{\widetilde{L}^\infty_{T_0}(B^{\f d 2}_{2,1})}+\|\vVe-\vV\|_{\widetilde{L}^\infty_{T_0}(B^{\f d 2-1}_{2,1})\cap \widetilde{L}^2_{T_0}(B^{\f d 2}_{2,1})}+\|\mathbb{P}\vue-\vv\|_{\widetilde{L}^\infty_{T_0}(B^{\f d 2-1}_{2,1})\cap L^1_{T_0}(\underline{B}^{\f d 2+1}_{2,1})}\Bigg)=0.
\]
\end{enumerate}

\end{Theorem}
\begin{Remark} The two-dimensional case of Theorem \ref{main} is of particular significance. Specifically, assume that $(a_0,\vu_0)$ satisfy \eqref{cond1} and $f\in L^1_\infty(B^{\f d 2-1}_{2,1}(\mathbb{T}_b^d))$ with $d=2$,
 then by the  global well-posedness result established in Theorem 6.3 of \cite{D3}, system \eqref{equ3} admits a global solution $\vv\in \left(\widetilde{C}_{\infty}(B^{0}_{2,1}(\mathbb{T}_b^2)) \cap  L^1_\infty(\underline{B}^{2}_{2,1}(\mathbb{T}_b^2))\right)^2$.
Furthermore, under the additional regularity assumptions,
 \[
 \mathbb{Q}f\in C_{\infty}(H^{-S}(\mathbb{T}_b^2)) \text{ and }
 \p_t\mathbb{Q}f\in L^1_\infty(H^{-S}(\mathbb{T}_b^2)) \text{ for some } S>0,
\]
 we can also establish existence, uniform existence, uniform energy estimates, uniform bounds on density, decay estimates and vanishing limits, with $(d,T_0)$ replaced by $(2,\infty)$.
\end{Remark}

\begin{Remark} Regarding the vanishing limits, we would like to emphasize two key findings\footnote{Fujii \cite{F} also established similar vanishing limit results for $\vae$, $\mathbb{Q}\vue$, and $\mathbb{P}\vue - \vv$ in the whole-space framework.}:
\begin{itemize}
  \item Danchin \cite{D1} established, for small initial data with critical regularity, the convergence to zero of the quantity $\|\vVe-\vV\|_{\widetilde{L}^\infty_{T_0}(B^{\f d 2-1-\theta}_{2,1})\cap L^2_{T_0}(B^{\f d 2-\theta}_{2,1})},$ and the quantity
  \[
  \|\mathbb{P}\vue-\vv\|_{\widetilde{L}^\infty_{T_0}(B^{\f d 2-1-\theta}_{2,1})\cap L^1_{T_0}(\underline{B}^{\f d 2+1-\theta}_{2,1})},
  \]
   where $\vartheta>0$. Our results thus extend those of Danchin in terms of both the size of initial data and spatial regularity.  \item  In the whole-space framework, Ogino \cite{O} recently established the vanishing property of the quantity $\ep\|\vae\|_{\widetilde{L}^\infty_{T_0}(B^{\f d 2}_{2,1})}$ for small initial data with critical regularity. However, whether this property holds for large initial data with critical regularity remains unknown. Our result provides an affirmative answer to this question in the periodic domain case.
\end{itemize}
\end{Remark}

\begin{Remark} Although we primarily consider initial data and external force $(a_0,\vu_0,f)$ that are independent of $\ep$, our method remains applicable to $\ep$-dependent cases $(a_0^\ep,\vue_0,f^\ep)$, provided they satisfy certain convergence conditions. Specifically, assume that $(a_0^\ep, \vue_0, f^\ep)$  satisfy the same hypotheses as $(a_0,\vu_0, f)$
in Theorem \ref{main}, and that $\vV_0\in (\text{Ker}L)^\perp\cap \Big(B_{2,1}^{\f{d} 2-1}(\mathbb{T}_b^d)\Big)^{d+1}$.
Suppose also that system \eqref{equ3}, with initial data $\vv_0\in \Big(B_{2,1}^{\f{d} 2-1}(\mathbb{T}_b^d)\Big)^d$ satisfying $\Div\vv_0=0$ and external forces  $h\in L^1_{T_0}(B^{\f d 2-1}_{2,1}(\mathbb{T}_b^d))$ satisfying $\Div h=0$, admits a solution $\vv\in \left(\widetilde{C}_{T_0}(B^{\f d 2-1}_{2,1}(\mathbb{T}_b^d)) \cap  L^1_{T_0}(\underline{B}^{\f d 2+1}_{2,1}(\mathbb{T}_b^d))\right)^d$. Furthermore, assume  following convergence conditions hold:
\[
\limsup\limits_{\ep\rightarrow0}\Big( \|\mathbb{Q}f^\ep\|_{L_{T_0}^\infty(H^{-S})\cap L_{T_0}^1(B^{\f d 2-1}_{2,1})}+\|\p_t\mathbb{Q}f^\ep\|_{L_{T_0}^1(H^{-S})}\Big)<\infty,
\]
\[
\lim_{(\zeta,\ep)\rightarrow(\infty,0)}\|\mathbb{Q}f^\ep\|^{h;\zeta}_{ L_{T_0}^1(B^{\f d 2-1}_{2,1})}=0,
\]
\[
\lim_{\ep\rightarrow0} \ep\|a_0^\ep\|_{B^{\f d 2}_{2,1}}+\|(a_0^\ep,\mathbb{Q}\vue_0)-\vV_0\|_{B^{\f d 2-1}_{2,1}}+\|\mathbb{P}\vue_0-\vv_0\|_{B^{\f d 2-1}_{2,1}}+\|\mathbb{P}f^\ep-h\|_{L_{T_0}^1(B^{\f d 2-1}_{2,1})}=0
\]
and for any $\theta\in (0,\f1 2)$,
  \[
\limsup\limits_{\ep\rightarrow0} \Bigg(\f{\|(a_0^\ep,\mathbb{Q}\vue_0)-\vV_0\|_{H^{\f d 2 -1-\theta}}}{\tau_{\theta}(\ep)}+\f{\|\mathbb{P}\vue_0-\vv_0\|_{B^{\f d 2-1-\theta}_{2,1}}}{\ep^{\f{\theta}{1+\theta}}}+\f{\|\mathbb{P}f^\ep-h\|_{L_{T_0}^1(B^{\f d 2-1-\theta}_{2,1})}}{\ep^{\f{\theta}{1+\theta}}}\Bigg)<\infty.
\]
Then the existence, uniform existence, uniform energy estimates, uniform bounds on density, decay estimates and vanishing limits all remain valid.
\end{Remark}

\section{Proof strategy and  difficulties} \label{newsec}
\medskip
\medskip
The proof of Theorem \ref{main} relies on two fundamental components: demonstrating the convergence of
$\ep\|\vae\|_{\widetilde{L}^\infty_T(B^{\f d 2}_{2,1})}$ to zero (as $\ep \to 0$) and establishing uniform estimates in  $\ep$ and $T$ for the
 quantity
 $X^\ep(T)+P^\ep(T)$, where
 \[
X^\ep(T):=\ep\|\vae\|^{h;\f{\eta_0}\ep}_{\widetilde{L}^\infty_T(B^{\f d 2}_{2,1})}+\f1 \ep \|\vae\|^{h;\f{\eta_0}\ep}_{L^1_T(B^{\f d 2}_{2,1})}+\|\vae\|^{l;\f{\eta_0}\ep}_{\widetilde{L}^\infty_T(B^{\f d 2-1}_{2,1})\cap L^1_T(B^{\f d 2+1}_{2,1})}+\|\mathbb{Q}\vue\|_{\widetilde{L}^\infty_T(B^{\f d 2-1}_{2,1})\cap L^1_T(B^{\f d 2+1}_{2,1})},
\]
 \[
 P^\ep(T):=\|\mathbb{P}\vue\|_{\widetilde{L}^\infty_T(B^{\f d 2-1}_{2,1})\cap L^1_T(\underline{B}^{\f d 2+1}_{2,1})}.
 \]
 A natural approach to handle $\ep\|\vae\|_{\widetilde{L}^\infty_T(B^{\f d 2}_{2,1})}$ is  to decompose it into two parts (see Danchin \cite{D1}),  namely:
 \[
\ep\|\vae\|_{\widetilde{L}^\infty_T(B^{\f d 2}_{2,1})}\leq \ep\|\vae\|^{h;\f{\eta_0}\ep}_{\widetilde{L}^\infty_T(B^{\f d 2}_{2,1})}+\ep\|\vae\|^{l;\f{\eta_0}\ep}_{\widetilde{L}^\infty_T(B^{\f d 2}_{2,1})}\lesssim \ep\|\vae\|^{h;\f{\eta_0}\ep}_{\widetilde{L}^\infty_T(B^{\f d 2}_{2,1})}+\|\vae\|^{l;\f{\eta_0}\ep}_{\widetilde{L}^\infty_T(B^{\f d 2-1}_{2,1})}.
\]
However, since the initial data is large, the quantity $\|\vae\|^{l;\f{\eta_0}\ep}_{\widetilde{L}^\infty_T(B^{\f d 2-1}_{2,1})}$ also  becomes large, which renders the  high-low frequency decomposition ineffective. Therefore, we adopt the high-medium-low frequency analysis framework introduced by Fujii \cite{F}.
Specifically, for any $0<\zeta<\f{\eta_0}\ep$, we have
\begin{align*}
\ep\|\vae\|_{\widetilde{L}^\infty_T(B^{\f d 2}_{2,1})}&\leq  \ep\|\vae\|^{h;\f{\eta_0}\ep}_{\widetilde{L}^\infty_T(B^{\f d 2}_{2,1})}+\ep\|\vae\|^{m;\zeta,\f{\eta_0}\ep}_{\widetilde{L}^\infty_T(B^{\f d 2}_{2,1})}+\ep\|\vae\|^{l;\zeta}_{\widetilde{L}^\infty_T(B^{\f d 2}_{2,1})}\\
&\lesssim  \ep\|\vae\|^{h;\f{\eta_0}\ep}_{\widetilde{L}^\infty_T(B^{\f d 2}_{2,1})}+\|\vae\|^{m;\zeta,\f{\eta_0}\ep}_{\widetilde{L}^\infty_T(B^{\f d 2-1}_{2,1})}+\zeta\ep \|\vae\|^{l;\zeta}_{\widetilde{L}^\infty_T(B^{\f d 2-1}_{2,1})}
\end{align*}
An intuitive advantage of this decomposition is that for any $\delta>0$, we can choose sufficiently large  $\zeta$ and sufficiently small $\ep$ such that
\[
\ep\|a_0\|^{h;\f{\eta_0}\ep}_{B^{\f d 2}_{2,1}}+\|a_0\|^{m;\zeta,\f{\eta_0}\ep}_{B^{\f d 2-1}_{2,1}}+ \zeta\ep\leq\delta.
\]
Based on our decomposition strategy, to establish the vanishing property of $\ep\|\vae\|_{\widetilde{L}^\infty_T(B^{\f d 2}_{2,1})}$, we must study the following quantity (see Lemma \ref{lema4.2} and Lemma \ref{lema4.3}):
\[
\ep\|\vae\|^{h;\f{\eta_0}\ep}_{\widetilde{L}^\infty_T(B^{\f d 2}_{2,1})}+\f1 \ep \|\vae\|^{h;\f{\eta_0}\ep}_{L^1_T(B^{\f d 2}_{2,1})}+\|\mathbb{Q}\vue\|^{h;\f{\eta_0}\ep}_{\widetilde{L}^\infty_T(B^{\f d 2-1}_{2,1})\cap L^1_T(B^{\f d 2+1}_{2,1}) }+\|(\vae,\mathbb{Q}\vue)\|^{m;\zeta,\f{\eta_0}\ep}_{\widetilde{L}^\infty_T(B^{\f d 2-1}_{2,1})\cap L^1_T(B^{\f d 2+1}_{2,1}) }.
\]
Building on this analysis, we construct the quantity
  $D^{\zeta,\ep}(T)\footnote{ Fujii \cite{F} introduced a similar quantity $A_{q,e}^{\ep,\alpha}$,  but with consistent variables for both high-medium and low frequencies. This work also inspired our construction of the quantity $D^{\zeta,\ep}(T)$}:=[\vae,\vue]^{\zeta,\ep}_{h,m}(T)+ [\vVe-\vV,\mathbb{P}\vue-\vv]^{\zeta,\ep}_{l}(T)$,   whose low-frequency component is required to remain small. Accordingly, in this part we work with the difference variables $\vVe - \vV$ and $\mathbb{P}\vue - \vv$, where $\vV$ and $\vv$ denote the solutions to systems \eqref{equ4} and \eqref{equ3},  respectively. Here
\begin{align*}
     [\vae,\vue]^{\zeta,\ep}_{h,m}(T):=&\ep\|\vae\|_{\widetilde{L}^\infty_T(B^{\f d 2}_{2,1})}+\ep\|\vae\|^{h;\f{\eta_0}\ep}_{\widetilde{L}^\infty_T(B^{\f d 2}_{2,1})}+\f1 \ep \|\vae\|^{h;\f{\eta_0}\ep}_{L^1_T(B^{\f d 2}_{2,1})}+\|\mathbb{Q}\vue\|^{h;\f{\eta_0}\ep}_{\widetilde{L}^\infty_T(B^{\f d 2-1}_{2,1})\cap L^1_T(B^{\f d 2+1}_{2,1}) }\\
     &+\|(\vae,\mathbb{Q}\vue)\|^{m;\zeta,\f{\eta_0}\ep}_{\widetilde{L}^\infty_T(B^{\f d 2-1}_{2,1})\cap L^1_T(B^{\f d 2+1}_{2,1}) }+ \|\mathbb{P}\vue\|^{h;\zeta}_{\widetilde{L}^\infty_T(B^{\f d 2-1}_{2,1})\cap L^1_T(B^{\f d 2+1}_{2,1}) },
  \end{align*} and
for $W\in \mathbb{R}^{d+1}$ and $\vu\in \mathbb{R}^d$  \begin{align*}
 [W,\vu]^{\zeta,\ep}_{l}(T):=\|W\|^{l;\zeta}_{\widetilde{L}^\infty_T(B^{\f d 2-1}_{2,1})\cap \widetilde{L}^2_T(B^{\f d 2}_{2,1}) }
 +\|\vu\|^{l;\zeta}_{\widetilde{L}^\infty_T(B^{\f d 2-1}_{2,1})\cap L^1_T(\underline{B}^{\f d 2+1}_{2,1}) }.
  \end{align*}
 To bridge the quantities $D^{\zeta,\ep}(T)$ and  $X^\ep(T)+P^\ep(T)$, we introduce the auxiliary quantity
     $Y^{\zeta,\ep}(T):=[\vae,\vue]^{\zeta,\ep}_{h,m}(T)+ [(\vae,\mathbb{Q}\vue),\mathbb{P}\vue]^{\zeta,\ep}_{l}(T)$. From \eqref{cons2}, we obtain
      \[
       Y^{\zeta,\ep}(T)\lesssim D^{\zeta,\ep}(T)+\|(\vv,\vV)\|^{l;\zeta}_{\widetilde{L}^\infty_T(B^{\f d 2-1}_{2,1})\cap L^1_T(\underline{B}^{\f d 2+1}_{2,1}) }.
       \]
Furthermore, Proposition \ref{pro3} establishes that
\[
X^\ep(T)+P^\ep(T)\lesssim \|(a_0,\vu_0)\|+f+\Big(Y^{\zeta,\ep}(T)\Big)^2+\Big(Y^{\zeta,\ep}(T)\Big)^3.
\]
These two inequalities together complete the connection between $D^{\zeta,\ep}(T)$ and  $X^\ep(T)+P^\ep(T)$.

The proof of the smallness of the quantity  $D^{\zeta,\ep}(T)$  (for sufficiently large $\zeta$ and $\ep$) now occupies a central position in the argument.  Due to the complex behavior of acoustic waves, it is difficult to directly establish the smallness of the quantity $\|\vVe-\vV\|^{l;\zeta}_{\widetilde{L}^\infty_T(B^{\f d 2-1}_{2,1})\cap L^1_T(B^{\f d 2+1}_{2,1}) }$. To overcome this difficulty, we employ an indirect approach by exploiting the decay properties of $\vVe-\vV$ in low-regularity spaces (This is also why the $L^2$-norm, rather than the $L^1$-norm, is chosen in the time variable):
  \begin{align}\label{key}
   \|\vVe-\vV\|^{l;\zeta}_{\widetilde{L}^\infty_T(B^{\f d 2-1}_{2,1})\cap \widetilde{L}^2_T(B^{\f d 2}_{2,1}) }&\lesssim \zeta^{2\vartheta}\|\vVe-\vV\|^{l;\zeta}_{\widetilde{L}^\infty_T(B^{\f d 2-1-2\vartheta}_{2,1})\cap \widetilde{L}^2_T(B^{\f d 2-2\vartheta}_{2,1}) }\nonumber\\
 &\lesssim C(\theta)\zeta^{2\vartheta}\|\vVe-\vV\|_{\widetilde{L}^\infty_T(H^{\f d 2-1-\vartheta})\cap L^2_T(H^{\f d 2-\vartheta})},
 \end{align}
 where  the decay properties of $\|\vVe-\vV\|_{\widetilde{L}^\infty_T(H^{\f d 2-1-\vartheta})\cap L^2_T(H^{\f d 2-\vartheta})}$ with respect to $\ep$ is established in Proposition \ref{pro1} and  we employed \eqref{minski2} at the second inequality.
 The same method is applicable to $\|\mathbb{P}\vue-\vv\|^{l;\zeta}_{\widetilde{L}^\infty_T(B^{\f d 2-1}_{2,1})\cap L^1_T(\underline{B}^{\f d 2+1}_{2,1}) }$ and decay properties in low-regularity spaces is shown in Proposition \ref{pro2}.

 In the process of analysis $ [\vae,\vue]^{\zeta,\ep}_{h,m}(T)$, one of the key tasks is to establish corresponding estimates for nonlinear terms of the form  $
   \int_0^T\|fg\|^{h;\zeta}_{B^{\f d 2-1}_{2,1}} \dt. $ During this process, the coupling of frequencies between $f$ and $g$ introduces several undesirable terms.
   According to Lemma \ref{lema4.1}, this integral satisfies the following estimate:
   \begin{align*} 
        \int_0^T\|fg\|^{h;\zeta}_{B^{\f d 2-1}_{2,1}} \dt &\leq  \int_0^T\|T_f g\|^{h;\zeta}_{B^{\f d 2-1}_{2,1}} \dt+ \int_0^T\|T_g f\|^{h;\zeta}_{B^{\f d 2-1}_{2,1}} \dt+ \int_0^T\|R(f, g)\|^{h;\zeta}_{B^{\f d 2-1}_{2,1}} \dt\\
        &\leq R_1+R_2+R_3,
   \end{align*}
   \[
   R_1:=\int_0^T \|f\|_{B^{s_1}_{p_1,r_1}} \|g\|^{h;\zeta_2}_{B^{s_2}_{p_2,r_2}}\dt, R_2:=\int_0^T \|f\|^{h;\zeta_3}_{B^{s_3}_{p_3,r_3}} \|g\|_{B^{s_4}_{p_4,r_4}}\dt, R_3:=\int_0^T \|f\|^{h;\zeta_5}_{B^{s_5}_{p_5,r_5}} \|g\|^{h;\zeta_6}_{B^{s_6}_{p_6,r_6}}\dt.
   \]
 Due to the large values of  $\|f\|^{l;\zeta_1}_{B^{s_1}_{p_1,r_1}}$ and  $\|g\|_{B^{s_4}_{p_4,r_4}}$, the terms $R_1$ and $R_2$
 have adverse effects. Notably, the term $\int_0^T\|R(f, g)\|^{h;\zeta}_{B^{\f d 2-1}_{2,1}} \dt$ admits a better estimate (namely, $\int_0^T\|R(f, g)\|^{h;\zeta}_{B^{\f d 2-1}_{2,1}} \dt \leq R_3$ ) compared to those of $\int_0^T\|T_f g\|^{h;\zeta}_{B^{\f d 2-1}_{2,1}} \dt$ and $\int_0^T\|T_g f\|^{h;\zeta}_{B^{\f d 2-1}_{2,1}} \dt$---a result that is somewhat unexpected. By fully exploiting the high-medium-low frequency structure, we perform a refined frequency decomposition of the terms  $R_1$ and $R_2$, and thereby derive a suitable estimate for $[\vae,\vue]^{\zeta,\ep}_{h,m}(T)$ in Proposition \ref{pro4}.  By further combining the decay properties provided in Propositions \ref{pro1}-\ref{pro2}, the smallness of $[\vae,\vue]^{\zeta,\ep}_{h,m}(T)$  can be confirmed.

\section{Full demonstration of Theorem \ref{main}}\label{complete}
We divide the proof into the following six steps, and the relevant notations are also listed below.
\[
X:=\|\vV\|_{\widetilde{L}^\infty_T(B^{\f d 2-1}_{2,1})\cap L^1_T(B^{\f d 2+1}_{2,1})},\ \ X_0:=\|a_0\|^{l;\f{\eta_0}\ep}_{B^{\f d 2-1}_{2,1}}+\ep\|a_0\|^{h;\f{\eta_0}\ep}_{B^{\f d 2}_{2,1}}+\|\mathbb{Q}\vu_0\|_{B^{\f d 2-1}_{2,1}}, P_0:=\|\mathbb{P}\vu_0\|_{B^{\f d 2-1}_{2,1}},
\]
\[
\mathbb{Q}_f:=\|\mathbb{Q}f\|_{L^1_{T_0}(B^{\f d 2-1}_{2,1})}+\|\mathbb{Q}f\|_{L^\infty_{T_0}(H^{-S})}+\|\p_t\mathbb{Q}f\|_{L^1_{T_0}(H^{-S})}, \ \ \mathbb{Q}_f(0):=\|\mathbb{Q}f(0)\|_{H^{-S}}
\]
\[
P^\ep(T):=\|\mathbb{P}\vue\|_{\widetilde{L}^\infty_T(B^{\f d 2-1}_{2,1})\cap L^1_T(\underline{B}^{\f d 2+1}_{2,1})},\ \ P:=\|\vv\|_{\widetilde{L}^\infty_T(B^{\f d 2-1}_{2,1})\cap L^1_T(\underline{B}^{\f d 2+1}_{2,1})},\ \  \mathbb{P}_f:=\|\mathbb{P}f\|_{L^1_{T_0}(B^{\f d 2-1}_{2,1})},
\]
\[
\vZe_\vartheta(T):=\|\vVe-\vV\|_{\widetilde{L}^\infty_T(H^{\f d 2-1-\vartheta})\cap L^2_T(H^{\f d 2-\vartheta})}, \ \
\vWe_\vartheta(T):=\|\mathbb{P}\vue-\vv\|_{\widetilde{L}^\infty_T(B^{\f d 2-1-\vartheta}_{2,1})\cap L^1_T(\underline{B}^{\f d 2+1-\vartheta}_{2,1})}.\]

{\bf{First step. The existence spaces of $\vV$.}} Since  system \eqref{equ3} admits a solution
\begin{align} \label{new1}
   \vv\in \left(\widetilde{C}_{T_0}(B^{\f d 2-1}_{2,1}(\mathbb{T}_b^d)) \cap  L^1_{T_0}(\underline{B}^{\f d 2+1}_{2,1}(\mathbb{T}_b^d))\right)^d,
\end{align}
it follows from Theorem 8.2 in \cite{D1} that system \eqref{equ4}  has a solution
\begin{align}\label{new2} 
   \vV\in (\text{Ker}L)^\perp\cap\left(\widetilde{C}_{T_0}(B^{\f d 2-1}_{2,1}(\mathbb{T}_b^d)) \cap  L^1_{T_0}(B^{\f d 2+1}_{2,1}(\mathbb{T}_b^d))\right)^{d+1}.
\end{align}

{\bf{Second step. The decay properties of $\vVe-\vV$.}}  This result is established in the following proposition, with a detailed proof provided in Subsection \ref{newsec1}.
\begin{Proposition}\label{pro1} Let $0<T\leq T_0$, $0<\vartheta<\f1 2$ and $0<\zeta<\f{\eta_0}\ep$. Suppose the system \eqref{equ2} admits a  solution $(\vae,\vue)$ such that $\ep\|\vae\|_{L^\infty((0,T)\times \mathbb{R}^d)}\leq \f1 2$ and
\begin{align}\label{new3} 
    (\vae,\vue)\in \widetilde{C}_{T}(\underline{B}^{\f d 2}_{2,1}(\mathbb{T}_b^d))\times \left(\widetilde{C}_{T}(B^{\f d 2-1}_{2,1}(\mathbb{T}_b^d)) \cap  L^1_{T}(\underline{B}^{\f d 2+1}_{2,1}(\mathbb{T}_b^d))\right)^d,
\end{align}

Then the following estimate holds:
\begin{align*}
  \vZe_\vartheta(T)\leq& C(\vartheta,d,\mu,\nu)e^{C(\vartheta,d,\mu,\nu)(P^2+X^2)}\Bigg(\Big(D^{\zeta,\ep}(T)+\|\vV\|^{h;\zeta}_{\widetilde{L}^2_{T_0}(B^{\f d 2}_{2,1})} \Big) \vZe_\vartheta(T)\\
  &+\vWe_\vartheta(T)(P+Y^{\zeta,\ep}(T))+\widetilde{\tau}_{\vartheta}(\ep)\times\Big(\mathbb{Q}_f(0)+X_0+\mathbb{Q}_f+(\mathbb{Q}_f)^2+X^2_0+P_0^2+(\mathbb{P}_f)^2\\
  &+Y^{\zeta,\ep}(T)+X^2+P^2+(Y^{\zeta,\ep}(T)^2+X^3+P^3+(Y^{\zeta,\ep}(T)^3+X^4+(Y^{\zeta,\ep}(T))^4\Big)\Bigg),
\end{align*}
where the function $\widetilde{\tau}_{\vartheta}(\ep)$ is monotonically increasing with respect to $\ep$ and satisfies  $\widetilde{\tau}_{\vartheta}(\ep)\rightarrow0$ as $\ep\rightarrow0$.
\end{Proposition}

{\bf{Third step. The decay properties of $\mathbb{P}\vue-\vv$.}}  This result is established in the following proposition, with a detailed proof provided in Subsection \ref{newsec2}.
\begin{Proposition}\label{pro2}Let $0<T\leq T_0$, $0<\vartheta<\f1 2$ and $0<\zeta<\f{\eta_0}\ep$. Suppose the system \eqref{equ2} admits a  solution $(\vae,\vue)$ satisfying \eqref{new3} with $\ep\|\vae\|_{L^\infty((0,T)\times \mathbb{R}^d)}\leq \f1 2$. Then the following estimate holds:
\begin{align*}
\vWe_\vartheta(T)\leq& C(\vartheta,d,\mu,\nu)e^{C(\vartheta,d,\mu,\nu)(P+X^\ep(T))}\Bigg( \Big(D^{\zeta,\ep}(T)+\|\vv\|^{h;\zeta}_{\widetilde{L}^\infty_{T_0}(B^{\f d 2-1}_{2,1})} \Big)\vWe_\vartheta(T)+\ep^{\f{\vartheta}{1+\vartheta}}\Big((X_0)^2\\
&+(P_0)^2+(\mathbb{P}_f)^2+(\mathbb{Q}_f)^2+P^2+(Y^{\zeta,\ep}(T))^2+P^3+(Y^{\zeta,\ep}(T))^3+P^4+(Y^{\zeta,\ep}(T))^4\Big)\Bigg).
\end{align*}
\end{Proposition}

{\bf{Fourth step. The control relation
between quantity $Y^{\zeta,\ep}(T)$  and quantity $X^\ep(T)+P^\ep(T)$.}} This result is established in the following proposition, with a detailed proof provided in Subsection \ref{newsec3}
\begin{Proposition} \label{pro3}Let $0<T\leq T_0$ and $0<\zeta<\f{\eta_0}\ep$.  Suppose the system \eqref{equ2} admits a  solution  $(\vae,\vue)$ satisfying \eqref{new3} with $\ep\|\vae\|_{L^\infty((0,T)\times \mathbb{R}^d)}\leq \f1 2$. Then the following estimate holds:
\begin{align*}
 X^\ep(T)+P^\ep(T)\leq C(d,\mu,\nu)\Big(X_0+P_0+\mathbb{P}_f+\mathbb{Q}_f+(Y^{\zeta,\ep}(T))^2+(Y^{\zeta,\ep}(T))^3\Big).
\end{align*}
\end{Proposition}

{\bf{ Fifth step. The estimates for quantity $ [\vae,\vue]^{\zeta,\ep}_{h,m}(T)$.}} This result is established in the following proposition, with a detailed proof provided in Subsection \ref{newsec4}.
\begin{Proposition} \label{pro4}Let $0<T\leq T_0$ and $0<\zeta<\f{\eta_0}\ep$. Suppose the system \eqref{equ2} admits a  solution  $(\vae,\vue)$ satisfying \eqref{new3} with $\ep\|\vae\|_{L^\infty((0,T)\times \mathbb{R}^d)}\leq \f1 2$. Then the following estimate holds:
\begin{align*}
 [\vae,\vue]^{\zeta,\ep}_{h,m}(T)\leq& C(d,\mu,\nu)e^{C(d,\mu,\nu)( Y^{\zeta,\ep}(T)+ (Y^{\zeta,\ep}(T))^2)}\Bigg(\ep\|\va_0\|^{h;\f{\eta_0}\ep}_{B_{2,1}^{\f d 2}}+
  \|(a_0,\vu_0)\|^{h;\zeta}_{B_{2,1}^{\f d 2-1}}\nonumber\\
  &+Y^{\zeta,\ep}(T)\|(\vv,\vV)\|^{h;\f\zeta{16}}_{\widetilde{L}_{T_0}^2(B^{\f d 2}_{2,1})}+\max\{\f 1 \zeta,\ep,\zeta\ep\}Y^{\zeta,\ep}(T) [\vae,\vue]^{\zeta,\ep}_{h,m}(T)+\|f\|^{h;\zeta}_{L_{T_0}^1(B_{2,1}^{\f d 2-1})}\nonumber\\
  &+(\zeta\ep)(Y^{\zeta,\ep}(T)+(Y^{\zeta,\ep}(T))^2)+  [\vVe-\vV,\mathbb{P}\vue-\vv]^{\zeta,\ep}_{l}(T)Y^{\zeta,\ep}(T)\Bigg).
\end{align*}
\end{Proposition}

{\bf{Sixth step. End of the proof.}}  By synthesizing the results of Propositions \ref{pro1}-\ref{pro4}, we appropriately select the parameters $\zeta$ and $\varepsilon$, and conclude the proof through a bootstrap argument.

  \medskip
  \medskip

{\bf{Proof:}} Using  \eqref{new1}, \eqref{new2} and the dominated convergence theorem, we obtain that
\begin{align}\label{energy}
\lim_{\zeta\rightarrow\infty}\|(\vv,\vV)\|^{h;\zeta}_{\widetilde{L}^\infty_{T_0}(B^{\f d 2-1}_{2,1})\cap L^1_{T_0}(B^{\f d 2+1}_{2,1}) }=0.
\end{align}

There exists a small $\ep^\ast>0$ such that for any $\ep\leq\ep^\ast$,
\begin{align}\label{cons8}
 \ep\|a_0\|_{L^{\infty}}\leq \f 1 8, \ \ \ep\|a_0\|_{B^{\f d 2}_{2,1}}+X_0+P_0\leq C_1 (X+P),
\end{align}
where the constant $C_1$ is independent of $\ep$. Furthermore,
there also exista a constant  $C_2$, independent of $\zeta,\ep,T$, satisfying that
\begin{align}
& \ep\|\vae\|_{L^\infty((0,T)\times \mathbb{R}^d)}\leq C_2\ep\|\vae\|_{\widetilde{L}^\infty_T(B^{\f d 2}_{2,1})},\label{cons1}\\
&Y^{\zeta,\ep}(T)\leq C_2\Big(D^{\zeta,\ep}(T)+\|(\vv,\vV)\|^{l;\zeta}_{\widetilde{L}^\infty_T(B^{\f d 2-1}_{2,1})\cap L^1_T(\underline{B}^{\f d 2+1}_{2,1}) }\Big),\label{cons2}\\
& [\vVe-\vV,\mathbb{P}\vue-\vv]^{\zeta,\ep}_{l}(T)\leq  C_2\zeta^{2\vartheta}(\vWe_\theta(T)+\vZe_\theta(T)),\label{cons7}
\end{align}
where \eqref{cons2} and \eqref{cons7} follows from Lemma \ref{LeD3} and \eqref{key} respectively. Henceforth, throughout this section, the notation $C\geq1$ will denote the positive constant introduced in Propositions \ref{pro1}-\ref{pro4}.

 Let
\begin{align*}
 T^\ast_\ep:=\sup\Big\{0<T<\infty: &\text{ system }\eqref{equ2} \text{ admits a solution } (\vae,\vue)\text{ satisfying \eqref{new3}  } \\
 & \text{with} \inf_{(t,x)\in [0,T]\times \mathbb{R}^d}1+\ep\vae(t,x)>0\Big\}.
\end{align*}
We now proceed to select appropriate parameters $\zeta$ and $\varepsilon$. In doing so, it is necessary to introduce several essential constants. We define
\begin{align*}
  &\tau_\theta(\ep):=\max\{\widetilde{\tau}_\theta(\ep),\ \ \ep^{\f{\theta}{1+\theta}}\},\\
 &C_3:=(C_1+C_2)(X+P),\ \  C_4:=C\Big(C_3+\mathbb{P}_f+\mathbb{Q}_f+(2C_3)^2+(2C_3)^3\Big),\\
 &C_5:=Ce^{C(P+C_4)}\Big((C_3)^2+(C_3)^2+(\mathbb{P}_f)^2+(\mathbb{Q}_f)^2+P^2+(2C_3)^2\\
 &\qquad+P^3+(2C_3)^3+P^4+(2C_3)^4\Big)\Big),\\
 &C_6:=Ce^{C(P^2+X^2)}\Big((P+2C_3)2C_5+\mathbb{Q}_f(0)+C_3+\mathbb{Q}_f+(\mathbb{Q}_f)^2+(C_3)^2+(C_3)^2+(\mathbb{P}_f)^2\\
  &+2C_3+X^2+P^2+(2C_3)^2+X^3+P^3+(2C_3)^3+X^4+(2C_3)^4\Big),
\end{align*}
where the constants $C_4$, $C_5$ and $C_6$ are constructed based on Proposition \ref{pro3}, Proposition \ref{pro2}, and Proposition \ref{pro1}, respectively.

Let $\delta_0\leq1$ be a small positive constant such that
\begin{align}\label{cons}
 C_2\delta_0\leq \f {C_3} 4,\ \ C_2\delta_0\leq \f 1 4,\ \  Ce^{C(P^2+X^2)}3\delta_0\leq\f 1 4,\ \ Ce^{C(P+C_4)}3\delta_0\leq\f 1 4.
\end{align}
Thanks to \eqref{energy}, there exists a large $\zeta_0>0$ such that
\begin{align}
&\|(\vv,\vV)\|^{h;\zeta_0}_{\widetilde{L}^\infty_{T_0}(B^{\f d 2-1}_{2,1})\cap \widetilde{L}^2_{T_0}(B^{\f d 2}_{2,1})}\leq \delta_0,\ \ Ce^{C(2C_3+(2C_3)^2)}\f1{\zeta_0}2C_3\leq \f{1}8, \label{cons3}\\
&Ce^{C(2C_3+(2C_3)^2)}\Big(\|(a_0,\vu_0)\|^{h;\zeta_0}_{B_{2,1}^{\f d 2-1}}+2C_3\|(\vv,\vV)\|^{h;\f{\zeta_0}{16}}_{\widetilde{L}_{T_0}^2(B^{\f d 2}_{2,1})}+\|f\|^{h;\zeta_0}_{L_{T_0}^1(B_{2,1}^{\f d 2-1})}\Big)\leq \f{\delta_0} 8,\label{cons4}
\end{align}
Now we choose a small $\ep_0<\min\{\ep^\ast,\f{\eta_0}{\zeta_0}\}$ such that
\begin{align}
&Ce^{C(2C_3+(2C_3)^2)}(\ep_0+\ep_0\zeta_0)2C_3\leq \f{1}8, \label{cons6}\\
&Ce^{C(2C_3+(2C_3)^2)}\Big(\ep_0\|a_0\|_{B_{2,1}^{\f d 2}}+\ep_0\zeta_0(2C_3+(2C_3)^2)+C_2(\zeta_0)^2(2C_5+2C_6)2C_3\tau_\theta(\ep_0)\Big)\nonumber\\
&+C_2(\zeta_0)^2(2C_5+2C_6)\tau_\theta(\ep_0)\leq \f{\delta_0}8.\label{cons5}
\end{align}
We define
\begin{align*}
T^{\ast\ast}_\ep:=\sup\Big\{0<T<T^{\ast}_{\ep}:& Y^{\zeta_0,\ep}(T)\leq 2C_3,\ \ \vWe_\vartheta(T)\leq 2C_5\ep^{\f{\theta}{1+\theta}},\\
&\vZe_\vartheta(T)\leq 2C_6\tau_\theta(\ep),\ \ D^{\zeta_0,\ep}(T)\leq 2\delta_0. \Big\}.
\end{align*}

From \eqref{cons8}, \eqref{cons4} and \eqref{cons5}, we obtain that for any $\ep\leq\ep_0$, $T^\ast_\ep\geq T^{\ast\ast}_\ep>0$.
We assert that for every $\ep\leq\ep_0$, the inequality $T^\ast_\ep\geq T_0$ holds. We employ proof by contradiction to demonstrate the above inference. Assuming the inference is false, it follows that $T^{\ast}_\ep<T_0$. Then we claim $T^{\ast\ast}_\ep=T^{\ast}_\ep$. Let $0<T<T^{\ast\ast}_\ep$. From \eqref{cons1} and \eqref{cons}, we have that
\begin{align}\label{fina5}
  \ep\|\vae\|_{L^\infty((0,T)\times \mathbb{R}^d)}\leq C_2\ep\|\vae\|_{\widetilde{L}^\infty_T(B^{\f d 2}_{2,1})}\leq C_2D^{\zeta_0,\ep}(T)\leq \f 1 2.
\end{align}

Thanks to \eqref{cons2} and \eqref{cons}, we infer that
\begin{align}\label{fina1}
 Y^{\zeta_0,\ep}(T)\leq C_2D^{\zeta_0,\ep}(T)+C_2(P+X)\leq \f 3 2 C_3.
\end{align}

Using Proposition \ref{pro1} and  \eqref{cons}-\eqref{cons3}, we obtain that
\[
\vZe_\vartheta(T)\leq Ce^{C(P^2+X^2)}3\delta_0\vZe_\vartheta(T)+C_6\tau_\theta(\ep)\leq \f 1 4\vZe_\vartheta(T)+C_6\tau_\theta(\ep),
\]
which implies
\begin{align}\label{fina2}
 \vZe_\vartheta(T)\leq \f 4 3C_6\tau_\theta(\ep).
\end{align}

According to Proposition \ref{pro2} and \eqref{cons}-\eqref{cons3}, we get that
\[
\vWe_\vartheta(T)\leq Ce^{C(P+C_4)}3\delta_0\vWe_\vartheta(T)+C_5\ep^{\f{\vartheta}{1+\vartheta}}\leq \f 1 4 \vWe_\vartheta(T)+C_5\ep^{\f{\vartheta}{1+\vartheta}},
\]
which yields
\begin{align}\label{fina3}
  \vWe_\vartheta(T)\leq \f 4 3C_5\ep^{\f{\vartheta}{1+\vartheta}}.
\end{align}

From Proposition \ref{pro4}, \eqref{cons7} and \eqref{cons3}-\eqref{cons5}, we infer that
\[
D^{\zeta_0,\ep}(T)\leq \f1 4D^{\zeta_0,\ep}(T)+\f{\delta_0}4,
\]
which implies
\begin{align}\label{fina4}
 D^{\zeta_0,\ep}(T)\leq \f{\delta_0}3.
\end{align}
Then \eqref{fina1}-\eqref{fina4} implies $T^{\ast\ast}_\ep=T^\ast_\ep$. Using Proposition \ref{pro3} and \eqref{fina5}, we obtain that
\[
\|\vue\|_{L^1_{T^\ast_\ep}(\underline{B}^{\f d 2+1}_{2,1})}<\infty,\ \  \ep\|\vae\|_{L^\infty_{T^\ast_\ep}(B^{\f d 2}_{2,1})}<\infty, \inf_{(t,x)\in[0,T^\ast_\ep)\times\mathbb{R}^d}1+\ep\vae(t,x)\geq \f1 2.
\]

The  continuation criterion (see Theorem 2 in \cite{D5}) ensures that the solution $(\vae,\vue)$ can be continued beyond $T^{\ast}_\ep$. This contradicts the definition of $T^{\ast}_\ep$. Consequently, we obtain $T^{\ast}_\ep\geq T_0$, and our argument also guarantees  $T^{\ast\ast}_\ep\geq T_0$.

Next we show  vanishing limits
\begin{align}\label{fina6}
 \lim_{\ep\rightarrow0}\Bigg(\ep\|\vae\|_{\widetilde{L}^\infty_{T_0}(B^{\f d 2}_{2,1})}+\|\vVe-\vV\|_{\widetilde{L}^\infty_{T_0}(B^{\f d 2-1}_{2,1})\cap \widetilde{L}^2_{T_0}(B^{\f d 2}_{2,1})}+\|\mathbb{P}\vue-\vv\|_{\widetilde{L}^\infty_{T_0}(B^{\f d 2-1}_{2,1})\cap L^1_{T_0}(\underline{B}^{\f d 2+1}_{2,1})}\Bigg)=0.
\end{align}
For any $\delta\leq\delta_0$, we first select a large $\zeta_\delta>0$ such that  \eqref{cons3} and \eqref{cons4} remain valid
when substituting $(\delta_0, \zeta_0)$ with $(\delta, \zeta_\delta)$, and (using \eqref{energy})
\[
 \|\vV\|^{h;\zeta_\delta}_{\widetilde{L}^\infty_{T_0}(B^{\f d 2-1}_{2,1})\cap \widetilde{L}^2_{T_0}(B^{\f d 2}_{2,1})}+\|\vv\|^{h;\zeta_\delta}_{\widetilde{L}^\infty_{T_0}(B^{\f d 2-1}_{2,1})\cap L^1_{T_0}(B^{\f d 2+1}_{2,1})}\leq \f 2 3\delta.
\]

 Then, we choose a small $\ep_\delta<\min\{\ep_0,\f{\eta_0}{\zeta_\delta}\}$ such that  \eqref{cons6} and \eqref{cons5} hold with $(\delta_0, \zeta_0, \ep_0)$ replaced by $(\delta, \zeta_\delta, \ep_\delta)$. Employing an argument completely analogous to the preceding analysis, we obtain that for any $\ep\leq\ep_\delta$,
\[
D^{\zeta_\delta,\ep}(T_0)\leq \f{\delta}3.
\]
Notice that
\begin{align*}
 &\|\vVe-\vV\|_{\widetilde{L}^\infty_{T_0}(B^{\f d 2-1}_{2,1})\cap \widetilde{L}^2_{T_0}(B^{\f d 2}_{2,1})}+\|\mathbb{P}\vue-\vv\|_{\widetilde{L}^\infty_{T_0}(B^{\f d 2-1}_{2,1})\cap L^1_{T_0}(\underline{B}^{\f d 2+1}_{2,1})}\\
 &\quad\leq C^{\ast}\Big(D^{\zeta,\ep}(T_0)+ \|\vV\|^{h;\zeta}_{\widetilde{L}^\infty_{T_0}(B^{\f d 2-1}_{2,1})\cap \widetilde{L}^2_{T_0}(B^{\f d 2}_{2,1})}+\|\vv\|^{h;\zeta}_{\widetilde{L}^\infty_{T_0}(B^{\f d 2-1}_{2,1})\cap L^1_{T_0}(B^{\f d 2+1}_{2,1})}\Big),
\end{align*}
where the constant $C^{\ast}$ is independent of $\ep,\zeta$.
Hence,  for any $\ep\leq\ep_\delta$,
\begin{align*}
 &\ep\|\vae\|_{\widetilde{L}^\infty_{T_0}(B^{\f d 2}_{2,1})}+\|\vVe-\vV\|_{\widetilde{L}^\infty_{T_0}(B^{\f d 2-1}_{2,1})\cap \widetilde{L}^2_{T_0}(B^{\f d 2}_{2,1})}+\|\mathbb{P}\vue-\vv\|_{\widetilde{L}^\infty_{T_0}(B^{\f d 2-1}_{2,1})\cap L^1_{T_0}(\underline{B}^{\f d 2+1}_{2,1})}\\
 &\quad\leq (1+C^{\ast})\Big(D^{\zeta_\delta,\ep}(T_0)+ \|\vV\|^{h;\zeta_\delta}_{\widetilde{L}^\infty_{T_0}(B^{\f d 2-1}_{2,1})\cap \widetilde{L}^2_{T_0}(B^{\f d 2}_{2,1})}+\|\vv\|^{h;\zeta_\delta}_{\widetilde{L}^\infty_{T_0}(B^{\f d 2-1}_{2,1})\cap L^1_{T_0}(B^{\f d 2+1}_{2,1})}\Big)\\
 &\quad\leq (1+C^{\ast})\delta,
\end{align*}
which implies \eqref{fina6}. \ \ $\Box$
\section{The decay properties of $\vVe-\vV$ and $\mathbb{P}\vue-\vv$ }\label{Sec3}
 We recall that $(\text{Ker}L)^\perp=\{(a, \Grad g): a\in L^2(\mathbb{T}_b^d), \intTd a=0 \text{ and } g\in H^1(\mathbb{T}_b^d)\}$. Follow Danchin \cite{D1} and Masmoudi \cite{M}, we introduce the following Hilbert basis $(\Phi_k^\alpha)_{\alpha\in\{-1,1\}, k\in\widetilde{\mathbb{Z}}^d \backslash\{0\}}$ in $(\text{Ker}L)^\perp$,
\[
(\Phi_k^\alpha)(x):=c_d\left(
\begin{array}{l}
\ \ \ \ \ \ \ \ 1\\
-\alpha sg(k)k/|k|
\end{array}
\right)e^{ik\cdot x} \ \ \text{ with } c_d:=(2\T)^{-\f 1 2}.
\]
The notation $sg(k)$ stands for a generalized sign function on $\mathbb{R}^d\backslash\{0\}$: its value
is 1 if and only if the first nonzero component of $k$ is positive, $-1$ elsewhere. It
is clear that $(\Phi_k^\alpha)(x)$ is an eigenvector of $L$ for the eigenvalue $-i\alpha sg (k) |k|$.

Since $\mathcal{L}(\tau)=e^{-\tau L}$, for function $A=\sum\limits_{\alpha,k}\widehat{A}_k^{\alpha}\Phi_k^\alpha\in (\text{Ker}L)^\perp$, we have
\[
\mathcal{L}(\tau)A=\sum\limits_{\alpha,k}\widehat{A}_k^{\alpha}\Phi_k^\alpha e^{i\alpha sg (k) |k|\tau}.
\]
If $\vu=\sum\limits_{l\in\widetilde{\mathbb{Z}}^d} {{\hat \vu}_l}\,\f{e^{ il \cdot x}}{\sqrt{\T}} \in \Big(\mathcal{S}'(\mathbb{T}_b^d)\Big)^d$ satisfy that $l\cdot {\hat \vu}_l=0$ for any $l\in \widetilde{\mathbb{Z}}^d$, we get
\begin{align*}
\mathcal{Q}_1^\ep(\vu,B)=\f i {2\sqrt{\T}} &\sum\limits_{m,\gamma}\Bigg(\sum\limits_{\alpha, k+l=m}\widehat{B}_k^\alpha k\cdot \widehat{\vu}_l\left(1+\f{\alpha\gamma sg(k) sg(m)} {|k||m|}(l+m)\cdot k\right)\\
&\qquad\times e^{i\f t \ep(\alpha sg(k) |k|-\gamma sg(m)|m|)}\Bigg)\Phi_m^\gamma.
\end{align*}
We also have
\begin{align*}
\mathcal{Q}_2^\ep(A,B)=-\f{ic_d}2\sum\limits_{m,\gamma}&\Bigg[\sum\limits_{\alpha, \beta, k+l=m}sg(m)|m|\left(\f{\widehat{A}_k^\alpha \widehat{B}_l^\beta+\widehat{B}_k^\alpha \widehat{A}_l^\beta} 2\right)\\
&\times\left(\beta sg(l)sg(m)\f{l\cdot m}{|l||m|}+\f {\gamma \kappa} 2+\f {\alpha\beta\gamma}2 sg(k) sg(l) \f{k\cdot l} {|k||l|}\right)\\
&\times e^{i\f t \ep (\alpha sg(k)|k|+\beta sg(l)|l|-\gamma sg(m)|m|)}\Bigg]\Phi_m^\gamma,
\end{align*}
\[
\mathcal{A}_2^\ep(D)B=-\f 1 2\sum\limits_{\alpha,\gamma,k}\alpha\gamma|k|^2 \widehat{B}_k^\alpha e^{i\f t \ep(\alpha-\gamma)sg(k)|k|}\Phi_k^\gamma.
\]
Danchin \cite{D1} calculated the limits of $\mathcal{Q}_1^\ep(\vu,B)$, $\mathcal{Q}_2^\ep(A,B)$ and $\mathcal{A}_2^\ep(D)B$, and his calculations were based on the Riemann-Lebesgue Lemma, i.e.
\[
\lim_{\ep\rightarrow 0}\int_0^T \chi(t)e^{i\f {at} \ep}\dt=0,  \text{ if } 0<T\leq\infty, a\in \mathbb{R}\backslash\{0\}, \chi\in L^1(0,T).
\]
He showed
\[
 \lim_{\ep\rightarrow 0}\mathcal{Q}_1^\ep(\vu,B)=\mathcal{Q}_1(\vu,B):=\f {i}{\sqrt{\T}}\sum\limits_{\gamma,m}\left(\sum\limits_{\begin{array}{c}
                                                                                           k+l=m\\
                                                                                            \alpha sg(k)|k|=\gamma sg(m)|m|
                                                                                         \end{array}
}\widehat{B}_k^\alpha \f{k\cdot \widehat{\vu}_l k\cdot m}{|k||m|}\right)\Phi_m^\gamma,
\]
\begin{align*}
\lim_{\ep\rightarrow 0}\mathcal{Q}_2^\ep(A,B)&=\mathcal{Q}_2(A,B)\\
&:=-ic_d\f {\kappa+3} 4 \sum\limits_{\gamma,m}\gamma sg(m)|m|\left(\sum\limits_{\begin{array}{c}
  k+l=m\\
  sg(k)|k|+sg(l)|l|=\gamma sg(m)|m|
  \end{array}}
  \widehat{A}_k^\gamma \widehat{B}_l^\gamma
                                                                                         \right)\Phi_m^\gamma,
\end{align*}

\[
\lim_{\ep\rightarrow 0^{+}}\mathcal{A}_2^\ep(D)B= \f1 2 \Delta B.
\]
Finally, we note that the operator $\mathcal{L}(\tau)$ is  an isomorphism on $H^{s}$ or $B^{s}_{2,1}$. Moreover, we have:
\begin{Lemma}\emph{(\cite{D1})}\label{LeD3}. For all $\alpha,s\in \mathbb{R}$ and $q\in[1,\infty]$, the operator
\begin{equation*}
   \left\{\bali
  \mathcal{S}'((0,T)\times \mathbb{T}^{d+1})&\rightarrow  \mathcal{S}'((0,T)\times \mathbb{T}^{d+1})\\
  &u\mapsto\mathcal{L}(\alpha t)u
  \eali
  \right.
\end{equation*}
is an isomorphism on $L^q(0,T;H^{s})$, $\widetilde{L}^q(0,T;H^{s})$, $L^q(0,T;B^{s}_{2,1})$ and $\widetilde{L}^q(0,T;B^{s}_{2,1})$.
\end{Lemma}
From Lemma \ref{LeD3} and the definition of $Y^{\zeta,\ep}(T)$, we infer that for $q\in[2,\infty]$,
\begin{align}\label{estr1}
 \|\vae\|_{ \widetilde{L}^\infty_T(B^{\f d 2-1}_{2,1}) }&+\|(\vae,\mathbb{Q}\vue)\|_{ \widetilde{L}^2_T(B^{\f d 2}_{2,1}) }\nonumber\\
 &+\|\mathbb{P}\vue\|_{\widetilde{L}_T^q(\underline{B}^{\f d 2-1+\f 2 q}_{2,1})}+\|\vVe\|_{\widetilde{L}_T^q(B^{\f d 2-1+\f 2 q}_{2,1})} \lesssim  Y^{\zeta,\ep}(T).
\end{align}
\subsection{The proof of Proposition \ref{pro1}} \label{newsec1}
Let $\vZe=\vVe-\vV, \vWe=\mathbb{P}\vue-\vv$. Subtracting \eqref{equ4} from \eqref{equ5}, we find
\begin{align}\label{equ6}
 \p_t \vZe+\mathcal{Q}_1^\ep(\vv,\vZe)+2\mathcal{Q}_2^\ep(\vV,\vZe)+&\mathcal{Q}_2^\ep(\vZe,\vZe)-\f 1 2 \De\vZe=F^\ep+R^{1,\ep}+R^{2,\ep}+R^{3,\ep}+S^\ep
\end{align}
with
\begin{align*}
 F^\ep=&-\mathcal{L}\Big(-\f t {\ep}\Big)\left(
\begin{array}{l}
\qquad \qquad \qquad \qquad  \qquad \qquad 0 \\
\mathbb{Q}(\vv\cdot\Grad \vWe+\vWe\cdot\mathbb{P}\vue+I(\ep\vae)\mathcal{A}\vue+\vae \widetilde{K}(\ep\vae)\Grad \vae)
\end{array}
\right)\\
&-\mathcal{Q}_1^\ep(\underline{\vWe},\vVe)-\mathcal{Q}_1^\ep(\f{\widehat{\vWe}_0}{\sqrt{\Td}},\vVe)
\end{align*}
\[
R^{1,\ep}=\mathcal{L}\Big(-\f t {\ep}\Big)\left(
\begin{array}{l}
\qquad  \quad 0\\
\mathbb{Q}f-\mathbb{Q}(\vv\cdot\Grad \vv)
\end{array}
\right),
\]
\[
R^{2,\ep}=(\mathcal{Q}_1-\mathcal{Q}^\ep_1)(\vv,\vV),\ \  R^{3,\ep}=(\mathcal{Q}_2-\mathcal{Q}^\ep_2)(\vV,\vV),\ \ S^\ep=\nu(\mathcal{A}^\ep_2(D)- {\De}/2)\vVe.
\]
We decompose $(0,\mathbb{Q}f)^T$ and  $(0,\mathbb{Q}(\vv\cdot\Grad \vv))^T$ on the basis $(\Phi_k^\alpha)_{\alpha\in\{-1,1\}, k\in\widetilde{\mathbb{Z}}^d \backslash\{0\}}$,
\[
(0,\mathbb{Q}f)^T=\sum\limits_{k,\alpha}\widehat{f}_k^{\alpha}\Phi_k^\alpha,\ \ (0,\mathbb{Q}(\vv\cdot\Grad \vv))^T=\sum\limits_{k,\alpha}\widehat{\Lambda}_k^{\alpha}\Phi_k^\alpha
\]
then
\[
R^{1,\ep}=\sum\limits_{k,\alpha}(\widehat{f}_k^{\alpha}-\widehat{\Lambda}_k^{\alpha})e^{-i\f t \ep \alpha sg(k)|k|}\Phi_k^\alpha.
\]
We also have
\[
S^{\ep}=\f{\nu} 2\sum\limits_{k,\alpha} |k|^2 \widehat{\vV}^{\alpha,\ep}_k e^{2\f {i\alpha t}\ep sg(k)|k|}\Phi_k^{-\alpha}.
\]
For $\gamma\in\{-1,1\}$ and $m\in \widetilde{\mathbb{Z}}^d\backslash{\{0\}}$, we introduce
\[
\mathcal{B}_m^{1,\gamma}=\{(\alpha,k,l)\in\{-1,1\}\times\widetilde{\mathbb{Z}}^d\times\widetilde{\mathbb{Z}}^d:\alpha sg(k)|k|\neq\gamma sg(m)|m| \text{ and } k+l=m\},
\]
\[
\mathcal{B}_m^{2,\gamma}=\{(\alpha,\beta,k,l)\in\{-1,1\}^2\times\widetilde{\mathbb{Z}}^d\times\widetilde{\mathbb{Z}}^d:\alpha sg(k)|k|+\beta sg(l)|l|\neq\gamma sg(m)|m| \text{ and } k+l=m\},
\]
then
\begin{align*}
R^{2,\ep}=-\f i {2\sqrt{\T}} &\sum\limits_{m,\gamma,(\alpha,k,l)\in \mathcal{B}_m^{1,\gamma}}\widehat{\vV}_k^\alpha k\cdot \widehat{\vv}_l\left(1+\f{\alpha\gamma sg(k) sg(m)} {|k||m|}(l+m)\cdot k\right)\\
&\qquad\times e^{i\f t \ep(\alpha sg(k) |k|-\gamma sg(m)|m|)}\Phi_m^\gamma,
\end{align*}
\begin{align*}
R^{3,\ep}=\f{ic_d}2\sum\limits_{m,\gamma,(\alpha,k,l)\in \mathcal{B}_m^{2,\gamma}}&sg(m)|m|\widehat{\vV}_k^\alpha \widehat{\vV}_l^\beta \\
&\times\left(\beta sg(l)sg(m)\f{l\cdot m}{|l||m|}+\f {\gamma \kappa} 2+\f {\alpha\beta\gamma}2 sg(k) sg(l) \f{k\cdot l} {|k||l|}\right)\\
&\times e^{i\f t \ep (\alpha sg(k)|k|+\beta sg(l)|l|-\gamma sg(m)|m|)}\Phi_m^\gamma.
\end{align*}

We follow Danchin \cite{D1} and Schochet \cite{SS2} and introduce the function $\phi_{M}^\ep:=\vZe-\ep \widetilde{R}_M^\ep$ with $\widetilde{R}_M^\ep:=\widetilde{R}_M^{1,\ep}+\widetilde{R}_M^{2,\ep}+\widetilde{R}_M^{3,\ep}+\widetilde{S}_M^\ep$, $M>0$ and
\[
\widetilde{S}^{\ep}_M=-\f{i} 4\nu\sum\limits_{\alpha,|k|\leq M} \alpha sg(k)|k| \widehat{\vV}^{\alpha,\ep}_k e^{2\f {i\alpha t}\ep sg(k)|k|}\Phi_k^{-\alpha},
\]
\begin{align*}
 \widetilde{R}^{1,\ep}_M=i\sum\limits_{\alpha,|k|\leq M}\alpha sg(k)|k|^{-1}(\widehat{f}_k^{\alpha}-\widehat{\Lambda}_k^{\alpha})e^{-i\f t \ep \alpha sg(k)|k|}\Phi_k^\alpha
\end{align*}
\begin{align*}
\widetilde{R}^{2,\ep}_M=-\f 1 {2\sqrt{\T}} &\sum\limits_{m,\gamma,(\alpha,k,l)\in \mathcal{B}_{m,M}^{1,\gamma}}\f{\widehat{\vV}_k^\alpha k\cdot \widehat{\vv}_l}{\alpha sg(k) |k|-\gamma sg(m)|m|}\left(1+\f{\alpha\gamma sg(k) sg(m)} {|k||m|}(l+m)\cdot k\right)\\
&\qquad\times e^{i\f t \ep(\alpha sg(k) |k|-\gamma sg(m)|m|)}\Phi_m^\gamma,
\end{align*}
\begin{align*}
\widetilde{R}^{3,\ep}_M=\f{c_d}2\sum\limits_{m,\gamma,(\alpha,k,l)\in \mathcal{B}_{m,M}^{2,\gamma}}&\f{\widehat{\vV}_k^\alpha \widehat{\vV}_l^\beta}{\alpha sg(k)|k|+\beta sg(l)|l|-\gamma sg(m)|m|} sg(m)|m|\\
&\times\left(\beta sg(l)sg(m)\f{l\cdot m}{|l||m|}+\f {\gamma \kappa} 2+\f {\alpha\beta\gamma}2 sg(k) sg(l) \f{k\cdot l} {|k||l|}\right)\\
&\times e^{i\f t \ep (\alpha sg(k)|k|+\beta sg(l)|l|-\gamma sg(m)|m|)}\Phi_m^\gamma.
\end{align*}
Here,  $\mathcal{B}_{m,M}^{1,\gamma}:=\{(\alpha,k,l)\in \mathcal{B}_{m}^{1,\gamma}: |k|\leq M, |l|\leq M\}$ and $\mathcal{B}_{m,M}^{2,\gamma}:=\{(\alpha,\beta,k,l)\in \mathcal{B}_{m}^{2,\gamma}: |k|\leq M, |l|\leq M\}$.

We denote by $R_M^\ep$ \textquotedblleft the low frequency part \textquotedblright of $R^\ep:=R^{1,\ep}+R^{2,\ep}+R^{3,\ep}+S^\ep$ obtained by keeping only the indices $(k,l)$ such that $|k|,|l|\leq M$ in the summations, and  $R^{\ep,M}:=R^\ep-R_M^\ep$. It is clear that
\begin{align}\label{equ7}
  \ep\p_t\widetilde{R}_M^\ep=R^\ep_M+\ep \widetilde{R}^{t,\ep}_M
\end{align}
with $\widetilde{R}^{t,\ep}_M:=\widetilde{R}^{1,t,\ep}_M+\widetilde{R}^{2,t,\ep}_M+\widetilde{R}^{3,t,\ep}_M+\widetilde{S}^{t,\ep}_M$ and
\[
\widetilde{S}^{t,\ep}_M=-\f{i} 4\nu\sum\limits_{\alpha,|k|\leq M} \alpha sg(k)|k| \p_t\widehat{\vV}^{\alpha,\ep}_k e^{2\f {i\alpha t}\ep sg(k)|k|}\Phi_k^{-\alpha},
\]
\begin{align*}
 \widetilde{R}^{1,t, \ep}_M=i\sum\limits_{\alpha,|k|\leq M}\alpha sg(k)|k|^{-1}(\p_t\widehat{f}_k^{\alpha}-\p_t\widehat{\Lambda}_k^{\alpha})e^{-i\f t \ep \alpha sg(k)|k|}\Phi_k^\alpha
\end{align*}
\begin{align*}
\widetilde{R}^{2,t,\ep}=-\f 1 {2\sqrt{\T}} &\sum\limits_{m,\gamma,(\alpha,k,l)\in \mathcal{B}_{m,M}^{1,\gamma}}\f{\p_t(\widehat{\vV}_k^\alpha k\cdot \widehat{\vv}_l)}{\alpha sg(k) |k|-\gamma sg(m)|m|}\left(1+\f{\alpha\gamma sg(k) sg(m)} {|k||m|}(l+m)\cdot k\right)\\
&\qquad\times e^{i\f t \ep(\alpha sg(k) |k|-\gamma sg(m)|m|)}\Phi_m^\gamma,
\end{align*}
\begin{align*}
\widetilde{R}^{3,t,\ep}=\f{c_d}2\sum\limits_{m,\gamma,(\alpha,k,l)\in \mathcal{B}_{m,M}^{2,\gamma}}&\f{\p_t(\widehat{\vV}_k^\alpha \widehat{\vV}_l^\beta)}{\alpha sg(k)|k|+\beta sg(l)|l|-\gamma sg(m)|m|} sg(m)|m|\\
&\times\left(\beta sg(l)sg(m)\f{l\cdot m}{|l||m|}+\f {\gamma \kappa} 2+\f {\alpha\beta\gamma}2 sg(k) sg(l) \f{k\cdot l} {|k||l|}\right)\\
&\times e^{i\f t \ep (\alpha sg(k)|k|+\beta sg(l)|l|-\gamma sg(m)|m|)}\Phi_m^\gamma.
\end{align*}
From \eqref{equ6} and \eqref{equ7}, we have that
\begin{align}\label{equ8}
  &\p_t\phi_M^\ep+ \widecheck{\mathcal{Q}}_1^\ep(\underline{\vv},\phi_M^\ep)+\widetilde{\mathcal{Q}}_1^\ep(\underline{\vv},\phi_M^\ep)+\mathcal{Q}_1^\ep(\f{\widehat{\vv}_0}{\sqrt{\Td}},\phi_M^\ep)+2\mathcal{Q}_2^\ep(\vV,\phi_M^\ep)-\f{\nu}2\De\phi_M^\ep \nonumber\\
  &=R^{\ep,M}-\ep\left(\widecheck{\mathcal{Q}}_1^\ep(\vv,\widetilde{R}_M^\ep)+\widetilde{\mathcal{Q}}_1^\ep(\vv,\widetilde{R}_M^\ep)+2\mathcal{Q}_2^\ep(\vV,\widetilde{R}_M^\ep)-\f \nu 2\De \widetilde{R}_M^\ep+\widetilde{R}_M^{t,\ep}\right)\nonumber\\
  &\quad-\mathcal{Q}_2^\ep(z^\ep,z^\ep)+F^\ep,
\end{align}
where
\[
\widecheck{\mathcal{Q}}_1^\ep(\vu,B):=\mathcal{L}\Big(-\f t {\ep}\Big)\left(
\begin{array}{l}
\quad\Div\left(\vu\mathcal{L}^1\Big(\f t {\ep}\Big)B\right)\\
\mathbb{Q}\Div\left(\vu\otimes\mathcal{L}^2\Big(\f t {\ep}\Big)B\right)
\end{array}
\right),\ \  \widetilde{\mathcal{Q}}_1^\ep(\vu,B):= \mathcal{Q}_1^\ep(\vu,B)-\widecheck{\mathcal{Q}}_1^\ep(\vu,B).
\]
Since $\vv\in \left(F_{T_0}^{\f d 2}(\mathbb{T}_b^d)\right)^d$, then by  interpolation, we infer that
\begin{align}\label{inter}
 \vv\in \widetilde{L}_{T_0}^2(\underline{B}^{\f d 2}_{2,1}(\mathbb{R}^d))\ \ (rather\ \ than\ \  \vv\in \widetilde{L}_{T_0}^2(B^{\f d 2}_{2,1}(\mathbb{R}^d))),
\end{align}
where
$\widetilde{L}_{T_0}^2(\underline{B}^{\f d 2}_{2,1}(\mathbb{R}^d)):=\{f: \|f\|_{\widetilde{L}^{2}_T(\underline{B}_{2,1}^{\f d 2})}<\infty\}$ and $\widetilde{L}_{T_0}^2(B^{\f d 2}_{2,1}(\mathbb{R}^d)):=\{f: \|f\|_{\widetilde{L}^{2}_T(B_{2,1}^{\f d 2})}<\infty\}$.

Due to \eqref{inter}, we decompose $\widecheck{\mathcal{Q}}_1^\ep(\vv,\phi_M^\ep)+\widetilde{\mathcal{Q}}_1^\ep(\vv,\phi_M^\ep)$ into three components in \eqref{equ8} as follows:
\begin{align*}
 \widecheck{\mathcal{Q}}_1^\ep(\vv,\phi_M^\ep)+\widetilde{\mathcal{Q}}_1^\ep(\vv,\phi_M^\ep)
&=\widecheck{\mathcal{Q}}_1^\ep(\underline{\vv},\phi_M^\ep)+\widetilde{\mathcal{Q}}_1^\ep(\underline{\vv},\phi_M^\ep)
+\widecheck{\mathcal{Q}}_1^\ep(\f{\widehat{\vv}_0}{\sqrt{\Td}},\phi_M^\ep)+\widetilde{\mathcal{Q}}_1^\ep(\f{\widehat{\vv}_0}{\sqrt{\Td}},\phi_M^\ep)\\
&=\widecheck{\mathcal{Q}}_1^\ep(\underline{\vv},\phi_M^\ep)+\widetilde{\mathcal{Q}}_1^\ep(\underline{\vv},\phi_M^\ep)+
\mathcal{Q}_1^\ep(\f{\widehat{\vv}_0}{\sqrt{\Td}},\phi_M^\ep).
\end{align*}
Through straightforward yet tedious calculations, we find that\footnote{For a function $f$ defined on $\Td$, $f$ is real-valued if and only if $\widehat{f}_k=\overline{\widehat{f}_{-k}}$ for any $k \in \widetilde{\mathbb{Z}}^d$. } for any $\ep, M>0$,
\begin{align}\label{def1}
  \text{ $\phi_M^\ep$ is real-valued } \text{ and } \mathcal{Q}_1^\ep(\f{\widehat{\vv}_0}{\sqrt{\Td}},\phi_M^\ep)=\mathcal{Q}_1(\f{\widehat{\vv}_0}{\sqrt{\Td}},\phi_M^\ep).
\end{align}
We recall the vanishing integral property established in \cite{D1}, which states that if $\vV$ is a real-valued function and $\Div \vu=0$, then
\begin{align}\label{eli}
\int_{\Td}\De_j\mathcal{Q}_1(\vu,\vV)\cdot\De_j\vV dx=0 \text{ for any } j\in\mathbb{Z}.
\end{align}

By employing the proof technique of  Proposition 4.1 in \cite{D1}, in combination with
 \eqref{def1}-\eqref{eli} and  \eqref{estQ1}-\eqref{estQ2},  we establish following estimates,
\begin{align}\label{est1}
 &\|\phi_M^\ep\|_{\widetilde{L}_T^\infty(H^{\f d 2-\vartheta-1})}+ \|\phi_M^\ep\|_{L_T^2(H^{\f d 2-\vartheta})}\nonumber\\
 &\lesssim e^{C(\|\underline{\vv}\|^2_{L^2_T(B^{\f d 2}_{2,1})}+\|\vV\|^2_{L^2_T(B^{\f d 2}_{2,1})})}\Bigg(\|R^{\ep,M}\|_{\widetilde{L}^1_T(H^{\f d 2 -\vartheta-1})+L^2_T(H^{\f d 2 -\vartheta})}+\|\mathcal{Q}_2^\ep(\vZe,\vZe)\|_{L^2_T(H^{\f d 2 -\vartheta-2})}\nonumber\\
 &\quad+\|F^\ep\|_{\widetilde{L}^1_T(H^{\f d 2 -\vartheta-1})+L^2_T(H^{\f d 2 -\vartheta})}+\ep\Bigg(\|\widetilde{R}_M^\ep(0)\|_{H^{\f d 2 -\vartheta-1}}+\|\widetilde{\mathcal{Q}}_1^\ep(\vv,\tilde{R}_M^\ep)\|_{\widetilde{L}^1_T(H^{\f d 2 -\vartheta-1})}\nonumber\\
 &\quad+\|\widetilde{R}_M^{t,\ep}\|_{\widetilde{L}^1_T(H^{\f d 2 -\vartheta-1})\cap L_T^2(H^{\f d 2-\vartheta})}+\|\widecheck{\mathcal{Q}}_1^\ep(\vv,\widetilde{R}_M^\ep)\|_{L^2_T(H^{\f d 2 -\vartheta-2})}\nonumber\\
 &\quad+\|\mathcal{Q}_2^\ep(\vV,\widetilde{R}_M^\ep\|_{L^2_T(H^{\f d 2 -\vartheta-2})}+\|\De\widetilde{R}_M^\ep\|_{L^2_T(H^{\f d 2 -\vartheta-2})}\Bigg)\Bigg)
\end{align}

Now we prove that the discrepancy between $\phi_M^\ep$ and $\vZe$ is small. Hence we need to estimate $\|\widetilde{R}_M^\ep\|_{\widetilde{L}_T^\infty(H^{\f d 2 -\vartheta-1})\cap L_T^2(H^{\f d 2 -\vartheta})}$. We claim that
\begin{align}\label{estR1}
 \|\widetilde{R}_M^\ep\|_{\widetilde{L}_T^\infty(H^{\f d 2 -\vartheta-1})\cap L_T^2(H^{\f d 2 -\vartheta})}\leq  C_M(Q_f+Y^{\zeta,\ep}(T)+X^2+P^2),
\end{align}
\begin{align}\label{estR2}
 \|\widetilde{R}_M^\ep(0)\|_{H^{\f d 2 -\vartheta-1}}\leq  C_M(Q_f(0)+X_0+X_0^2+P_0^2),
\end{align}
where $C_M$ denotes a  monotonically increasing positive function of $M$ to be precise hereafter.

First, we present some truncated estimates. For any time-dependent function $\Psi$ taking values in $\mathcal{S}'(\Td)$, any $g\in \mathcal{S}'(\Td)$, $q\in [1,\infty]$, $s\in \mathbb{R}$, $r\in\{1,2\}$ and $M\geq1,\sigma>0$, we have
\begin{align}
  &\|\Psi_M\|_{\widetilde{L}_T^q(H^s)}\lesssim M^{\sigma} \|\Psi\|_{\widetilde{L}_T^q(H^{s-\sigma})},\ \  \|\Psi^M\|_{\widetilde{L}_T^q(H^s)}\lesssim M^{-\sigma} \|\Psi\|_{\widetilde{L}_T^q(H^{s+\sigma})},\label{estM1}\\
   &\|g_M\|_{B^s_{2,r}}\lesssim M^{\sigma} \|g\|_{B^{s-\sigma}_{2,r}},\qquad\qquad  \|g^M\|_{B^s_{2,r}}\lesssim M^{-\sigma} \|g\|_{B^{s+\sigma}_{2,r}},\label{estM3}\\
& \|\mathcal{F}g_M\|_{\textit{l}^1(\widetilde{\mathbb{Z}}^d)} \lesssim M^{\sigma} \|g\|_{H^{\f d 2-\sigma}},\qquad\quad \|\mathcal{F}g^M\|_{\textit{l}^1(\widetilde{\mathbb{Z}}^d)} \lesssim M^{-\sigma} \|g\|_{H^{\f d 2+\sigma}}\label{estM2},
\end{align}
where
\[
f_M:=\sum\limits_{|k|\leq M} {{\hat f}_k}\,\f{e^{ ik \cdot x}}{\sqrt{\T}},\ \ f^M:=\sum\limits_{|k|> M} {{\hat f}_k}\,\f{e^{ ik \cdot x}}{\sqrt{\T}}.
\]
Next, let us estimate $\widetilde{S}_M^\ep, \widetilde{R}_M^{1,\ep}, \widetilde{R}_M^{2,\ep}, \widetilde{R}_M^{3,\ep}$ term by term.
From \eqref{estM1}-\eqref{estM3},
\begin{align}\label{estR3}
 \|\widetilde{S}_M^\ep\|_{\widetilde{L}_T^\infty(H^{\f d 2 -\vartheta-1})\cap L_T^2(H^{\f d 2 -\vartheta})}\lesssim M^{1-\vartheta}Y^{\zeta,\ep}(T),
\end{align}
\begin{align}\label{estR4}
 \|\widetilde{S}_M^\ep(0)\|_{H^{\f d 2 -\vartheta-1}}\lesssim M^{1-\vartheta}X_0.
\end{align}
Recall that $\mathbb{Q}f\in C_{T_0}(H^{-S})\cap L^1_{T_0}(B^{\f d 2-1}_{2,1})$, by interpolation, $\mathbb{Q}f\in L^2_{T_0}(H^{\f {\f d 2-1-S}{2}})$. From \eqref{minski2} and \eqref{estM1}, we have
\begin{align}\label{estR5}
  \| \widetilde{R}_M^{1,\ep}\|_{L_T^2(H^{\f d 2 -\vartheta})}&\lesssim M^{\f{\f d 2+S-2\vartheta+1}2}\|\mathbb{Q}f\|_{L^2_T(H^{\f {\f d 2-1-S}{2}-1})}+M^{1-\vartheta}\|\vv\otimes\vv\|_{L^2_T(\underline{H}^{\f d 2 -1})}\nonumber\\
&\lesssim  M^{\f{\f d 2+S-2\vartheta+1}2}\mathbb{Q}f+M^{1-\vartheta}P^2,
\end{align}
\begin{align}\label{estR6}
  \| \widetilde{R}_M^{1,\ep}\|_{\widetilde{L}_T^\infty(H^{\f d 2 -\vartheta-1})}&\lesssim M^{\f d 2+S-1} \|\mathbb{Q}f\|_{\widetilde{L}_T^\infty(H^{-S-\vartheta-1})}+M\|\vv\otimes\vv\|_{L^\infty_T(\underline{B}^{\f d 2 -2}_{2,\infty})}\nonumber\\
  &\lesssim M^{\f d 2+S-1}\mathbb{Q}f+MP^2,
\end{align}
\begin{align}\label{estR7}
  \| \widetilde{R}_M^{1,\ep}(0)\|_{H^{\f d 2 -\vartheta-1}}\lesssim M^{\f d 2+S-1}\mathbb{Q}f(0)+MP_0^2(T),
\end{align}
where the spatial norm $\underline{B}^{\f d 2 -2}_{2,\infty}$ in \eqref{estR6} is specifically introduced to provide a unified treatment for both cases of $d=2$ and $d>2$.

Next we turn to the estimates of $\widetilde{R}_M^{2,\ep}, \widetilde{R}_M^{3,\ep}$ and introduce the notation
\begin{align}\label{nota1}
C_M^1:=\max_{\begin{array}{c}
               |m|\leq 2M, \gamma\in\{-1,1\} \\
               (\alpha,k,l)\in\mathcal{B}^{1,\gamma}_{m,M}
             \end{array}
}|\alpha sg(k)|k|-\gamma sg(m)|m||^{-1},
\end{align}
\begin{align}\label{nota2}
C_M^2:=\max_{\begin{array}{c}
               |m|\leq 2M, \gamma\in\{-1,1\} \\
               (\alpha,\beta,k,l)\in\mathcal{B}^{2,\gamma}_{m,M}
             \end{array}
}|\alpha sg(k)|k|+\beta sg(l)|l|-\gamma sg(m)|m||^{-1}.
\end{align}
 Notice that $\Div \vv=0$, we have $l\cdot \hat{\vv}_l=0$ for all $l\in \widetilde{\mathbb{Z}}^d$. Hence $\widetilde{R}_M^{2,\ep}$ rewrites
\begin{align*}
\widetilde{R}^{2,\ep}_M=-\f 1 {2\sqrt{\T}} &\sum\limits_{m,\gamma,(\alpha,k,l)\in \mathcal{B}_{m,M}^{1,\gamma}}\f{\alpha\gamma sg(k) sg(m)\widehat{\vV}_k^\alpha}{\alpha sg(k) |k|-\gamma sg(m)|m|}\\
&\qquad \times \Bigg[\f{m\cdot k}{|k||m|}k\cdot \hat{\vv}_l+m\cdot\hat{\vv}_l\left(\alpha\gamma sg(k) sg(m)+
\f{l\cdot k}{|k||m|}\right)\Bigg]\\
&\qquad\times e^{i\f t \ep(\alpha sg(k) |k|-\gamma sg(m)|m|)}\Phi_m^\gamma.
\end{align*}
From \eqref{minski2}, \eqref{estM3}-\eqref{estM2} and Lemma \ref{leD2},
\[
\|\widetilde{R}^{2,\ep}_M\|_{L^2_T(H^{\f d 2 -\vartheta})}\lesssim \|\widetilde{R}^{2,\ep}_M\|_{L^2_T(H^{\f d 2})} \lesssim C_M^1 M^2 \|\vV\|_{L^2_T(H^{\f d 2})\cap L^\infty_T(H^{\f d 2-1})}\|\vv\|_{L^2_T(\underline{H}^{\f d 2})\cap L^\infty_T(H^{\f d 2-1})},
\]
\[
\|\widetilde{R}^{2,\ep}_M\|_{\widetilde{L}^\infty_T(H^{\f d 2 -1-\vartheta})}\lesssim \|\widetilde{R}^{2,\ep}_M\|_{L^\infty_T(H^{\f d 2-1})} \lesssim C_M^1 M^{2+2\vartheta} \|\vV\|_{L^\infty_T(H^{\f d 2-1-\vartheta})}\|\vv\|_{L^\infty_T(H^{\f d 2-1-\vartheta})},
\]
\begin{align}\label{estR10}
\|\widetilde{R}^{2,\ep}_M(0)\|_{H^{\f d 2 -1-\vartheta}}\lesssim C_M^1 M^{2+2\vartheta} X_0P_0.
\end{align}
It follows that
\begin{align}\label{estR8}
\|\widetilde{R}^{2,\ep}_M\|_{\widetilde{L}^\infty_T(H^{\f d 2 -1-\vartheta})\cap L^2_T(H^{\f d 2 -\vartheta})}\lesssim C_M^1 M^{2+2\vartheta} XP.
\end{align}
A similar argument gives that
\begin{align}\label{estR9}
\|\widetilde{R}^{3,\ep}_M\|_{\widetilde{L}^\infty_T(H^{\f d 2 -1-\vartheta})\cap L^2_T(H^{\f d 2 -\vartheta})}\lesssim C_M^1 M^{2+2\vartheta}X^2,
\end{align}
\begin{align}\label{estR11}
\|\widetilde{R}^{3,\ep}_M(0)\|_{H^{\f d 2 -1-\vartheta}}\lesssim C_M^1 M^{2+2\vartheta} X^2_0.
\end{align}
Defining
\begin{align}\label{defC}
 C_M:= C \max \{M, M^{\f{\f d 2+S-2\vartheta+1}2}, M^{\f d 2+S-1}, (C_M^1+C_M^2)M^{2+2\vartheta}\}.
\end{align}
Then \eqref{estR1} and \eqref{estR2}  are derived from
\eqref{estR3}-\eqref{estR7} and  \eqref{estR10}-\eqref{estR11}.  It follows from \eqref{estR1} and \eqref{estR2} that
\begin{align}
 \|\phi_M^\ep-\vZe\|_{\widetilde{L}_T^\infty(H^{\f d 2-\vartheta-1})\cap L_T^2(H^{\f d 2-\vartheta})}&\leq \ep C_M\left(\mathbb{Q}_f+Y^{\zeta,\ep}(T)+X^2+P^2\right)\label{est2},\\
 \|\phi_M^\ep(0)\|_{H^{\f d 2-\vartheta-1}}&\leq \ep C_M\left(\mathbb{Q}_f(0)+X_0+X_0^2+P_0^2\right).
\end{align}
It's time to estimate the right-hand side of \eqref{est1}. From \eqref{estR1}, we have that
\begin{align}\label{estR35}
  \|\De\widetilde{R}_M^\ep\|_{L^2_T(H^{\f d 2 -\vartheta-2})}\lesssim  C_M\left(\mathbb{Q}_f+Y^{\zeta,\ep}(T)+X^2+P^2\right).
\end{align}
From Lemma \ref{LeD3} and \eqref{estQ2}, we infer that
\begin{align}\label{estZ}
 \|\mathcal{Q}_2^\ep(\vZe,\vZe)\|_{L^2_T(H^{\f d 2 -\vartheta-2})}&\lesssim \|\vZe\|_{L^2_T(B^{\f d 2}_{2,1})}  \|\vZe\|_{L^\infty_T(H^{\f d 2 -\vartheta-1})}\nonumber\\
 &\lesssim \left(D^{\zeta,\ep}(T)+\|\vV\|^{h;\zeta}_{\widetilde{L}^2_{T_0}(B^{\f d 2}_{2,1})} \right) \vZe_\vartheta(T).
\end{align}

\leftline{\emph{Estimate for $R^{\ep,M}$.}}From \eqref{estM1}-\eqref{estM3}, we have
\begin{align}\label{estR12}
 \|S^{\ep,M}\|_{L^2_T(H^{\f d 2 -\vartheta-2})} \lesssim  M^{-\vartheta}  \|\vVe\|_{L^2_T(H^{\f d 2})} \lesssim  M^{-\vartheta} Y^{\zeta,\ep}(T).
\end{align}
\begin{align}\label{estR13}
 \|R^{1,\ep,M}\|_{\widetilde{L}^1_T(H^{\f d 2 -\vartheta-1})} &\lesssim  M^{-\vartheta}  \left(\|\mathbb{Q}f\|_{\widetilde{L}^1_T(H^{\f d 2 -1})}+\|\vv\cdot \Grad \vv\|_{\widetilde{L}^1_T(H^{\f d 2 -1})}\right)\nonumber\\
 &\lesssim  M^{-\vartheta}\left(\mathbb{Q}_f+P^2\right).
\end{align}
Notice that
\[
R^{2,\ep,M}=(\mathcal{Q}_1-\mathcal{Q}_1^\ep)(\vv,\vV^M)+(\mathcal{Q}_1-\mathcal{Q}_1^\ep)(\vv^M,\vV_M).
\]
From the definition of $\mathcal{Q}_1^\ep$ and \eqref{estM2},
\begin{align}\label{estR15}
\|\mathcal{Q}_1^\ep(\vv,\vV^M)\|_{\widetilde{L}^1_T(H^{\f d 2 -\vartheta-1})} &\lesssim \|(\vv \Grad\mathcal{L}^1(\f t \ep) \vV^M, \vv\cdot \Grad\mathcal{L}^2(\f t \ep) \vV^M , \mathcal{L}^2(\f t \ep) \vV^M \Grad \vv)\|_{L^1_T(\underline{H}^{\f d 2 -\vartheta-1})}\nonumber\\
&\lesssim  \|\vv\|_{L^2_T(\underline{B}^{\f d 2}_{2,1})\cap L^\infty_T(B^{\f d 2-1}_{2,1})} \|\vV^M\|_{L^2_T(H^{\f d 2-\vartheta})\cap L^1_T(H^{\f d 2+1-\vartheta})}\nonumber\\
&\lesssim M^{-\vartheta}XP .
\end{align}
Similarly,
\begin{align}\label{estR16}
\|\mathcal{Q}_1^\ep(\vv^M,\vV_M)\|_{\widetilde{L}^1_T(H^{\f d 2 -\vartheta-1})} \lesssim M^{-\vartheta}XP.
\end{align}
Using \eqref{estM3}-\eqref{estM2} and \eqref{estQ4},
\begin{align}
  \|\mathcal{Q}_1(\vv,V^M)\|_{H^{\f d 2-\vartheta-1}}&\lesssim \|\vV^M\|_{H^{\f d 2-\f\vartheta 2}}\||k|^{-\f \vartheta 2}\mathcal{F}\underline{\vv}\|_{\textit{l}^1(\widetilde{\mathbb{Z}}^d)}+|\widehat{\vv}_0|\|\vV^M\|_{H^{\f d 2-\vartheta}}\nonumber\\
  &\lesssim M^{-\f \vartheta 2}\|\vV\|_{H^{\f d 2}}\|\underline{\vv}\|_{H^{\f d 2}}+M^{-\vartheta }\|\vV\|_{H^{\f d 2+1}}\|\vv\|_{H^{\f d 2-1}},\label{estR17}\\
  \|\mathcal{Q}_1(\vv^M,V_M)\|_{H^{\f d 2-\vartheta-1}}&\lesssim \|\vV_M\|_{H^{\f d 2}}\||k|^{- \vartheta }\mathcal{F}\vv^M\|_{\textit{l}^1(\widetilde{\mathbb{Z}}^d)}\lesssim M^{- \vartheta}\|\vV\|_{H^{\f d 2}}\|\underline{\vv}\|_{H^{\f d 2}}.\label{estR18}
\end{align}
Then \eqref{estR15}-\eqref{estR18} implies
\begin{align}\label{estR19}
\|R^{2,\ep,M}\|_{\widetilde{L}^1_T(H^{\f d 2 -\vartheta-1})}\lesssim M^{-\f \vartheta 2}XP.
\end{align}
We also observe that
\begin{align*}
R^{2,\ep,M}=(\mathcal{Q}_2-\mathcal{Q}_2^\ep)(\vV,\vV^M)+(\mathcal{Q}_2-\mathcal{Q}_2^\ep)(\vV^M,\vV_M).
\end{align*}
Similar to \eqref{estR15} and \eqref{estR16},
\begin{align}\label{estQ8}
 \|\Big(\mathcal{Q}_2^\ep(\vV,\vV^M),\mathcal{Q}_2^\ep(\vV^M,\vV_M)\Big)\|_{L^2_T(H^{\f d 2 -\vartheta-2})} \lesssim M^{-\vartheta} \|\vV\|_{L^\infty_T(B^{\f d 2-1}_{2,1})}\|\vV\|_{L^2_T(B^{\f d 2}_{2,1})}.
\end{align}

From \eqref{minski2}, \eqref{estM3} and \eqref{estQ5}-\eqref{estQ6},
\[
\|\mathcal{Q}_2(\vV,\vV^M)\|_{H^{\f d 2 -\vartheta-2}}\lesssim \|\vV\|_{B^{\f 1 2}_{2,1}}\|\vV^M\|_{B^{\f d 2-1-\theta}_{2,1}}\lesssim M^{-\theta}\|\vV\|_{B^{\f d 2}_{2,1}}\|\vV\|_{B^{\f d 2-1}_{2,1}}, \text{ if }\ \  \f1 2<\f d 2-\theta<\f 3 2,
\]
and
\begin{align*}
\|\mathcal{Q}_2(\vV,\vV^M)\|_{H^{\f d 2 -\vartheta-2}}&\lesssim \|\vV\|_{B^{\f d 2-1-\theta}_{2,1}}\|\vV^M\|_{B^{\f d 2-1-\theta}_{2,1}}
 \lesssim  M^{-\theta}\|\vV\|_{B^{\f d 2}_{2,1}}\|\vV\|_{B^{\f d 2-1}_{2,1}}, \text{ if }\ \  \f d 2-\theta\geq\f 3 2.
\end{align*}
Hence
\begin{align}\label{estQ9}
\|\Big(\mathcal{Q}_2(\vV,\vV^M),\mathcal{Q}_2(\vV^M,\vV_M)\Big)\|_{L^2_T(H^{\f d 2 -\vartheta-2})} \lesssim M^{-\vartheta} \|\vV\|_{L^\infty_T(B^{\f d 2-1}_{2,1})}\|\vV\|_{L^2_T(B^{\f d 2}_{2,1})}.
\end{align}
Then \eqref{estQ8} and \eqref{estQ9} yield
\begin{align}
 \|R^{3,\ep,M}\|_{L^2_T(H^{\f d 2 -\vartheta-2})} &\lesssim M^{-\vartheta} \|\vV\|_{L^\infty_T(B^{\f d 2-1}_{2,1})}\|\vV\|_{L^2_T(B^{\f d 2}_{2,1})}\nonumber\\
 &\lesssim M^{-\vartheta} X^2.\label{estR20}
\end{align}
We conclude that
\begin{align}\label{estR28}
  \|R^{\ep,M}\|_{\widetilde{L}^1_T(H^{\f d 2 -\vartheta-1})+L^2_T(H^{\f d 2 -\vartheta})}\lesssim M^{-\f\vartheta 2}\left( \mathbb{Q}_f+ Y^{\zeta,\ep}(T)+X^2+P^2\right).
\end{align}

\leftline{\emph{Estimate for $F^\ep$.}} From Lemma \ref{LeD3}, we have that
\begin{align}\label{estR20}
 \|F^\ep\|_{\widetilde{L}^1_T(H^{\f d 2 -\vartheta-1})+L^2_T(H^{\f d 2-\vartheta-2})} \lesssim & \|\vv\cdot\Grad \vWe\|_{\widetilde{L}^1_T(H^{\f d 2 -\vartheta-1})}+\|\vWe\cdot\Grad \mathbb{P}\vue\|_{\widetilde{L}^1_T(H^{\f d 2 -\vartheta-1})}\\
 &+\|I(\ep\vae)\mathcal{A}\vue\|_{\widetilde{L}^1_T(H^{\f d 2 -\vartheta-1})}+\|\vae \widetilde{K}(\ep\vae)\Grad \vae\|_{\widetilde{L}^1_T(H^{\f d 2 -\vartheta-1})}\nonumber\\
 &+\|\mathcal{Q}^\ep_1(\underline{\vWe},\vVe)\|_{\widetilde{L}^1_T(H^{\f d 2 -\vartheta-1})}+ \|\mathcal{Q}_1^\ep(\f{\widehat{\vWe}_0}{\sqrt{\Td}},\vVe)\|_{L^2_T(H^{\f d 2 -\vartheta-2})}\nonumber
\end{align}
It is clear that
\begin{align}\label{estR21}
  \|\vv\cdot\Grad \vWe\|_{\widetilde{L}^1_T(H^{\f d 2 -\vartheta-1})} &\lesssim   \|\vv\|_{L^2_T(\underline{B}^{\f d 2 }_{2,1})\cap L^\infty_T(B^{\f d 2 -1}_{2,1}) } \|\vWe\|_{L^2_T(\underline{B}^{\f d 2 -\vartheta}_{2,1})\cap L^1_T(\underline{B}^{\f d 2 +1-\vartheta}_{2,1})}\nonumber\\
  &\lesssim P\vWe_\vartheta(T)
\end{align}
\begin{align}\label{estR22}
  \|\vWe\cdot\Grad \mathbb{P}\vue\|_{\widetilde{L}^1_T(H^{\f d 2 -\vartheta-1})} \lesssim   \|\mathbb{P}\vue\|_{L^1_T(\underline{B}^{\f d 2 +1}_{2,1})} \|\vWe\|_{L^\infty_T(B^{\f d 2 -\vartheta-1}_{2,1})}\lesssim Y^{\zeta,\ep}(T)\vWe_\vartheta(T),
\end{align}
\begin{align}\label{estR23}
\|I(\ep\vae)\mathcal{A}\vue\|_{\widetilde{L}^1_T(H^{\f d 2 -\vartheta-1})} &\lesssim  \|\ep\|\vae\|_{B^{\f d 2}_{2,1}}\|\vue\|_{\underline{B}^{\f d 2+1-\vartheta}_{2,1}}\|_{L^1(0,T)}\nonumber\\
&\lesssim \ep^{\vartheta} \left(\|\ep\vae\|_{L^\infty_T(B^{\f d 2}_{2,1})}\|\vue\|^{h;\f {\beta_0} {\ep}}_{L^1_T(B^{\f d 2+1}_{2,1})}+\|\vae\|_{L^2_T(B^{\f d 2}_{2,1})}\|\vue\|^{l;\f {\beta_0} {\ep}}_{L^2_T(\underline{B}^{\f d 2}_{2,1})}\right)\nonumber\\
&\lesssim \ep^{\vartheta} \left(Y^{\zeta,\ep}(T)\right)^2,
\end{align}
\begin{align}\label{estR24}
\|\vae \widetilde{K}(\ep\vae)\Grad \vae\|_{\widetilde{L}^1_T(H^{\f d 2 -\vartheta-1})} &\lesssim  \|\ep\|\vae\|_{B^{\f d 2}_{2,1}}\|\vae\Grad \vae\|_{B^{\f d 2-1-\vartheta}_{2,1}}\|_{L^1(0,T)}\nonumber\\
&\lesssim \|\ep\|\vae\|_{B^{\f d 2}_{2,1}} \|\vae\|_{B^{\f d 2}_{2,1}}\|\vae\|_{B^{\f d 2-\vartheta}_{2,1}}\|_{L^1(0,T)}\nonumber\\
&\lesssim \ep^{\vartheta} \Bigg(\|\ep\vae\|_{L^\infty_T(B^{\f d 2}_{2,1})}\|\vae\|_{L^2_T(B^{\f d 2}_{2,1})}\|\vae\|^{h;\f {\beta_0} {\ep}}_{L^2_T(B^{\f d 2}_{2,1})}\nonumber\\
&\qquad\qquad+\|\vae\|_{L^2_T(B^{\f d 2}_{2,1})}\|\vae\|_{L^2_T(B^{\f d 2}_{2,1})}\|\vae\|^{l;\f {\beta_0} {\ep}}_{L^\infty_T(B^{\f d 2-1}_{2,1})}\Bigg)\nonumber\\
&\lesssim \ep^{\vartheta} \left(Y^{\zeta,\ep}(T)\right)^3.
\end{align}
\begin{align}\label{estR25}
 \|\mathcal{Q}^\ep_1(\underline{\vWe},\vVe)\|_{\widetilde{L}^1_T(H^{\f d 2 -\vartheta-1})}&\lesssim \|\vVe\|_{L^2_T(B^{\f d 2}_{2,1})} \|\underline{\vWe}\|_{L^2_T(H^{\f d 2-\vartheta})}\lesssim Y^{\zeta,\ep}(T)\vWe_\vartheta(T).
\end{align}
From  \eqref{minski2}  and the definition of $\mathcal{Q}_1^\ep$,
\[
 \|\mathcal{Q}_1^\ep(\f{\widehat{\vWe}_0}{\sqrt{\Td}},\vVe)\|_{H^{\f d 2-\vartheta-2}} \lesssim |\widehat{\vWe}_0|\|\vVe\|_{H^{\f d 2-\vartheta-1}}
 \lesssim \|\vWe\|_{B^{\f d 2-\vartheta-1}_{2,1}} \|\vVe\|_{H^{\f d 2 -\vartheta}},
\]
then
\begin{align}\label{estR26}
 \|\mathcal{Q}_1^\ep(\f{\widehat{\vWe}_0}{\sqrt{\Td}},\vVe)\|_{L^2_T(H^{\f d 2 -\vartheta-2})}&\lesssim \|\vWe\|_{L^\infty_T(B^{\f d 2-\vartheta-1}_{2,1})} \|\vVe\|_{L^2_T(H^{\f d 2 -\vartheta})} \lesssim Y^{\zeta,\ep}(T)\vWe_\vartheta(T).
\end{align}
Inserting \eqref{estR21}-\eqref{estR26} into \eqref{estR20}, we find that
\begin{align}\label{estR27}
 \|F^\ep\|_{\widetilde{L}^1_T(H^{\f d 2 -\vartheta-1})+L^2_T(H^{\f d 2-\vartheta-2})} \lesssim (Y^{\zeta,\ep}(T)+P)\vWe_\vartheta(T)+\ep^{\vartheta} \Big((Y^{\zeta,\ep}(T))^2+ (Y^{\zeta,\ep}(T))^3\Big).
\end{align}
\leftline{\emph{Estimate for $\widetilde{\mathcal{Q}}_1^\ep(\vv,\widetilde{R}_M^\ep), \widecheck{\mathcal{Q}}_1^\ep(\vv,\widetilde{R}_M^\ep) $ and $\mathcal{Q}_2^\ep(\vV,\widetilde{R}_M^\ep)$}.}
Similarly to \eqref{estR26}, we have
\begin{align}\label{estR29}
 \|\widecheck{\mathcal{Q}}_1^\ep(\f{\widehat{\vv}_0}{\sqrt{\Td}},\widetilde{R}_M^\ep)\|_{L^2_T(H^{\f d 2 -\vartheta-2})}\lesssim  \|\vv\|_{L^\infty_T(B^{\f d 2-\vartheta-1}_{2,1})} \|\widetilde{R}_M^\ep\|_{L^2_T(H^{\f d 2 -\vartheta})}.
\end{align}
From Lemma \ref{leD1}, we find that
\begin{align}
 \|\widetilde{\mathcal{Q}}_1^\ep(\vv,\tilde{R}_M^\ep)\|_{\widetilde{L}^1_T(H^{\f d 2 -\vartheta-1})}&\lesssim \|\underline{\vv}\|_{L^2_T(B^{\f d 2}_{2,1})}\|\widetilde{R}_M^\ep\|_{L^2_T(H^{\f d 2 -\vartheta})},\label{estR30}\\
\|\widecheck{\mathcal{Q}}_1^\ep(\underline{\vv},\widetilde{R}_M^\ep)\|_{L^2_T(H^{\f d 2 -\vartheta-2})}&\lesssim \|\underline{\vv}\|_{L^2_T(B^{\f d 2}_{2,1})}\|\widetilde{R}_M^\ep\|_{L^\infty_T(H^{\f d 2 -\vartheta-1})},\label{estR31}\\
\|\mathcal{Q}_2^\ep(\vV,\widetilde{R}_M^\ep\|_{L^2_T(H^{\f d 2 -\vartheta-2})}&\lesssim  \|\vV\|_{L^2_T(B^{\f d 2}_{2,1})}\|\widetilde{R}_M^\ep\|_{L^\infty_T(H^{\f d 2 -\vartheta-1})}.\label{estR32}
\end{align}
Using \eqref{estR1} and \eqref{estR29}-\eqref{estR32}, we obtain that
\begin{align}\label{estR33}
  \|\widetilde{\mathcal{Q}}_1^\ep(\vv,\tilde{R}_M^\ep)\|&_{\widetilde{L}^1_T(H^{\f d 2 -\vartheta-1})}+\|\widecheck{\mathcal{Q}}_1^\ep(\vv,\widetilde{R}_M^\ep)\|_{L^2_T(H^{\f d 2 -\vartheta-2})}+\|\mathcal{Q}_2^\ep(\vV,\widetilde{R}_M^\ep\|_{L^2_T(H^{\f d 2 -\vartheta-2})}\nonumber\\
  &\lesssim C_M(X+P)(Q_f+Y^{\zeta,\ep}(T)+X^2+P^2)
\end{align}

\leftline{\emph{Estimate for $\widetilde{R}^{t,\ep}_M$}.}
The estimation methods for $\widetilde{R}^{1,t,\ep}_M,\widetilde{R}^{2,t,\ep}_M,\widetilde{R}^{3,t,\ep}_M$ are similar to those for $\widetilde{R}^{1,\ep}_M,\widetilde{R}^{2,\ep}_M,\widetilde{R}^{3,\ep}_M$; hence, we will not reiterate the process here. The results obtained are as follows:
\begin{align*}
  \|\widetilde{R}^{1,t,\ep}_M\|_{\widetilde{L}^1_T(H^{\f d 2-\vartheta-1})}&\lesssim M^{\f d 2+S-1} \|\p_t\mathbb{Q}f\|_{L^1_T(H^{-S})}+M\|\p_t\vv\|_{L^1_T(B^{\f d 2-1}_{2,1})}\|\vv\|_{\widetilde{L}^\infty_T(B^{\f d 2-1}_{2,1})},\\
   \|\widetilde{R}^{2,t,\ep}_M\|_{\widetilde{L}^1_T(H^{\f d 2-\vartheta-1})}&\lesssim C_M^1M^{2+2\vartheta}\Bigg(\|\p_t\vV\|_{L^1_T(H^{\f d 2-1-\vartheta})}\|\vv\|_{L^\infty_T(H^{\f d 2-1-\vartheta})}\\
   &\qquad\qquad\qquad\quad+\|\p_t\vv\|_{L^1_T(H^{\f d 2-1-\vartheta})}\|\vV\|_{L^\infty_T(H^{\f d 2-1-\vartheta})}\Bigg),\\
   \|\widetilde{R}^{3,t,\ep}_M\|_{\widetilde{L}^1_T(H^{\f d 2-\vartheta-1})}&\lesssim C_M^2M^{2+2\vartheta}\|\p_t\vV\|_{L^1_T(H^{\f d 2-1-\vartheta})}\|\vV\|_{L^\infty_T(H^{\f d 2-1-\vartheta})}.
\end{align*}
Notice that $\mathbb{Q}f\in L^2_T(H^{\f {\f d 2-1-S}{2}})$. Let $s^{\ast}=\min\{\f d 2-2, \f {\f d 2-1-S}{2}\}$, then from \eqref{minski2} and \eqref{estM3}, we have that
\begin{align*}
  \|\widetilde{S}^{t,\ep}_M\|_{L^2_T(H^{\f d 2-\vartheta-2})}&\lesssim M^{\f d 2-1-s^{\ast}}  \|\p_t\vVe\|_{L^2_T(H^{s^{\ast}-\vartheta})}\lesssim M^{\f d 2-1-s^{\ast}}  \|\p_t\vVe\|_{L^2_T(B^{s^{\ast}}_{2,\infty})}.
\end{align*}
It is easy to see that $ M^{\f d 2-1-s^{\ast}}\leq \max\{M,M^{\f{\f d 2+S-1}{2}}\}$.
For the time derivatives, we need to utilize the following equations,
\[
\p_t\vv=\mathbb{P}f+\mu\Delta\vv-\mathbb{P}(\vv\cdot\Grad\vv),
\]
\[
\p_t\vV=\f\nu2\De\vV-\mathcal{Q}_1(\vv,\vV)-\mathcal{Q}_2(\vV,\vV),
\]
\begin{align*}
 \p_t\vVe=& \nu\mathcal{L}\Big(-\f t {\ep}\Big)\left(\begin{array}{l}
\quad 0 \\
\De\mathbb{Q}\vue
\end{array}\right)-\mathcal{L}\Big(-\f t {\ep}\Big)\left(\begin{array}{l}
\qquad\qquad\qquad\Div(\vae\vue) \\
\mathbb{Q}(\vue\cdot\Grad\vue+\widetilde{K}(\ep\vae)\vae\Grad\vae+\kappa\vae\Grad\vae)
\end{array}\right)\\
&-\mathcal{L}\Big(-\f t {\ep}\Big)\left(\begin{array}{l}
\qquad\quad 0 \\
\mathbb{Q}(I(\ep\vae)\mathcal{A}\vue)
\end{array}\right)+\mathcal{L}\Big(-\f t {\ep}\Big)\left(\begin{array}{l}
  0 \\
\mathbb{Q}f
\end{array}\right).
\end{align*}
By estimating the right-hand side of the aforementioned equations, one can obtain estimates for the time derivative.
From \eqref{estr1},
\[
\|\mathbb{P}(\vv\cdot\Grad\vv)\|_{L^1_T(B^{\f d 2-1}_{2,1})}\lesssim \|\vv\|_{L^1_T(\underline{B}^{\f d 2+1}_{2,1})}\|\vv\|_{L^\infty_T(B^{\f d 2-1}_{2,1})}\lesssim P^2,
\]
\[
\|\Div(\vae\vue)\|_{L^2_T(B^{\f d 2-2}_{2,\infty})}\lesssim \|\vae\|_{L^2_T(B^{\f d 2}_{2,1})\cap L^\infty_T(B^{\f d 2-1}_{2,1})}\|\vue\|_{L^\infty_T(B^{\f d 2-1}_{2,1})\cap L^2_T(\underline{B}^{\f d 2}_{2,1})}\lesssim (Y^{\zeta,\ep}(T))^2,
\]
\[
\|\vue\cdot\Grad\vue\|_{L^2_T(B^{\f d 2-2}_{2,\infty})}\lesssim\|\vue\|_{L^2_T(\underline{B}^{\f d 2}_{2,1})}\|\vue\|_{L^\infty_T(B^{\f d 2-1}_{2,1})}\lesssim (Y^{\zeta,\ep}(T))^2,
\]
\[
\|\vae\Grad\vae\|_{L^2_T(B^{\f d 2-2}_{2,\infty})}\lesssim \|\vae\|_{L^2_T(B^{\f d 2}_{2,1})}\|\vae\|_{L^\infty_T(B^{\f d 2-1}_{2,1})}\lesssim (Y^{\zeta,\ep}(T))^2,
\]
\[
\|\widetilde{K}(\ep\vae)\vae\Grad\vae\|_{L^2_T(B^{\f d 2-2}_{2,\infty})}\lesssim \ep \|\vae\|_{L^\infty_T(B^{\f d 2}_{2,1})}\|\vae\Grad\vae\|_{L^2_T(B^{\f d 2-2}_{2,\infty})}\lesssim (Y^{\zeta,\ep}(T))^3,
\]
\[
\|I(\ep\vae)\mathcal{A}\vue\|_{L^2_T(B^{\f d 2-2}_{2,\infty})}\lesssim \ep \|\vae\|_{L^\infty_T(B^{\f d 2}_{2,1})}\|\vue\|_{L^2_T(\underline{B}^{\f d 2}_{2,1})}\lesssim (Y^{\zeta,\ep}(T))^2.
\]
From \eqref{estQ4} and \eqref{estQ7},
\[
\|\mathcal{Q}_1(\vv,\vV)\|_{L^1_T(H^{\f d 2-\vartheta-1})}\lesssim \|\vV\|_{L^2_T(B^{\f d 2}_{2,1})\cap L^1_T(B^{\f d 2+1}_{2,1})}\|\vv\|_{L^\infty_T(B^{\f d 2-1}_{2,1})\cap L^2_T(\underline{B}^{\f d 2}_{2,1})} \lesssim PX,
\]
\[
\|\mathcal{Q}_2(\vV,\vV)\|_{L^1_T(B^{\f d 2-1}_{2,1})}\lesssim \|\vV\|_{L^2_T(B^{\f d 2}_{2,1})}\|\vV\|_{L^2_T(B^{\f d 2}_{2,1})}\lesssim X^2.
\]
It follows that
\begin{align}
 \|\p_t\vv\|_{L^1_T(B^{\f d 2-1}_{2,1})}&\lesssim \mathbb{P}_f+P+P^2,\label{estp15}\\
 \|\p_t\vV\|_{L^1_T(H^{\f d 2-1-\vartheta})} &\lesssim  X+PX+X^2,\label{estp16}\\
  \|\p_t\vVe\|_{L^2_T(B^{s^{\ast}}_{2,\infty})}&\lesssim \mathbb{Q}_f+Y^{\zeta,\ep}(T)+(Y^{\zeta,\ep}(T))^2+(Y^{\zeta,\ep}(T))^3.\label{estp17}
\end{align}
 Now we can conclude that
 \begin{align}\label{estR34}
\|\widetilde{R}_M^{t,\ep}\|_{\widetilde{L}^1_T(H^{\f d 2 -\vartheta-1})\cap L_T^2(H^{\f d 2-\vartheta})}\leq C_M&\Big(\mathbb{Q}_f+\mathbb{Q}_fX+\mathbb{P}_f(P+X)+Y^{\zeta,\ep}(T)+P^2+X^2\\
&+(Y^{\zeta,\ep}(T))^2+P^3+X^3+(Y^{\zeta,\ep}(T))^3+X^4+(Y^{\zeta,\ep}(T))^4\Big).\nonumber
 \end{align}

 Let us define $M=\mathcal{X}^{-1}(\ep^{-1})$, where $\mathcal{X}^{-1}$ is the inverse function of $\mathcal{X}$ with $\mathcal{X}(M):=C_M M^{\f \vartheta 2}$. Substituting \eqref{estR2}, \eqref{estR35}-\eqref{estZ}, \eqref{estR28}, \eqref{estR27}, \eqref{estR33}, \eqref{estR34} into \eqref{est1}, and then adding the resulting expression to \eqref{est2}, we establish Proposition \ref{pro1} with $\widetilde{\tau}(\ep):=\max\{(\mathcal{X}^{-1}(\ep^{-1}))^{-\f \vartheta 2}, \ep^\vartheta\}$.\ \ $\Box$

\subsection{The proof of Proposition \ref{pro2}}\label{newsec2}
Starting from systems \eqref{equ2} and \eqref{equ3}, and noting that $\mathbb{P}(\mathbb{Q}\vue\cdot\Grad\mathbb{Q}\vue)=0$, we derive the following evolution equation for $\vWe$:
\begin{equation}\label{equ10}
\bali
  \p_t\vWe&+\mathbb{P}(A^{\ep}\cdot\Grad \vWe+\vWe\cdot\Grad A^{\ep})-\mu\De\vWe\\
  &=-\mathbb{P}(\mathbb{Q}\vue\cdot\Grad \vv+\vv\cdot\Grad\mathbb{Q}\vue+\vWe\cdot\Grad \vWe+I(\ep\vae)\mathcal{A}\vue),
  \eali
\end{equation}
where $A^\ep:=\mathbb{Q}\vue+\vv$.

In analogy with the estimation method for $\vZe$, we introduce the function $\psi_M^\ep:=\vWe+\ep\mathbb{P}\widetilde{R}_M^\ep$ with
\begin{align*}
\widetilde{R}_M^\ep(t,x):=-c_d\sum\limits_{m}\left(\sum\limits_{\begin{array}{c}
                                                                  \alpha\in\{-1,1\},k+l=m\\
                                                                 |k|,|l|\leq M
                                                                \end{array}}
                                                                \left(\f{(l\cdot k)\widehat{\vv}_k+(l\cdot\widehat{\vv}_k)l}{|l|^2}\right)\widehat{\vV}_l^{\alpha,\ep}e^{i|l|sg(l)\f{\alpha t}{\ep}}
\right)e^{im\cdot x}.
\end{align*}

Noting that

\begin{align*}
\mathbb{Q}\vue\cdot\Grad \vv=-ic_d\sum\limits_{m}\left(\sum\limits_{k+l=m}\alpha sg(l)\f{(l\cdot k)\widehat{\vv}_k}{|l|}
\widehat{\vV}_l^{\alpha,\ep} e^{ilsg(l)\f{\alpha t}\ep}
\right)e^{im\cdot x},
\end{align*}

\begin{align*}
\vv \cdot\Grad \mathbb{Q}\vue=-ic_d\sum\limits_{m}\left(\sum\limits_{k+l=m}\alpha  sg(l)\f{(l\cdot \widehat{\vv}_k)l}{|l|}
\widehat{\vV}_l^{\alpha,\ep} e^{ilsg(l)\f{\alpha t}\ep}
\right)e^{im\cdot x},
\end{align*}
we obtain that
\begin{align}\label{equ11}
\ep\p_t\widetilde{R}_M^\ep=R_M^\ep+\ep\widetilde{R}_M^{t,\ep},
\end{align}
where $R_M^\ep$ is the low frequency part of $R^\ep:=\mathbb{Q}\vue\cdot\Grad \vv+\vv \cdot\Grad \mathbb{Q}\vue$ obtained by keeping only the indices $l$ and $k$ such that $|l|$, $|k|\leq M$ and
\begin{align*}
\widetilde{R}_M^{t,\ep}(t,x):=-c_d\sum\limits_{m}\left(\sum\limits_{\begin{array}{c}
                                                                  \alpha\in\{-1,1\},k+l=m\\
                                                                 |k|,|l|\leq M
                                                                \end{array}}
                                                               \p_t \left(\f{(l\cdot k)\widehat{\vv}_k+(l\cdot\widehat{\vv}_k)l}{|l|^2}\widehat{\vV}_l^{\alpha,\ep}\right)e^{i|l|sg(l)\f{\alpha t}{\ep}}
\right)e^{im\cdot x}.
\end{align*}
From equations \eqref{equ10} and \eqref{equ11}, we obtain that
\begin{align}
 &\p_t\psi_M^\ep+\mathbb{P}(A^\ep\cdot \Grad \psi_M^\ep+\psi_M^\ep\cdot\Grad A^\ep)-\mu\De\psi_M^\ep\nonumber\\
 &\quad =-\mathbb{P}R^{\ep,M}+\ep\mathbb{P}(A^\ep\Grad\mathbb{P}\widetilde{R}_M^\ep+\mathbb{P}\widetilde{R}_M^\ep\cdot\Grad A^\ep-\mu\De\widetilde{R}_M^\ep)+\ep\mathbb{P}\widetilde{R}_M^{t,\ep}\nonumber\\
 &\qquad -\mathbb{P}(\vWe\cdot\Grad\vWe+I(\ep\vae)\mathcal{A}\vue)\label{equ12}
\end{align}
where $R^{\ep,M}:=R^\ep-R_M^\ep$.

Applying Proposition 7.4 in \cite{D3} to system \eqref{equ12}, we obtain
\begin{align}\label{estp1}
 &\|\psi_M^\ep\|_{\widetilde{L}_T^\infty(B_{2,1}^{\f d 2 -1-\vartheta})\cap L^1_T(\underline{B}_{2,1}^{\f d 2 +1-\vartheta})}\nonumber\\
 &\quad\lesssim e^{C\|A^\ep\|_{L^1_T(\underline{B}_{2,1}^{\f d 2 +1})}}\nonumber\\
 &\quad\quad\times\Bigg(\|\mathbb{P}R^{\ep,M}\|_{L^1_T(B_{2,1}^{\f d 2 -1-\vartheta})}+\|\vWe\cdot\Grad \vWe\|_{L^1_T(B_{2,1}^{\f d 2 -1-\vartheta})}+\|I(\ep\vae)\mathcal{A}\vue\|_{L^1_T(B_{2,1}^{\f d 2 -1-\vartheta})}\nonumber\\
 &\quad\quad\quad+\ep\Big(\|\widetilde{R}^\ep_M(0)\|_{B_{2,1}^{\f d 2 -1-\vartheta}}+\|\mathbb{P}\widetilde{R}_M^{\ep}\|_{L^1_T(B_{2,1}^{\f d 2 +1-\vartheta})}+\|A^\ep\cdot\Grad\mathbb{P}\widetilde{R}_M^\ep\|_{L^1_T(B_{2,1}^{\f d 2 -1-\vartheta})}\nonumber\\
 &\quad\quad\quad\quad+\|\mathbb{P}\widetilde{R}_M^\ep\cdot\Grad A^\ep\|_{L^1_T(B_{2,1}^{\f d 2 -1-\vartheta})}+\|\mathbb{P}\widetilde{R}_M^{t,\ep}\|_{L^1_T(B_{2,1}^{\f d 2 -1-\vartheta})}\Big)\Bigg).
\end{align}

First, we claim that
\begin{align}\label{estp2}
\|\mathbb{P}\widetilde{R}_M^{\ep}\|_{\widetilde{L}^\infty_T(B_{2,1}^{\f d 2 -1-\vartheta})\cap L^1_T(B_{2,1}^{\f d 2 +1-\vartheta})}\lesssim MPY^{\zeta,\ep}(T)
\end{align}
To this end, we introduce the following notation:
\begin{align*}
  \widetilde{\vV}^\ep(t,x)&:=-c_d\sum\limits_{k,\alpha}\widehat{\vV}_k^{\alpha,\ep}e^{i|k|sg(k)\f{\alpha t}\ep}e^{ik\cdot x},\\
    \widecheck{\vV}^\ep(t,x)&:=ic_d\sum\limits_{k,\alpha}\alpha sg(k)\widehat{\vV}_k^{\alpha,\ep}e^{i|k|sg(k)\f{\alpha t}\ep}e^{ik\cdot x},\\
    \widetilde{\vV}^{t,\ep}(t,x)&:=-c_d\sum\limits_{k,\alpha}\p_t\widehat{\vV}_k^{\alpha,\ep}e^{i|k|sg(k)\f{\alpha t}\ep}e^{ik\cdot x}.\\
\end{align*}
It is clear that
\begin{align}\label{estp4}
\|\mathbb{P}\widetilde{R}_M^{\ep}\|_{\widetilde{L}^\infty_T(B_{2,1}^{\f d 2 -1-\vartheta})\cap L^1_T(B_{2,1}^{\f d 2 +1-\vartheta})}\lesssim M\|\mathbb{P}\widetilde{R}_M^{\ep}\|_{\widetilde{L}^\infty_T(B_{2,1}^{\f d 2 -2-\vartheta})\cap L^1_T(B_{2,1}^{\f d 2 -\vartheta})}
\end{align}
and
\begin{align}\label{estp3}
 \| \widetilde{\vV}_M^\ep\|_{\widetilde{L}^q_T(B_{2,r}^s)}+ \| \widecheck{\vV}_M^\ep\|_{\widetilde{L}^q_T(B_{2,r}^s)}\approx \|\vV_M^\ep\|_{\widetilde{L}^q_T(B_{2,r}^s)},\ \  \| \widetilde{\vV}_M^{t,\ep}\|_{\widetilde{L}^q_T(B_{2,r}^s)}\approx  \| \p_t{\vV}_M^\ep\|_{\widetilde{L}^q_T(B_{2,r}^s)},
 \end{align}
where  $1\leq q,r\leq\infty$ and $s\in \mathbb{R}$.

Note that $\widetilde{R}_M^\ep$ can be  expressed in the following form,
\[
\widetilde{R}_M^\ep=\Grad \vv_M\cdot \Grad\De^{-1}\widetilde{\vV}_M^\ep+\Grad^2\De^{-1}\widetilde{\vV}_M^\ep\cdot\vv_M,
\]
which implies that
\[
\mathbb{P}(\widetilde{R}_M^\ep)=\mathbb{P}\big(\Grad \vv_M\cdot \Grad\De^{-1}\widetilde{\vV}_M^\ep+\Grad^2\De^{-1}\widetilde{\vV}_M^\ep\cdot\underline{\vv_M}\big)
\]
From \eqref{minski2} and \eqref{estp3}, we have that
\begin{align}
  \|\Grad \vv_M\cdot \Grad\De^{-1}\widetilde{\vV}_M^\ep\|_{ L^1_T(B_{2,1}^{\f d 2-\vartheta})}&\lesssim \|\vv_M\|_{ L^1_T(\underline{B}_{2,1}^{\f d 2+1})}\|\vV^\ep_M\|_{ L^\infty_T(\underline{B}_{2,1}^{\f d 2-1})},\label{estp5}\\
  \|\Grad^2\De^{-1}\widetilde{\vV}_M^\ep\cdot\underline{\vv_M}\|_{ L^1_T(B_{2,1}^{\f d 2-\vartheta})}&\lesssim \|\vv_M\|_{ L^2_T(\underline{B}_{2,1}^{\f d 2})}\|\vV^\ep_M\|_{ L^2_T(\underline{B}_{2,1}^{\f d 2})},\label{estp6}\\
  \|\Grad \vv_M\cdot \Grad\De^{-1}\widetilde{\vV}_M^\ep\|_{ \widetilde{L}^\infty_T(B_{2,1}^{\f d 2-2-\vartheta})}&\lesssim   \|\Grad \vv_M\cdot \Grad\De^{-1}\widetilde{\vV}_M^\ep\|_{ L^\infty_T(B_{2,\infty}^{\f d 2-2})}\nonumber\\
  &\lesssim \|\vv_M\|_{ L^\infty_T(B_{2,1}^{\f d 2-1})}\|\vV^\ep_M\|_{ L^\infty_T(B_{2,1}^{\f d 2-1})},\label{estp7}\\
  \|\Grad^2\De^{-1}\widetilde{\vV}_M^\ep\cdot\underline{\vv_M}\|_{ \widetilde{L}^\infty_T(B_{2,1}^{\f d 2-2-\vartheta})}&\lesssim \|\Grad^2\De^{-1}\widetilde{\vV}_M^\ep\cdot\underline{\vv_M}\|_{ L^\infty_T(B_{2,\infty}^{\f d 2-2})}\nonumber\\
  &\lesssim \|\vv_M\|_{ L^\infty_T(B_{2,1}^{\f d 2-1})}\|\vV^\ep_M\|_{ L^\infty_T(B_{2,1}^{\f d 2-1})}.\label{estp8}
\end{align}
Then \eqref{estp2} follows from \eqref{estp4}, \eqref{estp5}, \eqref{estp6}, \eqref{estp7} and \eqref{estp8}.

Next we turn to the estimate of $\psi_M^\ep$. By calculations similar to those in \eqref{estp2}, we get
\begin{align}\label{estp9}
\|\mathbb{P}\widetilde{R}^{\ep,M}(0)\|_{B_{2,1}^{\f d 2 -1-\vartheta}}\lesssim MX_0P_0.
\end{align}

From \eqref{estp2}, we infer that
\begin{align}\label{estp13}
\|A^\ep\cdot\Grad\mathbb{P}\widetilde{R}_M^\ep\|_{L^1_T(B_{2,1}^{\f d 2 -1-\vartheta})}&\lesssim \|A^\ep\|_{ L^\infty_T(B_{2,1}^{\f d 2-1})}\|\mathbb{P}\widetilde{R}_M^\ep\|_{L^1_T(\underline{B}_{2,1}^{\f d 2 +1-\vartheta})}\nonumber\\
&\lesssim \Big(P+Y^{\zeta,\ep}(T)\Big)MPY^{\zeta,\ep}(T)
\end{align}
and
\begin{align}\label{estp14}
\|\mathbb{P}\widetilde{R}_M^\ep\cdot\Grad A^\ep\|_{L^1_T(B_{2,1}^{\f d 2 -1-\vartheta})}&\lesssim \|A^\ep\|_{ L^2_T(\underline{B}_{2,1}^{\f d 2})}\|\mathbb{P}\widetilde{R}_M^\ep\|_{L^2_T(B_{2,1}^{\f d 2 -\vartheta})}\nonumber\\
&\lesssim \Big(P+Y^{\zeta,\ep}(T)\Big)MPY^{\zeta,\ep}(T).
\end{align}

We observe that
\begin{align*}
\mathbb{P}\widetilde{R}_M^{t,\ep}=&\mathbb{P}\Big(\Grad \p_t\vv_M\cdot \Grad\De^{-1}\widetilde{\vV}_M^\ep+\Grad^2\De^{-1}\widetilde{\vV}_M^\ep\cdot\underline{\p_t \vv_M}\\
&+\Grad\vv_M\cdot\Grad\De^{-1}\widetilde{\vV}_M^{t,\ep}+\Grad^2\De^{-1}\widetilde{\vV}_M^{t,\ep}\cdot\underline{\vv_M}\Big).
\end{align*}
From \eqref{minski2}, \eqref{estp15}, \eqref{estp17} and \eqref{estp3},  we get that
\begin{align}
  \|\Grad \p_t\vv_M\cdot \Grad\De^{-1}\widetilde{\vV}_M^\ep\|_{L^1_T(B_{2,1}^{\f d 2 -1-\vartheta})}&\lesssim M\|\p_t\vv\|_{ L^1_T(\underline{B}_{2,1}^{\f d 2-1})}\|\vV^\ep\|_{L^\infty_T(B_{2,1}^{\f d 2 -1})}\nonumber\\
  &\lesssim MY^{\zeta,\ep}(T)(\mathbb{P}_f+P+P^2),\label{estp18}\\
  \|\Grad^2\De^{-1}\widetilde{\vV}_M^\ep\cdot\underline{\p_t \vv_M}\|_{L^1_T(B_{2,1}^{\f d 2 -1-\vartheta})}
  &\lesssim MY^{\zeta,\ep}(T)(\mathbb{P}_f+P+P^2),\label{estp19}
\end{align}
\begin{align}
  \|\Grad\vv_M\cdot\Grad\De^{-1}\widetilde{\vV}_M^{t,\ep}\|_{L^1_T(B_{2,1}^{\f d 2 -1-\vartheta})}&\lesssim M\|\vv\|_{ L^2_T(\underline{B}_{2,1}^{\f d 2})}\|\p_t\vV^\ep\|_{L^2_T(B_{2,1}^{\f d 2 -2-\vartheta})}\nonumber\\
  &\lesssim M\|\vv\|_{ L^2_T(\underline{B}_{2,1}^{\f d 2})}\|\p_t\vV^\ep\|_{L^2_T(B_{2,\infty}^{\f d 2 -2})}\nonumber\\
  &\lesssim MP\Big(\mathbb{Q}_f+Y^{\zeta,\ep}(T)+(Y^{\zeta,\ep}(T))^2+(Y^{\zeta,\ep}(T))^3\Big),\label{estp20}\\
   \|\Grad^2\De^{-1}\widetilde{\vV}_M^{t,\ep}\cdot\underline{\vv_M}\|_{L^1_T(B_{2,1}^{\f d 2 -1-\vartheta})}&\lesssim MP\Big(\mathbb{Q}_f+Y^{\zeta,\ep}(T)+(Y^{\zeta,\ep}(T))^2+(Y^{\zeta,\ep}(T))^3\Big)\label{estp21}.
\end{align}
It follows from \eqref{estp18}-\eqref{estp21} that
\begin{align}\label{estp22}
   \|\mathbb{P}\widetilde{R}_M^{t,\ep}\|_{L^1_T(B_{2,1}^{\f d 2 -1-\vartheta})}\lesssim &M\Big((\mathbb{P}_f)^2+(\mathbb{Q}_f)^2+P^2+(Y^{\zeta,\ep}(T))^2\nonumber\\
   &\quad+P^3+(Y^{\zeta,\ep}(T))^3+P^4+(Y^{\zeta,\ep}(T))^4\Big).
\end{align}
Through direct computation, we obtain  that
\begin{align}
 \|\vWe\cdot\Grad \vWe\|_{L^1_T(B_{2,1}^{\f d 2 -1-\vartheta})} &\lesssim  \|\vWe\|_{L^\infty_T(B_{2,1}^{\f d 2 -1})} \|\vWe\|_{L^1_T(\underline{B}_{2,1}^{\f d 2 +1-\vartheta})}\nonumber\\
  &\lesssim \left(D^{\zeta,\ep}(T)+\|\vv\|^{h;\zeta}_{\widetilde{L}^\infty_{T_0}(B^{\f d 2-1}_{2,1})} \right) \vWe_\vartheta(T).\label{estp10}
\end{align}
Similar to \eqref{estR23}, we have that
\begin{align}\label{estp11}
 \|I(\ep\vae)\mathcal{A}\vue\|_{L^1_T(B_{2,1}^{\f d 2 -1-\vartheta})}&\lesssim \ep^{\vartheta} \left(Y^{\zeta,\ep}(T)\right)^2.
\end{align}
Note that
\begin{align*}
\mathbb{P}(R^{\ep,M})&=\mathbb{P}\Big(\Grad \vv^M\cdot \Grad |D|^{-1}\widecheck{\vV}^\ep+\Grad \vv_M\cdot\Grad|D|^{-1}\widecheck{\vV}^{M,\ep}\\
&\quad\quad+\Grad^2|D|^{-1}\widecheck{\vV}^\ep\cdot\vv^M+\Grad^2|D|^{-1}\widecheck{\vV}^{M,\ep}\cdot\underline{\vv_M}\Big).
\end{align*}
Computing directly gives that
\begin{align*}
\|\Grad \vv^M\cdot \Grad |D|^{-1}\widecheck{\vV}^\ep\|_{L^1_T(B_{2,1}^{\f d 2 -1-\vartheta})}&\lesssim \|\vv^M\|_{ L^2_T(\underline{B}_{2,1}^{\f d 2-\vartheta})}\|\widecheck{\vV}^\ep\|_{ L^2_T(\underline{B}_{2,1}^{\f d 2})}\\
&\lesssim M^{-\vartheta} \|\vv\|_{ L^2_T(\underline{B}_{2,1}^{\f d 2})}\|\vV^\ep\|_{ L^2_T(\underline{B}_{2,1}^{\f d 2})},\\
\|\Grad \vv_M\cdot\Grad|D|^{-1}\widecheck{\vV}^{M,\ep}\|_{L^1_T(B_{2,1}^{\f d 2 -1-\vartheta})}&\lesssim M^{-\vartheta} \|\vv\|_{ L^2_T(\underline{B}_{2,1}^{\f d 2})}\|\vV^\ep\|_{ L^2_T(\underline{B}_{2,1}^{\f d 2})},\\
\end{align*}
\begin{align*}
\|\Grad^2|D|^{-1}\widecheck{\vV}^\ep\cdot\vv^M\|_{L^1_T(B_{2,1}^{\f d 2 -1-\vartheta})}&\lesssim \|\widecheck{\vV}^\ep\|_{ L^2_T(\underline{B}_{2,1}^{\f d 2})}\|\vv^M\|_{ L^2_T(B_{2,1}^{\f d 2-\vartheta})}\\
&\lesssim M^{-\vartheta} \|\vv\|_{ L^2_T(\underline{B}_{2,1}^{\f d 2})}\|\vV^\ep\|_{ L^2_T(\underline{B}_{2,1}^{\f d 2})},\\
\|\Grad^2|D|^{-1}\widecheck{\vV}^{M,\ep}\cdot\underline{\vv_M}\|_{L^1_T(B_{2,1}^{\f d 2 -1-\vartheta})}&\lesssim M^{-\vartheta} \|\vv\|_{ L^2_T(\underline{B}_{2,1}^{\f d 2})}\|\vV^\ep\|_{ L^2_T(\underline{B}_{2,1}^{\f d 2})}.\\
\end{align*}
It follows that
\begin{align}\label{estp12}
 \|\mathbb{P}R^{\ep,M}\|_{L^1_T(B_{2,1}^{\f d 2 -1-\vartheta})}\lesssim M^{-\vartheta}P Y^{\zeta,\ep}(T).
\end{align}

We select $M=\ep^{-\f{1}{1+\vartheta}}$ and substitute  \eqref{estp4}, \eqref{estp9}-\eqref{estp14}, and  \eqref{estp22}-\eqref{estp12} into \eqref{estp1}. By adding the resulting expression to \eqref{estp4}, we complete the proof of Proposition \ref{pro2}.\ \ $\Box$

\section{The relationship between  $Y^{\zeta,\ep}(T)$  and  $X^\ep(T)+P^\ep(T)$, as well as estimates for $ [\vae,\vue]^{\zeta,\ep}_{h,m}(T)$}\label{Sec4}
The proof of Proposition \ref{pro4} is based on the following three lemmas.

For $0<\beta<\infty$, $j\in\mathbb{Z}$, $\mathbf{1}_{h;\beta}(j)$ is defined as
\[
\mathbf{1}_{h;\beta}(j)=1, \text{ if } 2^j\geq\beta,\ \ \mathbf{1}_{h;\beta}(j)=0,\text{ if } 2^j<\beta.
\]
\begin{Lemma}\label{lema4.1}
Let $s,s_1,s_2\in\mathbb{R}$, $s+\f d 2=s_1+s_2$, $1\leq r\leq\infty$, $0<\zeta<\infty$, $f\in\dot{B}_{2,1}^{s_1}(\mathbb{R}^d)$ and $g\in\dot{B}_{2,r}^{s_2}(\mathbb{R}^d)$.
If $s_1\leq\f d 2$,  then
\begin{align}\label{estpri2}
\|T_{f}g\|^{h;\zeta}_{B_{2,r}^s}\lesssim\|f\|_{B_{2,1}^{s_1}}\|g\|^{h;\f\zeta 4}_{B_{2,r}^{s_2}};
\end{align}
If $s>-\f d 2$, then
\begin{align}\label{estpri3}
 \|R(f,g)\|^{h;\zeta}_{\dot{B}_{2,r}^s}\lesssim \|f\|^{h;\f\zeta {64}}_{\dot{B}_{2,1}^{s_1}}\|g\|^{h;\f\zeta {16}}_{\dot{B}_{2,r}^{s_2}}.
\end{align}
\end{Lemma}
\bProof
For the dyadic blocks  $\Delta_j(T_{f}g)$  satisfying $\zeta\leq2^j$, we deduce from \eqref{zero} and \eqref{ao1} that
\begin{align*}
&2^{js}\|\Delta_j(T_{f}g)\|_{L^2}\\
&\quad\lesssim 2^{js}\sum\limits_{|j-j'|\leq2}\|{S}_{j'-2}f\|_{L^{\infty}}\mathbf{1}_{h;2^{j_b}}(j')\|\Delta_{j'}g\|_{L^{2}}\\
&\quad\lesssim
\sum\limits_{|j-j'|\leq2}2^{(j-j')s}\mathbf{1}_{h;\f\zeta 4}(j')2^{j's_2}\|\Delta_{j'}g\|_{L^{2}}\mathbf{1}_{h;2^{j_b}}(j')2^{j'(s_1-\f d 2)}\|{S}_{j'-2}f\|_{L^{\infty}}
\end{align*}
Taking the $\textit{l}^r$-norm over the dyadic blocks  $\zeta\leq2^j$ and employing Young's inequality, we obtain that \footnote{We observe that $S_j f= \f{\hat{f}_0}{\sqrt{|\mathbb{T}_b^d|}}+\sum\limits_{j'\leq j-1}\Delta_{j'} f$. To handle the term $\f{\hat{f}_0}{\sqrt{|\mathbb{T}_b^d|}}$, we introduce the indicator function $\mathbf{1}_{h;2^{j_b}}(j)$.}
\begin{align}
\|T_{f}g\|^{h;\zeta}_{\dot{B}_{2,r}^s}&\lesssim  \Big\|\mathbf{1}_{h;2^{j_b}}(j)2^{j(s_1-\f d 2)}\|{S}_{j-2}f\|_{L^{\infty}}\Big \|_{\textit{l}^\infty}\|g\|^{h;\f\zeta 4}_{B_{2,r}^{s_2}}\nonumber\\
&\lesssim \|f\|_{B_{\infty,1}^{s_1-\f d 2}}\|g\|^{h;\f\zeta 4}_{B_{2,r}^{s_2}}\lesssim \|f\|_{B_{2,1}^{s_1}}\|g\|^{h;\f\zeta 4}_{B_{2,r}^{s_2}}\nonumber.
\end{align}
For the dyadic blocks $\Delta_j R(f,g)$ satisfying $\zeta\leq 2^j$, we find that
\begin{align*}
2^{js}\|\Delta_j R(f,g)\|_{L^2}&\lesssim 2^{j(s+\f d 2)}\|\Delta_j R(f,g)\|_{L^1} \lesssim 2^{j(s+\f d 2)}\sum\limits_{j-j'\leq4}\sum\limits_{|k-j'|\leq2}\|\Delta_{j'}g\|_{L^{2}}\|\Delta_k f\|_{L^{2}}\\
&\lesssim
\sum\limits_{j-j'\leq 4}2^{(s+\f d 2)(j-j')}\mathbf{1}_{h;\f\zeta {16}}(j')2^{j's_2}\|\Delta_{j'}g\|_{L^{2}}\sum\limits_{k\geq j-6}\mathbf{1}_{h;\f\zeta {64}}(k)2^{ks_1}\|\Delta_{k}f\|_{L^{2}}.
\end{align*}
Taking the $\textit{l}^r$-norm over these  dyadic blocks  and applying Young's inequality establishes \eqref{estpri3}. $\Box$

\begin{Lemma}\emph{(\cite{F})}.\label{lema4.2} Let $0<T\leq \infty$. Assume system \eqref{equ2} admits a  solution
\[
(\vae,\vue) \text{ satisfying  \eqref{new3} } \text{ with } \inf_{(t,x)\in (0,T)\times \mathbb{R}^d}1+\ep\vae(t,x)>0.
\]
Then there exist a positive constants $\eta_0= \eta_0(d,\nu)$ such that  for any $\ep\leq1$ and $t\in(0,T]$, the following estimates hold,
\begin{align}
   \ep&\|\vae\|^{h;\f{\eta_0}\ep}_{\widetilde{L}_t^\infty(B^{\f d 2}_{2,1})}+\f1 \ep \|\vae\|^{h;\f{\eta_0}\ep}_{L_t^1(B^{\f d 2}_{2,1})}
  +\|\mathbb{Q}\vue\|^{h;\f{\eta_0}\ep}_{\widetilde{L}_t^\infty(B^{\f d 2-1}_{2,1})\cap L_t^1(B^{\f d 2+1}_{2,1})}\nonumber\\
  &\lesssim\left(\ep\|\vae\|^{h;\f{\eta_0}\ep}_{B_{2,1}^{\f d 2}}(0)+
  \|\mathbb{Q}\vue\|^{h;\f{\eta_0}\ep}_{B_{2,1}^{\f d 2-1}}(0)\right)+\|(F,\mathbb{Q}G,\mathbb{Q}(\vue\cdot\Grad\vue)\|^{h;\f{\eta_0}\ep}_{L_t^1(B_{2,1}^{\f d 2-1})}\label{estpri5}\\
  &\quad+\ep\sum\limits_{2^j\geq\f{\eta_0}\ep}2^{j\f d 2}\left(\Big\|\|\Grad{S}_{j-2}\vue\|_{L^\infty}\|\De_j\vae\|_{L^2}\Big\|_{L^1(0,t)}+\|\Delta_j F+{S}_{j-2}\vue\Grad \De_j \vae\|_{L^1(0,t;L^2)}\right),\nonumber
\end{align}
where
\begin{align*}
F=-\Div(\vae\vue),\ \  G=-\kappa\vae\Grad\vae-\widetilde{K}(\ep\vae)\vae\Grad\vae-I(\ep\vae)\mathcal{A}\vue+f^\ep.
\end{align*}

\end{Lemma}

We emphasize that although Lemma \ref{lema4.2} was originally proved for the whole space $\mathbb{R}^d$, its conclusions extend verbatim to periodic domains without modification.

\begin{Lemma}\emph{(\cite{D1})}.\label{lema4.3} Let $0<T\leq\infty$. Assume system \eqref{equ2} admits a  solution
\[
(\vae,\vue) \text{ satisfying \eqref{new3} }   \text{ with } \inf_{(t,x)\in (0,T)\times \mathbb{R}^d}1+\ep\vae(t,x)>0.
\]
Then for any $0\leq\zeta<\f{\eta_0}\ep$ and $t\in(0,T]$, we have that
\begin{align}
&\|(\vae,\mathbb{Q}\vue)\|^{m;\zeta,\f{\eta_0}\ep}_{\widetilde{L}_t^\infty(B_{2,1}^{\f d 2-1})}+ \|(\vae,\mathbb{Q}\vue)\|^{m;\zeta,\f{\eta_0}\ep}_{L_t^1(B_{2,1}^{\f d 2+1})}\nonumber\\
&\lesssim\|(\vae,\mathbb{Q}\vue)(0)\|^{m;\zeta,\f{\eta_0}\ep}_{B_{2,1}^{\f d 2-1}}
+\|\mathbb{Q}G\|^{m;\zeta,\f{\eta_0}\ep}_{L_t^1(B_{2,1}^{\f d 2-1})}
+\sum\limits_{\zeta\leq2^j<\f{\eta_0}\ep}2^{j(\f d 2-1)}\Bigg(\Big\|\|\Grad S_{j-2}\vue\|_{L^\infty}\|\De_j\vae\|_{L^2}\Big\|_{L^1(0,t)}\nonumber\\
&\quad+\Big\|\|\Grad S_{j-2}\vue\|_{L^\infty}\|\De_j\mathbb{Q}\vue\|_{L^2}\Big\|_{L^1(0,t)}+\|\Delta_j F+{S}_{j-2}\vue\Grad \De_j \vae\|_{L^1(0,t;L^2)}\nonumber\\
&\quad+\|\Delta_j\mathbb{Q}(\vue\cdot\Grad\vue)-S_{j-2}\vue\cdot\Grad\De_j\mathbb{Q}\vue\|_{L^1(0,t;L^2)}\Bigg),
\end{align}
where
\begin{align*}
F=-\Div(\vae\vue),\ \  G=-\kappa\vae\Grad\vae-\widetilde{K}(\ep\vae)\vae\Grad\vae-I(\ep\vae)\mathcal{A}\vue+f^\ep.
\end{align*}
\end{Lemma}
\subsection{The proof of Proposition \ref{pro3}} \label{newsec3}
From lemma \ref{lema4.2} and Lemma \ref{lema4.3} (set $\zeta=0$), we have that
\begin{align}\label{estpri7}
X^\ep(T)\lesssim& X_0+\|\mathbb{Q}G\|_{L_T^1(B_{2,1}^{\f d 2-1})}+\|(F,\mathbb{Q}(\vue\cdot\Grad\vue)\|^{h;\f{\eta_0}\ep}_{L_T^1(B_{2,1}^{\f d 2-1})}+\ep\|(T'_{\vue}\Grad\vae,\vae\Div\vue)\|^{h;\f{\eta_0}\ep}_{L_T^1(B_{2,1}^{\f d 2})}\nonumber\\
&+\|(T'_{\vue}\Grad\vae,\vae\Div\vue,\mathbb{Q}(T'_{\vue}\Grad\vue))\|^{l;\f{\eta_0}\ep}_{L_T^1(\underline{B}_{2,1}^{\f d 2-1})}+R_1+R_2,
\end{align}
where
\begin{align*}
  R_1=&\ep\sum\limits_{2^j\geq\f{\eta_0}\ep}2^{j\f d 2}\Big\|\|\Grad{S}_{j-2}\vue\|_{L^\infty}\|\De_j\vae\|_{L^2}\Big\|_{L^1(0,T)}\\
  &+\sum\limits_{2^j<\f{\eta_0}\ep}2^{j(\f d 2-1)}\Bigg(\Big\|\|\Grad S_{j-2}\vue\|_{L^\infty}\|\De_j\vae\|_{L^2}\Big\|_{L^1(0,T)}
  +\Big\|\|\Grad S_{j-2}\vue\|_{L^\infty}\|\De_j\mathbb{Q}\vue\|_{L^2}\Big\|_{L^1(0,T)}\Bigg)
\end{align*}
and
\begin{align*}
  R_2=&\ep\sum\limits_{2^j\geq\f{\eta_0}\ep}2^{j\f d 2}\|\De_j(T_{\vue}\Grad \vae)-S_{j-2}\vue\Grad\De_j\vae\|_{L^1(0,T;L^2)}+
  \sum\limits_{2^j<\f{\eta_0}\ep}2^{j(\f d 2-1)}\\
  &\Bigg(\|\De_j(T_{\vue}\Grad \vae)-S_{j-2}\vue\Grad\De_j\vae\|_{L^1(0,T;L^2)}
  +\|\De_j\mathbb{Q}(T_{\vue}\Grad \vue)-S_{j-2}\vue\Grad\De_j\mathbb{Q}\vue\|_{L^1(0,T;L^2)}\Bigg).
\end{align*}
Direct calculation shows that
\begin{align}\label{estpri8}
\|\mathbb{Q}G\|_{L_T^1(B_{2,1}^{\f d 2-1})}&\lesssim (1+\ep\|\vae\|_{L_T^\infty(B_{2,1}^{\f d 2})})\|\vae\|^2_{L_T^2(B_{2,1}^{\f d 2})}+\|\mathbb{Q}f\|_{L_T^1(B_{2,1}^{\f d 2-1})}\nonumber\\
&\quad+\ep\|\vae\|_{L_T^\infty(B_{2,1}^{\f d 2})}\|\vue\|^{h;\f{\eta_0}\ep}_{L_T^1(B_{2,1}^{\f d 2+1})}
+\|\vae\|_{L_T^2(B_{2,1}^{\f d 2})}\|\vue\|^{l;\f{\eta_0}\ep}_{L_T^2(\underline{B}_{2,1}^{\f d 2})}\nonumber\\
&\lesssim (Y^{\zeta,\ep}(T))^2+(Y^{\zeta,\ep}(T))^3+\mathbb{Q}_f,
\end{align}
\begin{align}\label{estpri19}
  \|F\|^{h;\f{\eta_0}\ep}_{L_T^1(B_{2,1}^{\f d 2-1})}&\lesssim  \|\Div(\vae\underline{\vue})\|^{h;\f{\eta_0}\ep}_{L_T^1(B_{2,1}^{\f d 2-1})}+\|\Div(\vae\f{\widehat{\vue}_0}{\sqrt{|\mathbb{T}_b^d|}})\|^{h;\f{\eta_0}\ep}_{L_T^1(B_{2,1}^{\f d 2-1})}\nonumber\\
  &\lesssim \|\vae\|_{L_T^2(B_{2,1}^{\f d 2})}\|\vue\|_{L_T^2(\underline{B}_{2,1}^{\f d 2})}+\ep\f{1}\ep\|\vae\|^{h;\f{\eta_0}\ep}_{L_T^1(B_{2,1}^{\f d 2})}\|\vue\|_{L_T^\infty(B_{2,1}^{\f d 2-1})}\nonumber\\
   &\lesssim (Y^{\zeta,\ep}(T))^2,
\end{align}
\begin{align}\label{estpri20}
 \|\mathbb{Q}(\vue\cdot\Grad\vue)\|^{h;\f{\eta_0}\ep}_{L_T^1(B_{2,1}^{\f d 2-1})}&\lesssim \|\underline{\vue}\cdot\Grad\vue\|^{h;\f{\eta_0}\ep}_{L_T^1(B_{2,1}^{\f d 2-1})}+ \|\f{\widehat{\vue}_0}{\sqrt{|\mathbb{T}_b^d|}}\cdot\Grad\vue\|^{h;\f{\eta_0}\ep}_{L_T^1(B_{2,1}^{\f d 2-1})}\nonumber\\
 &\lesssim \|\vue\|^2_{L_T^2(\underline{B}_{2,1}^{\f d 2})}+\ep\|\vue\|^{h;\f{\eta_0}\ep}_{L_T^1(B_{2,1}^{\f d 2+1})}\|\vue\|_{L_T^\infty(B_{2,1}^{\f d 2-1})}\nonumber\\
   &\lesssim (Y^{\zeta,\ep}(T))^2,
\end{align}
\begin{align}\label{estpri21}
 \ep\|(T'_{\vue}\Grad\vae,\vae\Div\vue)\|_{L_T^1(\underline{B}_{2,1}^{\f d 2})}&\lesssim \ep\Big\|\|\vae\|_{B_{2,1}^{\f d 2}}\|\vue\|_{\underline{B}_{2,1}^{\f d 2+1}}\Big\|_{L^1(0,T)}\nonumber\\
 &\lesssim\ep\|\vae\|_{L_T^\infty(B_{2,1}^{\f d 2})}\|\vue\|^{h;\f{\eta_0}\ep}_{L_T^1(B_{2,1}^{\f d 2+1})}
+\|\vae\|_{L_T^2(B_{2,1}^{\f d 2})}\|\vue\|^{l;\f{\eta_0}\ep}_{L_T^2(\underline{B}_{2,1}^{\f d 2})}\nonumber\\
   &\lesssim (Y^{\zeta,\ep}(T))^2,
\end{align}
\begin{align}\label{estpri22}
\|(T'_{\vue}\Grad\vae,\vae\Div\vue,\mathbb{Q}(T'_{\vue}\Grad\vue))\|_{L_T^1(\underline{B}_{2,1}^{\f d 2-1})}&\lesssim (\|\vae\|_{L_T^2(B_{2,1}^{\f d 2})}+\|\vue\|_{L_T^2(\underline{B}_{2,1}^{\f d 2})})\|\vue\|_{L_T^2(\underline{B}_{2,1}^{\f d 2})}\nonumber\\
&\lesssim  (Y^{\zeta,\ep}(T))^2,
\end{align}
\begin{align}\label{estpri23}
 R_1\lesssim \ep\|\vae\|_{L_T^\infty(B_{2,1}^{\f d 2})}\|\vue\|^{h;\f{\eta_0}\ep}_{L_T^1(B_{2,1}^{\f d 2+1})}
+\|\vae\|_{L_T^2(B_{2,1}^{\f d 2})}\|\vue\|_{L_T^2(\underline{B}_{2,1}^{\f d 2})}\lesssim  (Y^{\zeta,\ep}(T))^2.
\end{align}

According to Lemma 7.5 in \cite{D3}, we get that
\begin{align}\label{estpri25}
  \|\De_j(T_{\vue}\Grad \vae)-S_{j-2}\vue\Grad\De_j\vae\|_{L^2}\lesssim \sum\limits_{|j-j'|\leq2}&\Big(\|S_{j'-2}\vue-S_{j-2}\vue\|_{L^\infty}2^{j'}\|\De_{j'}\vae\|_{L^2}\nonumber\\
&\quad+\|\Grad S_{j'-2}\vue\|_{L^\infty}\|\De_{j'}\vae\|_{L^2}\Big)
\end{align}
\begin{align}\label{estpri26}
\|\De_j\mathbb{Q}(T_{\vue}\Grad \vue)-S_{j-2}\vue\Grad\De_j\mathbb{Q}\vue\|_{L^2}\lesssim \sum\limits_{|j-j'|\leq2}&\Big(\|S_{j'-2}\vue-S_{j-2}\vue\|_{L^\infty}2^{j'}\|\De_{j'}\vue\|_{L^2}\nonumber\\
&\quad+\|\Grad S_{j'-2}\vue\|_{L^\infty}\|\De_{j'}\vue\|_{L^2}\Big)
\end{align}
It follows that
\begin{align}\label{estpri24}
R_2&\lesssim \ep\|\vae\|_{L_T^\infty(B_{2,1}^{\f d 2})}\|\vue\|^{h;\f{\eta_0}\ep}_{L_T^1(B_{2,1}^{\f d 2+1})}
+(\|\vae\|_{L_T^2(B_{2,1}^{\f d 2})}+\|\vue\|_{L_T^2(\underline{B}_{2,1}^{\f d 2})})\|\vue\|_{L_T^2(\underline{B}_{2,1}^{\f d 2})}\nonumber\\
&\lesssim  (Y^{\zeta,\ep}(T))^2.
\end{align}

Substituting \eqref{estpri8}-\eqref{estpri23} and \eqref{estpri24} into \eqref{estpri7}, we obtain that
\begin{align}\label{estpri27}
 X^\ep(T)\lesssim X_0+\mathbb{Q}f+(Y^{\zeta,\ep}(T))^2+(Y^{\zeta,\ep}(T))^3.
\end{align}

We now turn to the analysis of $\mathbb{P}\vue$. Applying operator $\mathbb{P}$ to system \eqref{equ2}$_2$, we find that
\begin{align*}
\p_t\mathbb{P}\vue+\mathbb{P}(\vue\cdot\Grad\vue)-\mu\De\mathbb{P}\vue=\mathbb{P}f-\mathbb{P}(I(\ep\vae)\mathcal{A}\vue).
\end{align*}

Employing the standard energy method, we derive that\footnote{Here we  use the property  $\widehat{\mathbb{P}f}_0=\widehat{f}_0$ which holds for any function $f$ defined on $\Td$.}
\begin{align}\label{estpri31}
P^\ep(T)&\lesssim P_0+\|(\mathbb{P}f,\mathbb{P}(I(\ep\vae)\mathcal{A}\vue)\|_{L_T^1(B_{2,1}^{\f d 2-1})}+\|\mathbb{P}(T'_{\vue}\Grad \vue)\|_{L_T^1(\underline{B}_{2,1}^{\f d 2-1})}+\Big\|\widehat{(\vue\cdot\Grad\vue)}_0\Big\|_{L^1(0,T)}\nonumber\\
&\quad+\sum\limits_{j\in\mathbb{Z}}2^{j(\f d 2-1)}
\Bigg(\Big\|\|\Grad S_{j-2}\vue\|_{L^\infty}\|\De_j\mathbb{P}\vue\|_{L^2}\Big\|_{L^1(0,T)}\nonumber\\
&\quad\qquad\qquad\qquad\qquad+\|\Delta_j\mathbb{P}(T_{\vue}\Grad \vue)-S_{j-2}\vue\cdot\Grad\De_j\mathbb{P}\vue\|_{L^1(0,T;L^2)}\Bigg).
\end{align}
Note that
\[
\widehat{(\vue\cdot\Grad\vue)}_0\lesssim \sum\limits_{k\neq0}|k||\widehat{\vue}_{k}||\widehat{\vue}_{-k}|\lesssim \|\vue\|_{\underline{H}^0}\|\vue\|_{\underline{H}^1}\lesssim \|\vue\|^2_{\underline{B}^{\f d 2}_{2,1}}.
\]
It follows that
\begin{align}\label{estpri29}
  \Big\|\widehat{(\vue\cdot\Grad\vue)}_0\Big\|_{L^1(0,T)}\lesssim \|\vue\|^2_{L_T^2(\underline{B}_{2,1}^{\f d 2})}.
\end{align}

The remaining computations follow analogously to the previous analysis. Thus, we obtain the following estimate
\begin{align}\label{estpri28}
 P^\ep(T)\lesssim P_0+\mathbb{P}f+(Y^{\zeta,\ep}(T))^2.
\end{align}

Proposition \ref{pro3} follows from \eqref{estpri27} and \eqref{estpri28}.\ \ $\Box$

\subsection{The proof of Proposition \ref{pro4}} \label{newsec4}
We first estimate the incompressible component. Analogous to \eqref{estpri31}, we derive that
\begin{align}\label{estpri32}
\|\mathbb{P}\vue\|^{h;\zeta}_{\widetilde{L}^\infty_t(B^{\f d 2-1}_{2,1})\cap L^1_t(B^{\f d 2+1}_{2,1}) }&\lesssim \|\mathbb{P}\vu_0\|^{h;\zeta}_{B^{\f d 2-1}_{2,1}}+R_1(t)\nonumber\\
&\quad+\|(\mathbb{P}f,\mathbb{P}(I(\ep\vae)\mathcal{A}\vue),\mathbb{P}(T'_{\vue}\Grad \vue))\|^{h;\zeta}_{L_t^1(B_{2,1}^{\f d 2-1})}
\end{align}
where
\begin{align*}
R_1(t)=\sum\limits_{2^j\geq \zeta}2^{j(\f d 2-1)}\Bigg(&\Big\|\|\Grad S_{j-2}\vue\|_{L^\infty}\|\De_j\mathbb{P}\vue\|_{L^2}\Big\|_{L^1(0,t)}\\
&+\|\Delta_j\mathbb{P}(T_{\vue}\Grad \vue)-S_{j-2}\vue\cdot\Grad\De_j\mathbb{P}\vue\|_{L^1(0,t;L^2)}\Bigg).
\end{align*}

Substituting $\mathbb{Q}$ for $\mathbb{P}$ in  \eqref{estpri26} yields
\begin{align}\label{estpri33}
\|\De_j\mathbb{P}(T_{\vue}\Grad \vue)-S_{j-2}\vue\Grad\De_j\mathbb{P}\vue\|_{L^2}\lesssim \sum\limits_{|j-j'|\leq2}&\Big(\|S_{j'-2}\vue-S_{j-2}\vue\|_{L^\infty}2^{j'}\|\De_{j'}\vue\|_{L^2}\nonumber\\
&\quad+\|\Grad S_{j'-2}\vue\|_{L^\infty}\|\De_{j'}\vue\|_{L^2}\Big).
\end{align}

It follows that

\begin{align*}
R_1(t)&\lesssim \int_0^t \|\vue(\tau)\|_{\underline{B}_{2,1}^{\f d 2}}\|\vue(\tau)\|^{h;\f \zeta 4}_{B_{2,1}^{\f d 2}}d\tau\nonumber\\
&\lesssim \int_0^t \|\vue(\tau)\|_{\underline{B}_{2,1}^{\f d 2}}\|(\vVe-\vV,\mathbb{P}\vue-\vv)(\tau)\|^{m;\f \zeta 4,\zeta}_{B_{2,1}^{\f d 2}}d\tau+\int_0^t \|\vue(\tau)\|_{\underline{B}_{2,1}^{\f d 2}}\|(\vv,\vV)(\tau)\|^{m;\f \zeta 4,\zeta}_{B_{2,1}^{\f d 2}}d\tau\nonumber\\
&\quad+\int_0^t\|\vue(\tau)\|_{\underline{B}_{2,1}^{\f d 2}}\Big(\|\vue(\tau)\|^{h;\zeta }_{B_{2,1}^{\f d 2-1}}\Big)^{\f1 2}\Big(\|\vue(\tau)\|^{h;\zeta }_{B_{2,1}^{\f d 2+1}}\Big)^{\f1 2}d\tau.
\end{align*}
Note that
\begin{align}\label{estpri34}
 C&\int_0^t\|\vue(\tau)\|_{\underline{B}_{2,1}^{\f d 2}}\Big(\|\vue(\tau)\|^{h;\zeta }_{B_{2,1}^{\f d 2-1}}\Big)^{\f1 2}\Big(\|\vue(\tau)\|^{h;\zeta }_{B_{2,1}^{\f d 2+1}}\Big)^{\f1 2}d\tau\nonumber\\
 &\leq 4C^2\int_0^t\|\vue(\tau)\|^2_{\underline{B}_{2,1}^{\f d 2}}\|\vue(\tau)\|^{h;\zeta }_{B_{2,1}^{\f d 2-1}}d\tau+\f 1 {16}\int_0^t\|\vue(\tau)\|^{h;\zeta }_{B_{2,1}^{\f d 2+1}}d\tau
\end{align}
and the term $\f 1 {16}\int_0^t\|\vue(\tau)\|^{h;\zeta }_{B_{2,1}^{\f d 2+1}}d\tau$ can ultimately be absorbed by the terms on the
left-hand side of \eqref{estpri32} and \eqref{estpri38}. Therefore, we can discard this term and obtain that
\begin{align}\label{estpri54}
 R_1(t)&\lesssim \int_0^t \|\vue(\tau)\|^2_{\underline{B}_{2,1}^{\f d 2}}[\vae,\vue]^{\zeta,\ep}_{h,m}(\tau)d\tau+\int_0^t \|\vue(\tau)\|_{\underline{B}_{2,1}^{\f d 2}}\|(\vv,\vV)(\tau)\|^{h;\f \zeta 4}_{B_{2,1}^{\f d 2}}d\tau\nonumber\\
 &\quad+\int_0^t \|\vue(\tau)\|_{\underline{B}_{2,1}^{\f d 2}}\|(\vVe-\vV,\mathbb{P}\vue-\vv)(\tau)\|^{l;\zeta}_{\underline{B}_{2,1}^{\f d 2}}d\tau
\end{align}
Straightforward computation yields
\begin{align}\label{estpri35}
 \|\mathbb{P}(I(\ep\vae)\mathcal{A}\vue)|_{L_t^1(B_{2,1}^{\f d 2-1})}&\lesssim \int_0^t \ep\|\vae(\tau)\|_{B_{2,1}^{\f d 2}}\|\vue(\tau)\|_{\underline{B}_{2,1}^{\f d 2+1}}d\tau\\
 &\lesssim \int_0^t \|\vue(\tau)\|^{h;\zeta}_{B_{2,1}^{\f d 2+1}}[\vae,\vue]^{\zeta,\ep}_{h,m}(\tau)d\tau+(\ep\zeta)\int_0^t \|\vae(\tau)\|_{B_{2,1}^{\f d 2}}\|\vue(\tau)\|^{l;\zeta}_{\underline{B}_{2,1}^{\f d 2}}d\tau\nonumber
\end{align}
Using Lemma \ref{lema4.1} and following a calculation similar to that in \eqref{estpri54}, we get that
\begin{align}\label{estpri36}
 \|\mathbb{P}(T'_{\vue}\Grad \vue)\|^{h;\zeta}_{L_t^1(B_{2,1}^{\f d 2-1})}&\lesssim \int_0^t \|\vue(\tau)\|_{\underline{B}_{2,1}^{\f d 2}}\|\vue(\tau)\|^{h;\f \zeta {16}}_{B_{2,1}^{\f d 2}}d\tau\nonumber\\
 &\lesssim\int_0^t \|\vue(\tau)\|^2_{\underline{B}_{2,1}^{\f d 2}}[\vae,\vue]^{\zeta,\ep}_{h,m}(\tau)d\tau+\int_0^t \|\vue(\tau)\|_{\underline{B}_{2,1}^{\f d 2}}\|(\vv,\vV)(\tau)\|^{h;\f \zeta {16}}_{B_{2,1}^{\f d 2}}d\tau\nonumber\\
 &\quad+\int_0^t \|\vue(\tau)\|_{\underline{B}_{2,1}^{\f d 2}}\|(\vVe-\vV,\mathbb{P}\vue-\vv)(\tau)\|^{l;\zeta}_{\underline{B}_{2,1}^{\f d 2}}d\tau
\end{align}

Substituting \eqref{estpri54}-\eqref{estpri36}  into \eqref{estpri32} yields
\begin{align}\label{estpri37}
 \|\mathbb{P}\vue\|^{h;\zeta}_{\widetilde{L}^\infty_t(B^{\f d 2-1}_{2,1})\cap L^1_t(B^{\f d 2+1}_{2,1}) }&\lesssim \|\mathbb{P}\vu_0\|^{h;\zeta}_{B^{\f d 2-1}_{2,1}}+\|\mathbb{P}f\|^{h;\zeta}_ {L^1_t(B^{\f d 2-1}_{2,1}) }+\int_0^t \|\vue(\tau)\|_{\underline{B}_{2,1}^{\f d 2}}\|(\vv,\vV)(\tau)\|^{h;\f \zeta {16}}_{B_{2,1}^{\f d 2}}d\tau\nonumber\\
 &\quad+\int_0^t \Big(\|\vue(\tau)\|^{h;\zeta}_{B_{2,1}^{\f d 2+1}}+\|\vue(\tau)\|^2_{\underline{B}_{2,1}^{\f d 2}}\Big)[\vae,\vue]^{\zeta,\ep}_{h,m}(\tau)d\tau\nonumber\\
 &\quad+\int_0^t \|\vue(\tau)\|_{\underline{B}_{2,1}^{\f d 2}}\|(\vVe-\vV,\mathbb{P}\vue-\vv)(\tau)\|^{l;\zeta}_{\underline{B}_{2,1}^{\f d 2}}d\tau\nonumber\\
 &\quad +(\ep\zeta)\int_0^t \|\vae(\tau)\|_{B_{2,1}^{\f d 2}}\|\vue(\tau)\|_{\underline{B}_{2,1}^{\f d 2}}d\tau
\end{align}

Next we  estimate the compressible component. From lemma \ref{lema4.2} and Lemma \ref{lema4.3}, we get that
\begin{align}\label{estpri38}
    &\ep\|\vae\|^{h;\f{\eta_0}\ep}_{\widetilde{L}_t^\infty(B^{\f d 2}_{2,1})}+\f1 \ep \|\vae\|^{h;\f{\eta_0}\ep}_{L_t^1(B^{\f d 2}_{2,1})}
  +\|\vae\|^{m;\zeta,\f{\eta_0}\ep}_{\widetilde{L}_t^\infty(B^{\f d 2-1}_{2,1})\cap L_t^1(B^{\f d 2+1}_{2,1})}+\|\mathbb{Q}\vue\|^{h;\zeta}_{\widetilde{L}_t^\infty(B^{\f d 2-1}_{2,1})\cap L_t^1(B^{\f d 2+1}_{2,1})}\\
  &\lesssim \left(\ep\|\vae\|^{h;\f{\eta_0}\ep}_{B_{2,1}^{\f d 2}}(0)+
  \|(\vae,\mathbb{Q}\vue)\|^{h;\zeta}_{B_{2,1}^{\f d 2-1}}(0)\right)+\ep\|(T'_{\vue}\Grad\vae,\vae\Div\vue)\|^{h;\f{\eta_0}\ep}_{L_t^1(B_{2,1}^{\f d 2})}\nonumber\\
  &\quad+\|(\mathbb{Q}G, \Div\vue\vae, T_{\vue}\Grad\vae, T'_{\vue}\Grad\vae,\mathbb{Q}(T'_{\vue}\Grad\vue), \mathbb{Q}(T_{\vue}\Grad\vue))\|^{h;\zeta}_{L_t^1(B_{2,1}^{\f d 2-1})}+R_2(t)+R_3(t),\nonumber
\end{align}
where
\[
R_2(t)=\ep\sum\limits_{2^j\geq\f{\eta_0}\ep}2^{j\f d 2}\left(\Big\|\|\Grad{S}_{j-2}\vue\|_{L^\infty}\|\De_j\vae\|_{L^2}\Big\|_{L^1(0,t)}+\|\De_j(T_{\vue}\Grad \vae)-{S}_{j-2}\vue\Grad \De_j \vae\|_{L^1(0,t;L^2)}\right)
\]
\begin{align*}
 R_3(t)&=\sum\limits_{\zeta\leq2^j<\f{\eta_0}\ep}2^{j(\f d 2-1)}\Bigg(\Big\|\|\Grad S_{j-2}\vue\|_{L^\infty}\|\De_j\vae\|_{L^2}\Big\|_{L^1(0,t)}+\Big\|\|\Grad S_{j-2}\vue\|_{L^\infty}\|\De_j\mathbb{Q}\vue\|_{L^2}\Big\|_{L^1(0,t)}\\
& +\|\Delta_j (T_{\vue}\Grad \vae)-{S}_{j-2}\vue\Grad \De_j \vae\|_{L^1(0,t;L^2)}+\|\Delta_j\mathbb{Q}(T_{\vue}\Grad \vue)-S_{j-2}\vue\cdot\Grad\De_j\mathbb{Q}\vue\|_{L^1(0,t;L^2)}\Bigg).
\end{align*}
Direct calculation shows that
\begin{align}\label{estpri39}
 \ep\|(T'_{\vue}\Grad\vae,\vae\Div\vue)\|_{L_t^1(\underline{B}_{2,1}^{\f d 2})}&\lesssim \int_0^t \ep\|\vae(\tau)\|_{B_{2,1}^{\f d 2}}\|\vue(\tau)\|_{\underline{B}_{2,1}^{\f d 2+1}}d\tau\\
 &\lesssim \int_0^t \|\vue(\tau)\|^{h;\zeta}_{B_{2,1}^{\f d 2+1}}[\vae,\vue]^{\zeta,\ep}_{h,m}(\tau)d\tau+(\ep\zeta)\int_0^t \|\vae(\tau)\|_{B_{2,1}^{\f d 2}}\|\vue(\tau)\|^{l;\zeta}_{\underline{B}_{2,1}^{\f d 2}}d\tau\nonumber
\end{align}
From Lemma \ref{lema4.1}, we infer that
\begin{align}
 \|\vae\Grad\vae\|^{h;\zeta}_{L_t^1(B_{2,1}^{\f d 2-1})}&\lesssim \int_0^t\|\vae(\tau)\|_{B_{2,1}^{\f d 2}}\|\vae(\tau)\|^{l;\f \zeta {16}}_{B_{2,1}^{\f d 2}}d\tau\nonumber\\
 &\lesssim \int_0^t\|\vae(\tau)\|_{B_{2,1}^{\f d 2}}\Bigg(\|(\vVe-\vV)(\tau)\|^{m;\f \zeta {16},\zeta}_{B_{2,1}^{\f d 2}}+\|\vV(\tau)\|^{m;\f \zeta {16},\zeta}_{B_{2,1}^{\f d 2}}\nonumber\\
 &\quad+\Big(\|\vae(\tau)\|^{m;\zeta,\f{\eta_0}\ep}_{B_{2,1}^{\f d 2-1}}\Big)^{\f 1 2}\Big(\|\vae(\tau)\|^{m;\zeta,\f{\eta_0}\ep}_{B_{2,1}^{\f d 2+1}}\Big)^{\f 1 2}+\|\vae(\tau)\|^{h;\f{\eta_0}\ep}_{B_{2,1}^{\f d 2}}\Bigg)d\tau\nonumber
\end{align}
Note that
\begin{align}\label{estpri46}
  &C\int_0^t\|\vae(\tau)\|_{B_{2,1}^{\f d 2}}\Big(\|\vae(\tau)\|^{m;\zeta,\f{\eta_0}\ep}_{B_{2,1}^{\f d 2-1}}\Big)^{\f 1 2}\Big(\|\vae(\tau)\|^{m;\zeta,\f{\eta_0}\ep}_{B_{2,1}^{\f d 2+1}}\Big)^{\f 1 2}d\tau\nonumber\\
&\quad\leq 4C^2\int_0^t\|\vae(\tau)\|^2_{B_{2,1}^{\f d 2}}\|\vae(\tau)\|^{m;\zeta,\f{\eta_0}\ep}_{B_{2,1}^{\f d 2-1}}d\tau+\f 1 {16}\int_0^t\|\vae(\tau)\|^{m;\zeta,\f{\eta_0}\ep}_{B_{2,1}^{\f d 2+1}}d\tau
\end{align}
and the term $\f 1 {16}\int_0^t\|\vae(\tau)\|^{m;\zeta,\f{\eta_0}\ep}_{B_{2,1}^{\f d 2+1}}d\tau$ can ultimately be absorbed by the term on the left-hand side of \eqref{estpri38}. Therefore, we can discard this term and obtain that
\begin{align}\label{estpri47}
  \|\vae\Grad\vae\|^{h;\zeta}_{L_t^1(B_{2,1}^{\f d 2-1})}&\lesssim \int_0^t \Big(\|\vae(\tau)\|^2_{B_{2,1}^{\f d 2}}+\f1 \ep\|\vae(\tau)\|^{h;\f{\eta_0}\ep}_{B_{2,1}^{\f d 2}}\Big)[\vae,\vue]^{\zeta,\ep}_{h,m}(\tau)d\tau\nonumber\\
  &\quad+\int_0^t\|\vae(\tau)\|_{B_{2,1}^{\f d 2}}\Big(\|(\vVe-\vV)(\tau)\|^{l;\zeta}_{B_{2,1}^{\f d 2}}+\|\vV(\tau)\|^{h;\f \zeta {16}}_{B_{2,1}^{\f d 2}}\Big)d\tau
\end{align}

We also have that
\begin{align}\label{estpri41}
   \|(\widetilde{K}(\ep\vae)\vae\Grad\vae,I(\ep\vae)\mathcal{A}\vue)\|_{L_t^1(B_{2,1}^{\f d 2-1})}&\lesssim \int_0^t \ep\|\vae(\tau)\|_{B_{2,1}^{\f d 2}}(\|\vae(\tau)\|^2_{B_{2,1}^{\f d 2}}+\|\vue(\tau)\|_{\underline{B}_{2,1}^{\f d 2+1}})d\tau\nonumber\\
    &\lesssim \int_0^t (\|\vae(\tau)\|^2_{B_{2,1}^{\f d 2}}+\|\vue(\tau)\|^{h;\zeta}_{\underline{B}_{2,1}^{\f d 2+1}})[\vae,\vue]^{\zeta,\ep}_{h,m}(\tau)d\tau\nonumber\\
     &\quad+(\ep\zeta)\int_0^t \|\vae(\tau)\|_{B_{2,1}^{\f d 2}}\|\vue(\tau)\|^{l;\zeta}_{\underline{B}_{2,1}^{\f d 2}}d
\end{align}
Then \eqref{estpri47} and \eqref{estpri41} yield that
\begin{align}\label{estpri42}
\|\mathbb{Q}G\|^{h;\zeta}_{L_t^1(B_{2,1}^{\f d 2-1})}&\lesssim \int_0^t \Big(\|\vae(\tau)\|^2_{B_{2,1}^{\f d 2}}+\|\vue(\tau)\|^{h;\zeta}_{\underline{B}_{2,1}^{\f d 2+1}}+\f 1 {\ep}\|\vae(\tau)\|^{h;\f {\eta_0} {\ep}}_{B_{2,1}^{\f d 2}}\Big)[\vae,\vue]^{\zeta,\ep}_{h,m}(\tau)d\tau+\|\mathbb{Q}f\|^{h;\zeta}_{L_t^1(B_{2,1}^{\f d 2-1})}\nonumber\\
&\quad+\int_0^t\|\vae(\tau)\|_{B_{2,1}^{\f d 2}}\|(\vVe-\vV)(\tau)\|^{l;\zeta}_{B_{2,1}^{\f d 2}}d\tau+\int_0^t\|\vae(\tau)\|_{B_{2,1}^{\f d 2}}\|\vV(\tau)\|^{h;\f \zeta {16}}_{B_{2,1}^{\f d 2}}d\tau\nonumber\\
&\quad+(\ep\zeta)\int_0^t \|\vae(\tau)\|_{B_{2,1}^{\f d 2}}\|\vue(\tau)\|^{l;\zeta}_{\underline{B}_{2,1}^{\f d 2}}d\tau.
\end{align}

Note that
\[
T_{\vue}\Grad\vae=T_{\underline{\vue}}\Grad\vae+\f{\widehat{\vue}_0}{\sqrt{|\mathbb{T}_b^d|}}\Grad\vae \ \ \text{ and } \ \ T_{\vue}\Grad\vue=T_{\underline{\vue}}\Grad\vue+\f{\widehat{\vue}_0}{\sqrt{|\mathbb{T}_b^d|}}\Grad\vue.
\]

It is clear that
\begin{align}\label{estpri48}
\|\f{\widehat{\vue}_0}{\sqrt{|\mathbb{T}_b^d|}}\Grad\vae\|^{h;\zeta}_{L_t^1(B_{2,1}^{\f d 2-1})}&\lesssim \int_0^t \|\vue(\tau)\|_{B_{2,1}^{\f d 2-1}}
\|\vae(\tau)\|^{h;\zeta }_{B_{2,1}^{\f d 2}}d\tau\\
&\lesssim \f{1}{\zeta} \int_0^t \|\vue(\tau)\|_{B_{2,1}^{\f d 2-1}}
\|\vae(\tau)\|^{m;\zeta,\f{\eta_0}\ep }_{B_{2,1}^{\f d 2+1}}d\tau+\ep \int_0^t \|\vue(\tau)\|_{B_{2,1}^{\f d 2-1}}
\f1 \ep\|\vae(\tau)\|^{h;\f{\eta_0}\ep }_{B_{2,1}^{\f d 2}}d\tau\nonumber
\end{align}
and
\begin{align}\label{estpri49}
\|\mathbb{Q}\Big(\f{\widehat{\vue}_0}{\sqrt{|\mathbb{T}_b^d|}}\Grad\vue\Big)\|^{h;\zeta}_{L_t^1(B_{2,1}^{\f d 2-1})}\lesssim \f{1}{\zeta} \int_0^t \|\vue(\tau)\|_{B_{2,1}^{\f d 2-1}}
\|\vue(\tau)\|^{h;\zeta}_{B_{2,1}^{\f d 2+1}}d\tau.
\end{align}
Similar to \eqref{estpri42}, we infer that

\begin{align}\label{estpri50}
&\|T_{\underline{\vue}}\Grad\vae\|^{h;\zeta}_{L_t^1(B_{2,1}^{\f d 2-1})}\nonumber\\
&\quad\lesssim \int_0^t \|\vue(\tau)\|_{\underline{B}_{2,1}^{\f d 2}}\|\vae(\tau)\|^{h;\f {\zeta} {16}}_{B_{2,1}^{\f d 2}}d\tau\nonumber\\
&\quad\lesssim\int_0^t\|\vue(\tau)\|_{\underline{B}_{2,1}^{\f d 2}}\Big(\|(\vVe-\vV)(\tau)\|^{m;\f {\zeta} {16},\zeta}_{B_{2,1}^{\f d 2}}+\|\vV(\tau)\|^{m;\f {\zeta} {16},\zeta}_{B_{2,1}^{\f d 2}}\Big)d\tau\nonumber\\
&\qquad+\int_0^t \|\vue(\tau)\|_{\underline{B}_{2,1}^{\f d 2}}\Big(\|\vae(\tau)\|^{m;\zeta,\f {\eta_0}\ep}_{B_{2,1}^{\f d 2-1}}\Big)^{\f 1 2}\Big(\|\vae(\tau)\|^{m;\zeta,\f {\eta_0}\ep}_{B_{2,1}^{\f d 2+1}}\Big)^{\f 1 2}d\tau\nonumber\\
&\qquad+(\ep\zeta)\int_0^t \f1 \ep\|\vae(\tau)\|^{h;\f {\eta_0} {\ep}}_{B_{2,1}^{\f d 2}}\|\vue(\tau)\|^{l;\zeta}_{\underline{B}_{2,1}^{\f d 2-1}}d\tau\\
&\qquad+\int_0^t \|\vae(\tau)\|^{h;\f {\eta_0} {\ep}}_{B_{2,1}^{\f d 2}}\Big(\|\vue(\tau)\|^{h;\zeta}_{B_{2,1}^{\f d 2-1}}\Big)^{\f 1 2}\Big(\|\vue(\tau)\|^{h;\zeta}_{B_{2,1}^{\f d 2+1}}\Big)^{\f 1 2}d\tau\nonumber\\
&\quad\lesssim \int_0^t \Big(\|\vae(\tau)\|^2_{B_{2,1}^{\f d 2}}+\|\vue(\tau)\|^2_{\underline{B}_{2,1}^{\f d 2}}\Big)[\vae,\vue]^{\zeta,\ep}_{h,m}(\tau)d\tau+\int_0^t\|\vue(\tau)\|_{\underline{B}_{2,1}^{\f d 2}}\|\vV(\tau)\|^{h;\f {\zeta} {16}}_{B_{2,1}^{\f d 2}}d\tau\nonumber\\
&\qquad+(\ep\zeta)\int_0^t \f1 \ep\|\vae(\tau)\|^{h;\f {\eta_0} {\ep}}_{B_{2,1}^{\f d 2}}\|\vue(\tau)\|_{B_{2,1}^{\f d 2-1}}d\tau+\int_0^t\|\vue(\tau)\|_{\underline{B}_{2,1}^{\f d 2}}\|(\vVe-\vV)(\tau)\|^{l;\zeta}_{B_{2,1}^{\f d 2}}d\tau\nonumber.
\end{align}

\begin{align}\label{estpri44}
&\|(\mathbb{Q}(T_{\underline{\vue}}\Grad\vue),\mathbb{Q}(T'_{\vue}\Grad\vue))\|^{h;\zeta}_{L_t^1(B_{2,1}^{\f d 2-1})}\nonumber\\
&\quad\lesssim \int_0^t \|\vue(\tau)\|_{\underline{B}_{2,1}^{\f d 2}}\|\vue(\tau)\|^{h;\f {\zeta} {16}}_{B_{2,1}^{\f d 2}}d\tau\nonumber\\
 &\quad \lesssim \int_0^t \|\vue(\tau)\|^2_{\underline{B}_{2,1}^{\f d 2}}[\vae,\vue]^{\zeta,\ep}_{h,m}(\tau)d\tau+\int_0^t\|\vue(\tau)\|_{\underline{B}_{2,1}^{\f d 2}}\|(\vV,\vv)(\tau)\|^{h;\f \zeta {16}}_{B_{2,1}^{\f d 2}}d\tau\nonumber\\
 &\qquad+\int_0^t\|\vue(\tau)\|_{\underline{B}_{2,1}^{\f d 2}}\|(\vVe-\vV,\mathbb{P}\vue-\vv)(\tau)\|^{l;\zeta}_{\underline{B}_{2,1}^{\f d 2}}d\tau,
\end{align}

\begin{align}\label{estpri43}
\|T'_{\vue}\Grad\vae\|^{h;\zeta}_{L_t^1(B_{2,1}^{\f d 2-1})}&\lesssim \int_0^t \|\vae(\tau)\|_{B_{2,1}^{\f d 2}}\|\vue(\tau)\|^{h;\f {\zeta} {16}}_{B_{2,1}^{\f d 2}}d\tau\nonumber\\
 &\lesssim \int_0^t \|\vae(\tau)\|^2_{B_{2,1}^{\f d 2}}[\vae,\vue]^{\zeta,\ep}_{h,m}(\tau)d\tau+\int_0^t\|\vae(\tau)\|_{B_{2,1}^{\f d 2}}\|(\vVe-\vV,\mathbb{P}\vue-\vv)(\tau)\|^{l;\zeta}_{\underline{B}_{2,1}^{\f d 2}}d\tau\nonumber\\
  &\quad+\int_0^t\|\vae(\tau)\|_{B_{2,1}^{\f d 2}}\|(\vV,\vv)(\tau)\|^{h;\f \zeta {16}}_{B_{2,1}^{\f d 2}}d\tau.
\end{align}

Using Lemma \ref{lema4.1}, \eqref{estpri43} and \eqref{estpri50}, we get that
\begin{align}\label{estpri56}
&\|\Div\vue\vae\|^{h;\zeta}_{L_t^1(B_{2,1}^{\f d 2-1})}\nonumber\\
&\quad\lesssim \int_0^t \|\vue(\tau)\|_{\underline{B}_{2,1}^{\f d 2}}\|\vae(\tau)\|^{h;\f {\zeta} {16}}_{B_{2,1}^{\f d 2}}d\tau+\int_0^t \|\vae(\tau)\|_{B_{2,1}^{\f d 2}}\|\vue(\tau)\|^{h;\f {\zeta} {16}}_{B_{2,1}^{\f d 2}}d\tau\nonumber\\
&\quad\lesssim \int_0^t \Big(\|\vae(\tau)\|^2_{B_{2,1}^{\f d 2}}+\|\vue(\tau)\|^2_{\underline{B}_{2,1}^{\f d 2}}\Big)[\vae,\vue]^{\zeta,\ep}_{h,m}(\tau)d\tau+(\ep\zeta)\int_0^t \f1 \ep\|\vae(\tau)\|^{h;\f {\eta_0} {\ep}}_{B_{2,1}^{\f d 2}}\|\vue(\tau)\|_{B_{2,1}^{\f d 2-1}}d\tau\nonumber\\
&\qquad+\int_0^t\Big(\|\vue(\tau)\|_{\underline{B}_{2,1}^{\f d 2}}+\|\vae(\tau)\|_{B_{2,1}^{\f d 2}}\Big)\|(\vV,\vv)(\tau)\|^{h;\f {\zeta} {16}}_{B_{2,1}^{\f d 2}}d\tau\\
&\qquad+\int_0^t\Big(\|\vue(\tau)\|_{\underline{B}_{2,1}^{\f d 2}}+\|\vae(\tau)\|_{B_{2,1}^{\f d 2}}\Big)\|(\vVe-\vV,\mathbb{P}\vue-\vv)(\tau)\|^{l;\zeta}_{B_{2,1}^{\f d 2}}d\tau.\nonumber
\end{align}

According to \eqref{estpri25}-\eqref{estpri26}, we obtain that
\begin{align}\label{estpri51}
R_2(t)&\lesssim\int_0^t \ep\|\vae(\tau)\|_{B_{2,1}^{\f d 2}}\|\vue(\tau)\|_{\underline{B}_{2,1}^{\f d 2+1}}d\tau\\
 &\lesssim \int_0^t \|\vue(\tau)\|^{h;\zeta}_{B_{2,1}^{\f d 2+1}}[\vae,\vue]^{\zeta,\ep}_{h,m}(\tau)d\tau+(\ep\zeta)\int_0^t \|\vae(\tau)\|_{B_{2,1}^{\f d 2}}\|\vue(\tau)\|^{l;\zeta}_{\underline{B}_{2,1}^{\f d 2}}d\tau\nonumber
\end{align}
and
\begin{align*}
 R_3&\lesssim\int_0^t \|\vue(\tau)\|_{\underline{B}_{2,1}^{\f d 2}}\Big(\|\vae(\tau)\|^{h;\f {\zeta} {4}}_{B_{2,1}^{\f d 2}}+\|\vue(\tau)\|^{h;\f {\zeta} {4}}_{B_{2,1}^{\f d 2}}\Big)d\tau\\
 &\lesssim\int_0^t \|\vue(\tau)\|_{\underline{B}_{2,1}^{\f d 2}}\Big(\|\vae(\tau)\|^{h;\f {\zeta} {16}}_{B_{2,1}^{\f d 2}}+\|\vue(\tau)\|^{h;\f {\zeta} {16}}_{B_{2,1}^{\f d 2}}\Big)d\tau.
\end{align*}

From \eqref{estpri44}-\eqref{estpri50}, we find that
\begin{align}\label{estpri52}
R_3\lesssim &\int_0^t \Big(\|\vae(\tau)\|^2_{B_{2,1}^{\f d 2}}+\|\vue(\tau)\|^2_{\underline{B}_{2,1}^{\f d 2}}\Big)[\vae,\vue]^{\zeta,\ep}_{h,m}(\tau)d\tau+\int_0^t\|\vue(\tau)\|_{\underline{B}_{2,1}^{\f d 2}}\|(\vV,\vv)(\tau)\|^{h;\f \zeta {16}}_{B_{2,1}^{\f d 2}}d\tau\\
&+\int_0^t\|\vue(\tau)\|_{\underline{B}_{2,1}^{\f d 2}}\|(\vVe-\vV,\mathbb{P}\vue-\vv)(\tau)\|^{l;\zeta}_{\underline{B}_{2,1}^{\f d 2}}d\tau+(\ep\zeta)\int_0^t \f1 \ep\|\vae(\tau)\|^{h;\f {\eta_0} {\ep}}_{B_{2,1}^{\f d 2}}\|\vue(\tau)\|_{B_{2,1}^{\f d 2-1}}d\tau\nonumber.
\end{align}

Note that
\begin{align}\label{estpri55}
\ep\|\vae\|_{\widetilde{L}_t^\infty(B^{\f d 2}_{2,1})}\lesssim  \ep\|\vae\|^{h;\f{\eta_0}\ep}_{\widetilde{L}_t^\infty(B^{\f d 2}_{2,1})}+\|\vae\|^{m;\zeta,\f{\eta_0}\ep}_{\widetilde{L}_t^\infty(B^{\f d 2-1}_{2,1})}+(\zeta\ep)\|\vae\|^{l;\zeta}_{\widetilde{L}_t^\infty(B^{\f d 2-1}_{2,1})}.
\end{align}
Substituting \eqref{estpri39} and \eqref{estpri42}-\eqref{estpri52}  into \eqref{estpri38} and adding the resulting expression to \eqref{estpri37} and \eqref{estpri55}, we get that
\begin{align}\label{estpri53}
   [\vae,\vue]^{\zeta,\ep}_{h,m}(t)
  &\lesssim \left(\ep\|\vae\|^{h;\f{\eta_0}\ep}_{B_{2,1}^{\f d 2}}(0)+
  \|(\vae,\vue)\|^{h;\zeta}_{B_{2,1}^{\f d 2-1}}(0)\right)+R_{\zeta,\ep}(t)\nonumber\\
  &\quad+\int_0^t\Big(\|\vae(\tau)\|_{B_{2,1}^{\f d 2}}+\|\vue(\tau)\|_{\underline{B}_{2,1}^{\f d 2}}\Big)\|(\vVe-\vV,\mathbb{P}\vue-\vv)(\tau)\|^{l;\zeta}_{\underline{B}_{2,1}^{\f d 2}}d\tau\\
 &\quad+\int_0^t \Big(\|\vae(\tau)\|^2_{B_{2,1}^{\f d 2}}+\|\vue(\tau)\|^2_{\underline{B}_{2,1}^{\f d 2}}+\|\vue(\tau)\|^{h;\zeta}_{B_{2,1}^{\f d 2+1}}+\f 1 {\ep}\|\vae(\tau)\|^{h;\f {\eta_0} {\ep}}_{B_{2,1}^{\f d 2}}\Big)[\vae,\vue]^{\zeta,\ep}_{h,m}(\tau)d\tau,\nonumber
\end{align}
where
\begin{align*}
 R_{\zeta,\ep}(t)&=\int_0^t \Big(\|\vue(\tau)\|_{\underline{B}_{2,1}^{\f d 2}}+\|\vae(\tau)\|_{B_{2,1}^{\f d 2}}\Big)\|(\vv,\vV)(\tau)\|^{h;\f \zeta {16}}_{B_{2,1}^{\f d 2}}d\tau+\|f\|^{h;\zeta}_{L_t^1(B_{2,1}^{\f d 2-1})}\\
 &+\max\{\f 1 \zeta,\ep,\zeta\ep\} \int_0^t \|\vue(\tau)\|_{B_{2,1}^{\f d 2-1}}
\Big(\f1 \ep\|\vae(\tau)\|^{h;\f{\eta_0}\ep }_{B_{2,1}^{\f d 2}}+\|\vae(\tau)\|^{m;\zeta,\f{\eta_0}\ep }_{B_{2,1}^{\f d 2+1}}+\|\vue(\tau)\|^{h;\zeta}_{B_{2,1}^{\f d 2+1}}\Big)d\tau\\
&+(\ep\zeta)\int_0^t \|\vae(\tau)\|_{B_{2,1}^{\f d 2}}\|\vue(\tau)\|_{\underline{B}_{2,1}^{\f d 2}}d\tau+(\zeta\ep)\|\vae\|^{l;\zeta}_{\widetilde{L}_t^\infty(B^{\f d 2-1}_{2,1})}
\end{align*}

By Gronwall's inequality, we have that
\begin{align*}
 &[\vae,\vue]^{\zeta,\ep}_{h,m}(t)\\
 &\; \leq Ce^{\int_0^t \|\vae(\tau)\|^2_{B_{2,1}^{\f d 2}}+\|\vue(\tau)\|^2_{\underline{B}_{2,1}^{\f d 2}}+\|\vue(\tau)\|^{h;\zeta}_{B_{2,1}^{\f d 2+1}}+\f 1 {\ep}\|\vae(\tau)\|^{h;\f {\eta_0} {\ep}}_{B_{2,1}^{\f d 2}}d\tau}\Bigg(\ep\|\vae\|^{h;\f{\eta_0}\ep}_{B_{2,1}^{\f d 2}}(0)+
  \|(\vae,\vue)\|^{h;\zeta}_{B_{2,1}^{\f d 2-1}}(0)\\
 &\quad+R_{\zeta,\ep}(t)+\int_0^t\Big(\|\vae(\tau)\|_{B_{2,1}^{\f d 2}}+\|\vue(\tau)\|_{\underline{B}_{2,1}^{\f d 2}}\Big)\|(\vVe-\vV,\mathbb{P}\vue-\vv)(\tau)\|^{l;\zeta}_{\underline{B}_{2,1}^{\f d 2}}d\tau\Bigg)\\
& \;\leq Ce^{C( Y^{\zeta,\ep}(t)+ (Y^{\zeta,\ep}(t))^2)}\Bigg(\ep\|\vae\|^{h;\f{\eta_0}\ep}_{B_{2,1}^{\f d 2}}(0)+
  \|(\vae,\vue)\|^{h;\zeta}_{B_{2,1}^{\f d 2-1}}(0)+Y^{\zeta,\ep}(t)\|(\vv,\vV)\|^{h;\f\zeta{16}}_{\widetilde{L}_t^2(B^{\f d 2}_{2,1})}\nonumber\\
  &\quad+\max\{\f 1 \zeta,\ep,\zeta\ep\}Y^{\zeta,\ep}(t)[\vae,\vue]^{\zeta,\ep}_{h,m}(t)+(\zeta\ep)(Y^{\zeta,\ep}(t)+(Y^{\zeta,\ep}(t))^2)\nonumber\\
  &\quad+ [\vVe-\vV,\mathbb{P}\vue-\vv]^{\zeta,\ep}_{l}(t)Y^{\zeta,\ep}(t)+\|f\|^{h;\zeta}_{L_t^1(B_{2,1}^{\f d 2-1})}\Bigg).\ \  \Box
\end{align*}

\section{Appendix}\label{Sec5}

\begin{Lemma}\emph{(\cite{D1}).}\label{Be2} Let $G$ which vanishes at $0$ be a smooth function on an open interval $I\subset\mathbb{R}$ satisfying $(-R,R)\subset I$ with some $R>0$. If $s>0$,  then for any $f\in\dot{B}_{{2},1}^{s}(\mathbb{T}^d)$   satisfying $\|f\|_{L^\infty}<R$, there exists a positive constant $C=C(d,s,R,G)$ such that
\[
\|G(f)\|_{\dot{B}_{2,1}^s}\leq C
\|f\|_{\dot{B}_{2,1}^s}.
\]
\end{Lemma}
\begin{Lemma}\emph{(\cite{D1})}\label{leD2}. Let A be a function on $\widetilde{\mathbb{Z}}^d\times \widetilde{\mathbb{Z}}^d$ such that for some $\alpha\in \mathbb{R}, \beta, \gamma\geq0 $, we have
\[
|A(k,m)|\leq K |m|^{\alpha} |k|^{\beta} |m-k|^{\gamma}.
\]
Let
\[
\mathcal{Q}(u,v)=A(k,m)\hat{u}_k \hat{v}_{m-k}e^{im\cdot x}.
\]
Then the following inequality holds true :
\[
\|\mathcal{Q}(u,v)\|_{H^{\sigma-\alpha}}\lesssim K\left(\|\mathcal{F}u\|_{\textit{l}^1(\widetilde{\mathbb{Z}}^d)}\|v\|_{H^{\sigma+\beta+\gamma}}
+\|\mathcal{F}v\|_{\textit{l}^1(\widetilde{\mathbb{Z}}^d)}\|u\|_{H^{\sigma+\beta+\gamma}}\right) \text{ if } \sigma\geq0,
\]
\[
\|\mathcal{Q}(u,v)\|_{H^{\sigma-\alpha}}\lesssim K\left(\|\mathcal{F}u\|_{\textit{l}^1(\widetilde{\mathbb{Z}}^d)}\|\underline{v}\|_{H^{\sigma+\beta+\gamma}}
+\|\mathcal{F}v\|_{\textit{l}^1(\widetilde{\mathbb{Z}}^d)}\|\underline{u}\|_{H^{\sigma+\beta+\gamma}}\right) \text{ if } \sigma\geq0 \text{ and } \sigma+\beta+\gamma>0.
\]
\end{Lemma}
\begin{Lemma}\emph{(\cite{D1})}\label{leD1}. Let $A,B\in (\text{Ker}L)^\perp$, $0<s<1+\f d 2$, and $0<\theta<1$. Then the following estimates hold:
\begin{align}
\|  \widecheck{\mathcal{Q}}_1^\ep(\vu,B)\|_{H^{s-2}}&\lesssim \|\vu\|_{L^\infty\cap H^{\f d 2}}\|B\|_{H^{s-1}}\label{estQ1}\\
 \|  \widetilde{\mathcal{Q}}_1^\ep(\vu,B)\|_{H^{\f d 2-1-\theta}}&\lesssim \|\vu\|_{ \underline{H}^{\f d 2}}\|B\|_{H^{\f d 2-\theta}}\label{estQ3},\\
  \|  \mathcal{Q}_2^\ep(A,B)\|_{H^{s-2}}&\lesssim \min\{\|A\|_{B^{\f d 2}_{2,1}}\|B\|_{H^{s-1}}, \|B\|_{B^{\f d 2}_{2,1}}\|A\|_{H^{s-1}}\}\label{estQ2}.
\end{align}
\end{Lemma}
\begin{Lemma}\emph{(\cite{D1})}\label{leD4}. Let $A,B\in (\text{Ker}L)^\perp$. Then following estimates are valid:
\begin{align}
  \|\mathcal{Q}_1(\vu,A)\|_{H^{\sigma-1}}&\lesssim \|A\|_{H^{\sigma+\gamma}}\||k|^{-\gamma}\mathcal{F}\underline{\vu}\|_{\textit{l}^1(\widetilde{\mathbb{Z}}^d)}+|\widehat{\vv}_0|\|A\|_{H^{\sigma}}\ \ \text{ if }\ \  \sigma, \gamma\geq0,\label{estQ4}\\
  \|\mathcal{Q}_2(A,B)\|_{H^{\sigma-2}}&\lesssim \min\{\|A\|_{B^{\f 1 2}_{2,1}}\|B\|_{B^{\sigma-1}_{2,1}}, \|B\|_{B^{\f 1 2}_{2,1}}\|A\|_{B^{\sigma-1}_{2,1}}\}, \ \ \text{ if }\ \  \f1 2<\sigma<\f3 2,\label{estQ5}\\
    \|\mathcal{Q}_2(A,B)\|_{H^{\sigma-2}}&\lesssim \|A\|_{B^{\sigma-1}_{2,1}}\|B\|_{B^{\sigma-1}_{2,1}}, \ \ \text{ if }\ \  \f3 2\leq\sigma,\label{estQ6}\\
    \|\mathcal{Q}_2(A,A)\|_{B^{\sigma-1}}&\lesssim \|A\|_{B^{\f 1 2}_{2,1}} \|A\|_{B^{\sigma}_{2,1}}, \ \ \text{ if }\ \  -\f1 2<\sigma.\label{estQ7}
\end{align}
\end{Lemma}
\centerline{\bf Acknowledgements}

\vspace{2mm}
I would like to extend my deepest gratitude to Professor Yongzhong Sun for his invaluable encouragement and mentorship throughout this research. I am also sincerely thankful to Professor Yang Li for his intellectual guidance and especially for recommending the seminal works \cite{F, F1}. Lastly, I wish to thank the referees for their insightful suggestions regarding the references, notation, and overall structure of the paper.
\vspace{5mm}

\centerline{\bf Conflict of interest}
\vspace{2mm}
The  author states that there is no conflict of interest.

\vspace{5mm}

\centerline{\bf Data Availability Statement}
\vspace{2mm}
Data sharing not applicable to this article as no datasets were generated or
analyzed during the current study.



\begin{thebibliography}{100}
\bibitem{BCD}   Bahouri, H., Chemin, J.-Y., Danchin, R.: Fourier analysis and nonlinear partial differential equations. Grundlehren der mathematischen Wissenschaften [Fundamental Principles of Mathematical Sciences], 343. Springer, Heidelberg, 2011.
\bibitem{CL} Chemin, J.-Y., Lerner, N.: Flow of non-Lipschitz vector fields and Navier-Stokes equations. {J. Differential Equations.} {\bf 121}, 314-328 (1995)
\bibitem{DG}  Desjardins, B., Grenier, E.: Low Mach number limit of viscous compressible flows in the whole space. {R. Soc. Lond. Proc. Ser. A Math. Phys. Eng. Sci.} {\bf 455}, 2271-2279 (1999)
\bibitem{DGLM}  Desjardins, B., Grenier, E., Lions, P.-L., Masmoudi, N.: Incompressible limit for solutions of the isentropic Navier-Stokes equations with Dirichlet boundary conditions, {J. Math. Pures Appl.} {\bf 78}, 461-471 (1999)
\bibitem{D1} Danchin, R.: Zero Mach number limit for compressible flows with periodic boundary conditions. {Amer. J. Math.} {\bf 124}, 1153-1219 (2002)
\bibitem{D2} Danchin, R., He, L.: The incompressible limit in Lp type critical spaces. {Math. Ann.} {\bf 366}, 1365-1402 (2016)
\bibitem{D3} Danchin, R.: Zero Mach number limit in critical spaces for compressible Navier-Stokes equations. {Ann. Sci. $\acute{E}$cole Norm. Sup. (4)}, {\bf 35}, 27-75 (2002)
\bibitem{D4} Danchin, R.: Global existence in critical spaces for compressible Navier-Stokes equations. {Invent. Math.} {\bf 141}, 579-614 (2000)
\bibitem{D5} Danchin, R.: Well-posedness in critical spaces for barotropic viscous fluids with truly not constant density. {Comm. Partial Differential Equations.} {\bf 32}, 1373-1397 (2007)
 \bibitem{EDG} Ebin, D. G.:The motion of slightly compressible fluids viewed as a motion with strong constraining force, {Annals of Mathematics}, {\bf 105}, 141-200 (1977)
\bibitem{Fe3} Feireisl, E., Karper, T., Kreml, O., Stebel, Jan.: Stability with respect to domain of the low Mach number limit of compressibl viscous fluids. {Math. Models Methods Appl. Sci.} {\bf 23}, 2465-2493 (2013)
\bibitem{Fe4} Feireisl, E., Novotn\'{y}, A.:  Singular limits in thermodynamics of viscous fluids, Second, Advances in Mathematical Fluid Mechanics, Birkh\"{a}user/Springer, Cham, 2017.
\bibitem{F} Fujii, M.: Low Mach number limit of the global solution to the compressible Navier-Stokes system for large data in the critical Besov space. {Math. Ann.} {\bf 388}, 4083-4134 (2024)
\bibitem{F1} Fujii, M., Li, Y.: Low Mach number limit for the global large solutions to the 2D Navier--Stokes--Korteweg system in the critical $\widehat{L^p}$ framework. {Calc. Var. Partial Differential Equations.} {\bf 64}, Paper No. 29 (2025)
\bibitem{G1} Gallagher, I.: Applications of Schochet's methods to parabolic equations. {J. Math. Pures Appl.} {\bf 77}, 989-1054 (1998)
\bibitem{Gr} Grenier, E.: Oscillatory perturbations of the Navier Stokes equations, {J. Math. Pures Appl.}, {\bf 76}, 477-498 (1997).
\bibitem{DH} Hoff, D.: The Zero-Mach limit of compressible flows. {Commun. Math. Phys.}  {\bf 192}, 543-554 (1998)
\bibitem{HLo} Hagstrom, T., Lorenz, J.: All-time existence of classical solutions for slightly compressible flows. {SIAM J. Math. Anal.} {\bf 29}, 652-672 (1998)
\bibitem{JN}  Jiang, N., Masmoudi, N.: Low Mach number limits and acoustic waves, Handbook of mathematical analysis in mechanics of viscous fluids, 2018, pp. 2721-2770.
\bibitem{KM} Klainerman, S., Majda, A.: Singular limits of quasilinear hyperbolic systems with large parameters and the incompressible limit of compressible fluids, {Comm. Pure Appl. Math.} {\bf 35}, 629-651 (1982)
\bibitem{KM1} Klainerman, S., Majda, A.: Compressible and incompressible fluids, {Comm. Pure Appl. Math.}  {\bf 35} (1982), 629-651 (1982)
\bibitem{LS} Li, S.: Incompressible limit of global solutions to the compressible magnetohydrodynamics equations for large initial data in critical Besov spaces. {J. Math. Phys.} {\bf 66}, 071508 (2025)
\bibitem{LM}  Lions, P.-L., Masmoudi, N.: Incompressible limit for a viscous compressible fluid. {J. Math. Pures Appl.} {\bf 77}, 585-627 (1998)
\bibitem{Lif} Li, F. C.,  Mu, Y. M.: Low Mach number limit for the compressible magnetohydrodynamic equations in a periodic domain.  {Discrete Contin. Dyn. Syst.} {\bf 38}, 1669-1705 (2018)
\bibitem{M} Masmoudi, N.: Incompressible, inviscid limit of the compressible Navier-Stokes system. {Ann. Inst. H. Poincar$\acute{e}$ Anal. Non Lin$\acute{e}$aire.} {\bf 18}, 199-224 (2001)
\bibitem{O} Ogino, S. Strong Convergence of Low Mach Number Limit for the Compressible Navier-Stokes Equations in the Scaling Critical Spaces. {J. Math. Fluid Mech.} {\bf 27}, 40 (2025)
\bibitem{SS2} Schochet, S.: Fast singular limits of hyperbolic PDEs. {J. Differential Equations} {\bf 114}, 476-512 (1994)
\bibitem{SS1} Schochet, S.: The mathematical theory of the incompressible limit in fluid dynamics, Handbook of mathe-matical fluid dynamics. Vol. IV, 2007, pp. 123-157.
\end{thebibliography}
\end{document}